\documentclass[12pt,a4paper,twoside,reqno]{amsart}

\title{On the geometry of moduli spaces of holomorphic chains over compact
Riemann surfaces}

\date{20 Dec 2005}

\author[L.\ \'Alvarez-C\'onsul]{L.\ \'Alvarez-C\'onsul}
\address{CSIC, Serrano 123, 28006 Madrid, Spain}
\email{lac@mat.csic.es}
\author[O.\ Garc\'{\i}a-Prada]{O.\ Garc\'{\i}a-Prada}
\address{CSIC, Serrano 121, 28006 Madrid, Spain}
\email{oscar.garcia-prada@uam.es}
\author[A.H.W.\ Schmitt]{A.H.W.\ Schmitt}
\address{ Universit\"at Duisburg-Essen, Campus Essen, FB6 Mathematik und Informatik,
D-45117 Essen, Germany}
\email{alexander.schmitt@uni-essen.de}


\usepackage{amsthm, amssymb, amscd} 
\theoremstyle{plain}
\newtheorem{theorem}{Theorem}[section]
\newtheorem{lemma}[theorem]{Lemma}
\newtheorem{corollary}[theorem]{Corollary}
\newtheorem{proposition}[theorem]{Proposition}
\newtheorem{conjecture}[theorem]{Conjecture}
\theoremstyle{definition}
\newtheorem{definition}[theorem]{Definition}
\newtheorem{definition-theorem}[theorem]{Definition-Theorem}
\newtheorem{example}[theorem]{Example}

\theoremstyle{remark}
\newtheorem{remark}[theorem]{Remark}

\newcommand{\secref}[1]{\S\ref{#1}}

\numberwithin{equation}{section}
\newcommand{\Hom}{\operatorname{Hom}}
\newcommand{\End}{\operatorname{End}}

\newcommand{\coker}{\operatorname{coker}}
\newcommand{\surj}{\to\kern-1.8ex\to}
\newcommand{\hra}{\hookrightarrow}
\newcommand{\wt}{\widetilde}
\newcommand{\id}{\operatorname{id}}
\newcommand{\im}{\operatorname{im}}
\newcommand{\Ext}{\operatorname{Ext}}
\newcommand{\CC}{{\mathbb C}}
\newcommand{\HH}{\mathbb{H}}
\newcommand{\lto}{\longrightarrow}

\usepackage[all]{xy}\CompileMatrices\SelectTips{cm}{12}
\DeclareMathOperator{\tr}{tr}

\newcommand{\C}{\mathbb{C}}
\newcommand{\cD}{\mathcal{D}}
\newcommand{\dbar}{\bar{\partial}}
\newcommand{\lsurj}{\longrightarrow\kern-1.8ex\to}
\newcommand{\cB}{\mathcal{B}}
\newcommand{\cT}{\mathcal{T}}
\newcommand{\cV}{\mathcal{V}}

\usepackage{pstricks, floatflt, pst-text}
\newtheorem*{claim}{Claim}

\newrgbcolor{rosa}{1 .8 .8}
\newrgbcolor{dunkelrosa}{.75 .55 .55}
\newrgbcolor{orange}{1 .5 .25}
\newrgbcolor{lightorange}{1 .75 .5}
\newrgbcolor{verylightorange}{1 .85 .75}
\newrgbcolor{lightgrey}{.65 .65 .65}
\newrgbcolor{vlightgrey}{.8 .8 .8}
\newrgbcolor{lightgreen}{.8 1 .8}
\newrgbcolor{darkgreen}{.25 .75 .25}
\newrgbcolor{mix}{1 1 .8}
\newrgbcolor{lila}{.75 .5 .75}
\newrgbcolor{bblue}{.1 .25 .5}
\newrgbcolor{lblue}{.8 .9 1}

\newcommand{\adeg}{\mathop{{\rm deg}_\alpha}}
\newcommand{\GL}{\mathop{\rm GL}}
\newcommand{\U}{\mathop{\rm U}}
\newcommand{\gr}{\mathop{\rm gr}}
\newcommand{\Id}{\mathop{\rm id}}
\newcommand{\la}{\lambda}
\newcommand{\lma}{\longmapsto}
\newcommand{\lra}{\longrightarrow}
\newcommand{\OO}{{\mathcal O}}
\newcommand{\ol}{\overline}
\newcommand{\p}{\prime}
\newcommand{\PP}{{\mathbb P}}
\newcommand{\q}{\quad}
\newcommand{\Q}{{\mathbb Q}}
\newcommand{\R}{{\mathbb R}}
\newcommand{\rk}{\mathop{\rm rk}}
\newcommand{\SL}{\mathop{\rm SL}}

\newcommand{\vol}{\mathop{\rm vol}}
\newcommand{\Z}{{\mathbb Z}}

\begin{document}

\begin{abstract}
We study holomorphic $(n+1)$-chains $E_n\rightarrow E_{n-1}
\rightarrow \cdots \rightarrow E_0$ consisting of holomorphic vector
bundles over a compact Riemann surface and homomorphisms between
them. A notion of stability depending on $n$ real parameters was
introduced in \cite{AG1} and moduli spaces were constructed in
\cite{Schtree,Schquiv}. In this paper we study the variation of the
moduli spaces with respect to the stability parameters. In particular
we characterize a parameter region where the moduli spaces are
birationally equivalent. A detailed study is given for the case of
$3$-chains, generalizing that of $2$-chains (triples) in \cite{BGG}.
Our work is motivated by the study of the topology of moduli spaces
of Higgs bundles and their relation to representations of the
fundamental group of the surface.
\end{abstract}

\maketitle

\pagestyle{myheadings}
\markboth{\sc L.\ \'Alvarez-C\'onsul, O.\ Garc\'{\i}a-Prada,
A.H.W.\ Schmitt}{\sc Holomorphic Chains}
\tableofcontents

\section{Introduction}                          
\label{sec:introduction}

Let $X$ be a compact Riemann surface of genus $g\geq 2$. 
A holomorphic $(n+1)$-chain over $X$ is an object
$$
E_n\stackrel{\phi_n}{\lra} E_{n-1} \stackrel{\phi_{n-1}}{\lra} \cdots 
\stackrel{\phi_1}{\lra}  E_0
$$
consisting of holomorphic vector bundles $E_j$ on $X$,
$j=0,...,n$, and homomorphisms $\phi_i\colon E_{i}\lra E_{i-1}$,
$i=1,...,n$. The ranks and degrees of $E_i$ define the type of the
chain. A notion of stability for $(n+1)$-chains, depending on $n$
real parameters $\alpha_i$, has been introduced in \cite{AG1} and
moduli spaces have been constructed in \cite{Schtree,Schquiv}. These
objects generalize the holomorphic triples $ E_1\stackrel{\phi}{\lra}
E_0$ introduced in \cite{G,BG}.  The variation of the moduli spaces
of holomorphic triples with respect to the stability parameter
$\alpha$ has been studied in \cite{BGG}, where a birational
characterization of the moduli space has been given. It turns out
that the moduli space of $\alpha$-stable triples is empty outside of
an interval $(\alpha_m,\alpha_M)$ where the bounds are determined by
the type of the triple, and $\alpha_M=\infty$ if the ranks of $E_0$
and $E_1$ are equal.  The main result in \cite{BGG} is that for
$\alpha\in (\alpha_m,\alpha_M)$ and $\alpha\geq 2g-2$, the moduli
space of $\alpha$-stable triples is non-empty, smooth and
irreducible.

In this paper we undertake a systematic study of holomorphic
$(n+1)$-chains for arbitrary $n$. We study the parameter region where
the moduli spaces may be non-empty. This region is partitioned into
chambers and we study the variations in the moduli space as we cross
a wall.  We show that the region is bounded by $n$ hyperplanes, which
play the role of $\alpha_m$ in the case of triples. The determination
of other bounding hyperplanes --- the analogues of $\alpha_M$ --- is
more difficult and is only done in some cases. However, we
characterize a region where the moduli spaces of chains of a given
type are birationally equivalent. It turns out that, similarly to the
case of triples, the stability parameters $\alpha_i$ must satisfy
$\alpha_i-\alpha_{i-1}\geq 2g-2$. After developing the general
theory, we study in more detail the case of $3$-chains and finish the
paper giving the birational characterization of the moduli spaces for
some special values of the ranks.

Our main motivation to study this problem comes from the theory of
Higgs bundles on a Riemann surface $X$ and its relation to
representations of the fundamental group of $X$.  By results of
Hitchin \cite{H}, Donaldson \cite{D}, Simpson \cite{S}, and Corlette
\cite{C}, the moduli space of reductive representations of the
fundamental group of $X$ in a non-compact reductive Lie group $G$ can
be identified with the moduli space of polystable $G$-Higgs bundles.
As shown by Hitchin, the $L^2$-norm of the Higgs field with respect
to the solution of the Hitchin equations defines a proper function on
the moduli space of $G$-Higgs bundles, which in many cases is a
perfect Morse-Bott function. Hence the study of the topology of the
moduli space of Higgs bundles, such as Betti numbers, reduces to the
study of the topology of the critical subvarieties of the Morse
function.  It turns out \cite{S,AG1} that for $G=\GL(n,\C)$ and
$G=\U(p,q)$ these critical subvarieties correspond precisely to
moduli spaces of $(n+1)$-chains for different values of $n$ and for
certain values of the stability parameters, namely
$\alpha_i-\alpha_{i-1}= 2g-2$ --- exactly the extremes of the region
where our results on the birationality of moduli spaces of chains
apply.

It is indeed the irreducibility of the moduli space of triples proved
in \cite{BGG} that has allowed to count in \cite{BGG2} the number of
connected components of the moduli space of representations of the
fundamental group (and its universal central extension) in $\U(p,q)$.
Also, the computation of the Betti numbers of the moduli spaces of
triples $\mathcal{O} \stackrel{\phi}{\lra} E$ with $\rk (E)=2$ by
Thaddeus \cite{Th} has enabled Gothen \cite{Go1,Go2} to compute the
Betti numbers of the moduli spaces of $\SL(3,\CC)$-Higgs bundles and
$\U(2,1)$-Higgs bundles, and the parabolic versions given in
\cite{GGM} and \cite{L}.

To carry out this programme, one has then to study the topology of
the moduli spaces of $(n+1)$-chains for arbitrary $n$ when
$\alpha_i-\alpha_{i-1}= 2g-2$. As it turns out in most of the cases
studied so far, it is easier to understand the moduli space for some
particular chamber, and then analyze the wall-crossings until we get
to $\alpha_i-\alpha_{i-1}= 2g-2$.  In this paper we give the first
steps in this direction beyond the $n=1$ case.  In fact, the examples
we consider include the case $ E_2\stackrel{\phi_2}{\lra} E_{1}
\stackrel{\phi_{1}}{\lra} E_0$ with $\rk(E_1)=\rk(E_0)=1$, whose
Betti numbers, together with the Betti numbers of the moduli space of
triples $E_{1} \stackrel{\phi}{\lra} E_0$ with $\rk (E_0)=1$ would
give the Betti numbers of the moduli space of representations of the
fundamental group in $\U(n,1)$ --- the group of isometries of the
$n$-dimensional hyperbolic space. This is indeed a computation that
we plan to undertake in a future paper.

\subsection*{Acknowledgments}   

This work has been partially supported by the Spanish Ministerio de
Educaci\'on y Ciencia (MEC) and the German DAAD via the ``Acciones
Integradas Hispano-Alemanas'' programme, Contract No.\ HA2004-0083
(Spain)\allowbreak /\allowbreak D/04/\allowbreak 4\allowbreak
2\allowbreak 2\allowbreak 5\allowbreak 7 (Germany).
LAC\ was partially supported by the Spanish ``Programa Ram\'on y
Cajal''. LAC\ and OGP\ were partially supported by MEC under Grant
MTM\allowbreak 2004-07090-C03-01.
AHWS\ acknowledges support by the DFG via a Heisenberg fellowship and
via the priority programme ``Globale Methoden in der Komplexen
Geometrie---Global Methods in Complex Geometry''.

Parts of this paper have been written while AHWS\ and OGP\ stayed at
the Institut des Hautes \'Etudes Scientifiques, whose hospitality and
support is thanked.
During these visits, AHWS\ benefitted from support of the European
Commission through its 6th Framework Programme ``Structuring the
European Research Area'' and the Contract No.\ RITA-CT-2004-505493
for the provision of Transnational Access implemented as Specific
Support Action.

\section{Definitions and basic facts}           
\label{sec:definitions}                
In this section, we will, on the one hand, review the formalism of
holomorphic chains from \cite{AG1} as well as existing results and,
on the other hand, prove some substantial new results concerning the
stability parameters for holomorphic chains. Recall that we do work
exclusively on a smooth projective curve $X$ of genus $g\geq 2$
defined over the complex numbers.
\par
The concepts that we shall explain below are natural generalizations
of notions and ideas from the theory of vector bundles on the curve
$X$. Let us thus pause to recall the relevant concepts from the
setting of vector bundles.  If $E$ is a vector bundle on $X$, we
write $\rk(E)$ for its \it rank\rm, i.e., the dimension of the fiber
vector spaces and denote by $\deg(E)$ its \it degree\rm. Using Chern
classes, we find
$$
c_1(E)=\deg(E)\cdot [{\rm pt}].
$$
The \it slope of $E$ \rm is the quotient $\mu(E):=\deg(E)/\rk(E)$.
A vector bundle is said to be \it (semi)stable\rm, if for every
non-zero proper subbundle $0\subseteq F\subseteq E$ the inequality
$$
\mu(F)(\le)\mu(E)
$$
is satisfied. The symbol ``$(\le)$" means that ``$<$" is used in
the definition of ``stable" and ``$\le$" in the definition of
``semistable". The semistable vector bundles of rank $r$ and degree
$d$ are classified by an irreducible normal projective variety
${\mathcal U}(r,d)$ of dimension $r^2(g-1)+1$. It contains the
smooth dense open subvariety ${\mathcal U}^{\rm s}(r,d)$ that
parameterizes stable vector bundles.
\subsection{Definitions}
A \it holomorphic $(n+1)$-chain \rm is a tuple
$C=(E_j,j=0,...,n;\phi_i, i=1,...,n)$, consisting of vector
bundles$E_j$ on $X$, $j=0,...,n$, and homomorphisms $\phi_i\colon
E_{i}\lra E_{i-1}$, $i=1,...,n$. The tuple
$t:=(\rk(E_j),j=0,...,n;\deg(E_j),j=0,...,n)$ will be referred to as
the \it type of the chain $(E_j,j=0,...,n;\phi_i, i=1,...,n)$\rm.  We
will often write a chain in the form
$$
C: E_n\stackrel{\phi_n}{\lra} E_{n-1} \stackrel{\phi_{n-1}}{\lra}
\cdots \stackrel{\phi_1}{\lra} E_0.
$$
A \it subchain of the holomorphic chain $(E_j,j=0,...,n;\phi_i,
i=1,...,n)$ \rm is a tuple $C':=(F_j, j=0,...,n)$ with $F_j$ a
subsheaf of $E_j$, $j=0,...,n$, such that $\phi_i(F_i)\subseteq
F_{i-1}$, $i=1,...,n$. The subchains $(0,...,0)$ and
$(E_j,j=0,...,n)$ are called the \it trivial subchains\rm.
\begin{remark}
Note that a subchain $(F_j, j=0,...,n)$ gives rise to the holomorphic
chain $(F_j, j=0,...,n;\phi_{i|F_{i}}\colon\allowbreak F_i\lra
F_{i-1}, i=1,...,n)$.
\end{remark}
Now, fix a tuple $\alpha=(\alpha_j,j=0,...,n)$ of real numbers.  For
a holomorphic $(n+1)$-chain $C=(E_j,j=0,...,n;\phi_i, i=1,...,n)$, we
define the \it $\alpha$-degree \rm as
$$
\adeg(C):=\sum_{j=0}^{n}\bigl(\deg(E_j)+\alpha_j \rk(E_j)\bigr)
$$
and the \it $\alpha$-slope \rm as
$$
\mu_{\alpha}(C):= \frac{\adeg(C)} {\sum_{j=0}^{n}\rk(E_j)}.
$$
A holomorphic $(n+1)$-chain $C$ is said to be \it
$\alpha$-(semi)sta\-ble\rm, if the inequality
$$
\mu_{\alpha}\bigl(F_j,j=0,...,n;\phi_{i|F_{i}}, i=1,...,n\bigr)
(\le) \mu_{\alpha}(C)
$$
is verified for any non-trivial subchain $C'=(F_j, j=0,...,n)$ of
$C$. Here, the convention for ``$(\le)$" is as before.  Last but not
least, we call a chain $C$ \it $\alpha$-polystable\rm, if it may be
written as a direct sum $C=C_1\oplus\cdots\oplus C_t$ where $C_k$ is
an $\alpha$-stable holomorphic chain with
$\mu_\alpha(C_k)=\mu_\alpha(C)$, $k=1,...,t$. Since holomorphic
chains form in a natural way an Abelian category, one easily derives
the following result.
\begin{proposition}[The Jordan-H\"older filtration]
\label{eq:JHF}
Let $C$ be an $\alpha$-semi\-stable holomorphic chain. Then, $C$
possesses a (in general non-unique) so-called Jordan-H\"older
filtration
$$
0=:C_0\subsetneq C_1\subsetneq  \cdots \subsetneq C_m:=C
$$
by holomorphic subchains, such that
$\mu_\alpha(C_i)=\mu_\alpha(C)$ and $C_i/C_{i-1}$ is $\alpha$-stable,
$i=1,...,m$.  The so-called graduation
$$
G:={\rm gr}(C):=\bigoplus_{i=1}^m G_i,\q  G_i:=C_i/C_{i-1},\ i=1,...,m,
$$
of $C$ is then $\alpha$-polystable. The equivalence class of ${\rm
gr}(C)$ does not depend on the Jordan-H\"older filtration of $C$.
\end{proposition}
Using the above proposition, we call two $\alpha$-semistable
holomorphic chains $C$ and $C'$ \it S-equivalent\rm, if their
graduations ${\rm gr}(C)$ and ${\rm gr}(C')$ are equivalent.
\begin{remark}
\label{quotientformulation}
i) Suppose $C'=(F_j, j=0,...,n)$ is a subchain of the holomorphic
chain $C=(E_j,j=0,...,n;\phi_i, i=1,...,n)$.  Let $\ol{F}_j$ be the
subbundle of $E_j$ generated by $F_j$, $j=0,...,n$.  Then,
$\ol{C}':=(\ol{F}_j, j=0,...,n)$ is still a subchain with
$\deg(F_j)\le \deg(\ol{F}_j)$, $j=0,...,n$. Thus, semistability has
to be checked only against subchains composed of subbundles.
\par
ii) Let $C=(E_j,j=0,...,n;\phi_i, i=1,...,n)$ be a holomorphic chain.
A holomorphic chain $(Q_j,j=0,...,n;\psi_i, i=1,...,n)$ is called a
\it quotient chain of $C$\rm, if there exist surjective homomorphisms
$\pi_j\colon E_j\lra Q_j$, $j=0,...,n$, such that $\psi_i\circ
\pi_i=\pi_{i-1}\circ \phi_i$, $i=1,...,n$. Note that
$(\ker(\pi_j),j=0,...,n)$ is then a holomorphic subchain of $C$ and
that we have the \it trivial quotients \rm $C$ and $(0,...,0;
0,...,0)$. Moreover, for any subchain $(F_j,j=0,...,n)$, we obtain
the induced quotient chain $(E_j/F_j,j=0,...,n;\ol{\phi}_i,
i=1,...,n)$.
\par
Standard arguments now show that a holomorphic chain $C$ is
$\alpha$-(semi)stable, if and only if the inequality
$$
\mu_{\alpha}(C) (\le) \mu_{\alpha}\bigl(Q_j,j=0,...,n;\psi_i,
i=1,...,n\bigr)
$$
holds for any non-trivial quotient $(Q_j,j=0,...,n;\psi_i,
i=1,...,n)$ of $C$.
\par
iii) Let $\alpha=(\alpha_j,j=0,...,n)$ be as above and $\beta\in
\R$.Set $\alpha^\p:=(\alpha_j+\beta, j=0,...,n)$. Then, it is obvious
that a holomorphic chain $C$ is $\alpha$-(semi)stable, if and only if
it is $\alpha^\p$-(semi)stable. Thus, we may assume that $\alpha_0$
be zero. In particular, the semistability concept for holomorphic
$(n+1)$-chains depends only on $n$ rational parameters.
\par
iv) If $C=(E_j,j=0,...,n;\phi_i, i=1,...,n)$ is a holomorphic chain,
we get the \it dual holomorphic chain \rm
$C^\vee:=(E^\p_j,j=0,...,n;\phi^\p_i, i=1,...,n)$ with
$E^\p_j:=E^\vee_{n-j}$, $j=0,...,n$, and
$\phi_i^\p:=\phi^\vee_{n+1-i}$, $i=1,...,n$. Then, $C$ is
$(\alpha_0,...,\alpha_{n+1})$-(semi)stable, if and only if $C^\vee$
is $(-\alpha_{n+1},...,\allowbreak-\alpha_{1})$-(semi)stable.
\end{remark}
\subsection{Moduli spaces}
Given a fixed type $t=(r_j,j=0,...,n; d_j, j=0,...,n)$ and a fixed
stability parameter $\alpha=(\alpha_j, j=0,...,n)$ consisting of \sl
rational \rm numbers, the S-equivalence classes of
$\alpha$-semistable holomorphic chains $C$ of type $t$ form a
projective moduli scheme ${\mathcal M}_\alpha(t)$. The GIT
construction for ${\mathcal M}_{\alpha}(t)$ is contained in
\cite{Schquiv} (see also \cite{Schtree} where the semistability
condition appears in a different guise).  The chamber structure of
the parameter region (see Remark \ref{rationality1}) will reveal
that, for any parameter $\alpha$, there is a rational parameter
$\alpha^\p$, such that the notion of $\alpha$-(semi)stability for
holomorphic chains of type $t$ is equivalent to the notion of
$\alpha^\p$-(semi)stability, so that the moduli spaces do indeed
exist for any stability parameter.
\subsection{The parameter region}
\label{sub:parameter-region}
We need to study the moduli spaces ${\mathcal M}_\alpha(t)$ for
fixed type $t$ in dependence of the parameter $\alpha$. As remarked
in \ref{quotientformulation}, iii), we may write
$\alpha=(0,\alpha_1,...,\alpha_n)$. Thus, we may view the
stability parameter $\alpha$ as an element of $\R^n$. The first
interesting problem is to determine \sl a priori \rm a region
$R(t)\subset \R^n$, such that the existence of an $\alpha$-semistable
holomorphic chain of the predetermined type $t$ implies that
$\alpha$ lies in the region $R(t)$. Here, we show that the
semistability condition for certain natural subchains gives some
halfspaces $H_1$,...,$H_n$, such that $R(t)$ lies in the intersection
of those halfspaces. The next natural and important question is
whether $R(t)$ is bounded or not. We will relate this question to the
existence or non-existence of semistable chains in the category of
$k$-vector spaces where $k=\CC(X)$ is the function field of the
curve $X$.
\par
Let $C=(E_j, j=0,...,n;\phi_i,i=1,...,n)$ be a holomorphic
$(n+1)$-chain. Define the \it $i$-th standard subchain \rm to be
$C_i:=(E_0,...,E_i,0,...,0)$, $i=0,...,n-1$. A straightforward
calculation gives the following result.
\begin{proposition}
\label{StandardBounds} Suppose $C$ is $\alpha$-(semi)stable with
$\alpha=(0,\alpha_1,...,\alpha_n)$. Then, the condition arising from
the $i$-th standard subchain $C_i$ is
$$
\begin{array}{rl}
&\bigl(\alpha_1r_1+\cdots+\alpha_ir_i\bigr)\bigl(r_{i+1}+\cdots+r_n\bigr)-\hfill
\\
&- \bigl(\alpha_{i+1}r_{i+1}+\cdots+\alpha_nr_n\bigr)\bigl(r_0+\cdots+r_i\bigr)\hfill
\\
(\le) &\bigl(r_0+\cdots+r_i\bigr) \bigl(d_{i+1}+\cdots+d_n\bigr)-
\\
&- \bigl(r_{i+1}+\cdots+r_n\bigr)\bigl(d_0+\cdots+d_i\bigr),\q i=0,...,n-1.
\end{array}
$$
\end{proposition}
Let $h_i$ be the hyperplane determining the halfspace $H_i$ from
Proposition \ref{StandardBounds}, i.e., $h_i$ is defined by the
equation
\begin{equation}
\label{eq:standard-hyperplanes}
\begin{array}{rl}
&(r_{i+1}+...+ r_n)(r_1\alpha_1+...+r_i\alpha_i)-
\\
&-
(r_0+...+r_i)(r_{i+1}\alpha_{i+1}+...+r_n\alpha_n)\\
=&(r_0+...+r_i)(d_{i+1}+...+d_n)-
\\
&-(r_{i+1}+...+ r_n)(d_0+...+d_i),\q
i=0,...,n.
\end{array}
\end{equation}
\begin{example}
\label{ex:StandardBounds2}
In the case $n=2$, we find the following two inequalities:
\begin{eqnarray*}
\alpha_1r_0r_1+\alpha_2r_0r_2 &\ge& (r_1+r_2)d_0-r_0d_1-r_0d_2
\\
-\alpha_1r_1r_2+\alpha_2(r_0+r_1)r_2 &\ge&  r_2d_0+r_2d_1-(r_0+r_1)d_2.
\end{eqnarray*}
Note that these two inequalities bound $\alpha_2$ from below.  The
region cut out by these inequalities is sketched in Figure 1.
\begin{figure}[h]
\begin{center}
\begin{pspicture}(-4,-.5)(6.2,5)
\psset{xunit=1cm, yunit=1cm}
\pspolygon*[linecolor=lblue](-4,-.5)(-4,5)(6.2,5)(6.2,-.5)
\pspolygon*[linecolor=lightorange](-3.33,5)(6.2,5)(6.2,4.275)(.167,-.25)
\psdot[linecolor=mix, linewidth=30pt](1,3)
\uput[0](.1,3){\scaleboxto(1.5,1){\lightorange $R(t)$}}
\psline[linecolor=red, linewidth=1.5pt](-3.33,5)(.33,-.5)
\psline[linecolor=red, linewidth=1.5pt](-.167,-.5)(6.2,4.275)
\pstextpath[c](.167,-.45){\psline[linestyle=none](-.167,-.7)(6.2,4.075)}
{\red$-\alpha_1r_1r_2+\alpha_2(r_0+r_1)r_2$}
\pstextpath[c](.167,-.75){\psline[linestyle=none](-.167,-.75)(6.2,3.775)}
{\red$= r_2d_0+r_2d_1-(r_0+r_1)d_2$}
\pstextpath[c](.13,-.5){\psline[linestyle=none](-3.52,5)(.13,-.5)}
{\red$\alpha_1r_0r_1+\alpha_2r_0r_2$}
\pstextpath[l](.13,-.75){\psline[linestyle=none](-3.52,4.75)(.13,-.75)}{\red\q$=(r_1+r_2)d_0-r_0d_1-r_0d_2$}
\end{pspicture}
\end{center}
\caption{The parameter region for 3-chains.}
\end{figure}
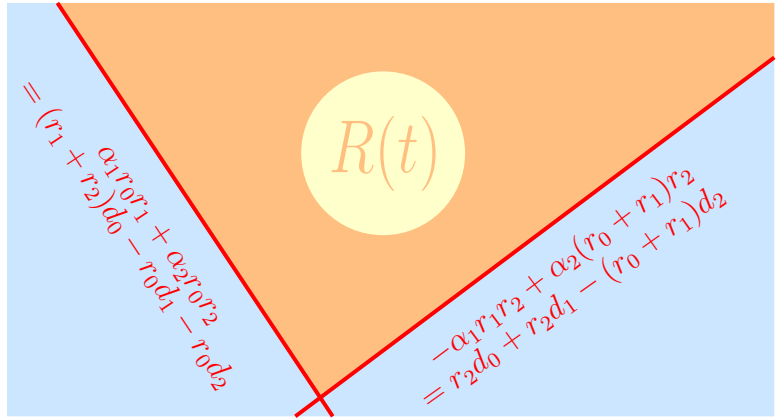
\end{example}
\begin{remark}[Degenerate holomorphic chains]
\label{MapsAreNonZero}
Fix the type $t$ and suppose we are given a stability parameter
$\alpha$ and an $\alpha$-semistable holomorphic chain
$C=(E_j,j=0,...,n;\phi_i, i=1,...,n)$, such that, say, $\phi_{i_0+1}$
is zero. Then, $(0,...,0, E_{i_0+1},...,E_n, 0,...,0,
\phi_{i_0+2},...,\phi_n)$ is both a subchain and a quotient chain. By
Remark \ref{quotientformulation}, ii), this implies that the
inequality arising from the $i_0$-th standard subchain must become an
equality, that is,
$$
\begin{array}{rl}
&\bigl(\alpha_1r_1+\cdots+\alpha_{i_0}r_{i_0}\bigr)\bigl(r_{i_0+1}+\cdots+r_n\bigr)
-
\\
& -
\bigl(\alpha_{i_0+1}r_{i_0+1}+\cdots+\alpha_nr_n\bigr)\bigl(r_0+\cdots+r_{i_0}\bigr)
\\
= & \bigl(r_0+\cdots+r_{i_0}\bigr) \bigl(d_{i_0+1}+\cdots+d_n\bigr)-
\\
&-
\bigl(r_{i_0+1}+\cdots+r_n\bigr)\bigl(d_0+\cdots+d_{i_0}\bigr).
\end{array}
$$
Therefore, the parameter $\alpha$ lies on the boundary of the
region cut out by the inequalities in Proposition
\ref{StandardBounds}. Moreover, the moduli space for $(\alpha_j,
j=0,...,n)$-semistable holomorphic chains of type$t$ for which the
inequality associated with the $i_0$-th standard subchain becomes an
equality can be easily seen to be a product of the moduli space of
$(\alpha_1,...,\alpha_{i_0})$-semistable $(i_0+1)$-chains of type
$(r_j, j=0,...,i_0; d_j, j=0,...,i_0)$ (which, for $i_0=0$, is the
moduli space of semistable vector bundles of rank $r_0$ and degree
$d_0$) and the moduli space of
$(\alpha_{i_0+1},...,\alpha_{n})$-semistable $(n+i_0)$-chains of type
$(r_j, j=i_0+1,...,n; d_j, j=i_0+1,...,n)$ (see Section
\ref{sub:FirstEx}).  The upshot is that, in our study, we may
restrict to holomorphic chains in which all homomorphisms are
non-trivial.
\end{remark}
It seems to be a very difficult problem to determine the ``exact"
shape of the parameter region in general (see Section
\ref{sec:3-chains} for the case $n=2$).  Thus, one should be more
modest and try to understand the behaviour along a half line in
$\R^n$. Here, we may offer the following result.
\begin{theorem}
\label{AlongALine}
Fix the type $t$.  Choose rational stability parameters $\beta$ and
$\gamma$ in $\Q^n$ and set
$$
\alpha^\la:= \beta+\la\cdot \gamma,\q \la\in\R_{\ge 0}.
$$
Then, there exists a value $\la_\infty$, such that for any
$\la>\la_\infty$, a holomorphic chain $C=(E_j, j=0,...,n;\phi_i,
i=1,...,n)$ is $\alpha^\la$-(semi)stable, if and only if it satisfies
the following conditions:
\begin{enumerate}
\item For every subchain $(F_j,j=0,...,n)$, the condition
\begin{equation*}
\frac{\gamma_1\rk(F_1)+\cdots+\gamma_n\rk(F_n)}
{\rk(F_0)+\cdots+\rk(F_n)}
\le
\frac{\gamma_1\rk(E_1)+\cdots+\gamma_n\rk(E_n)}
{\rk(E_0)+\cdots+\rk(E_n)}
\end{equation*}
is verified.
\item If we have equality above, then
$$
\mu_\beta(F_j,j=0,...,n; \phi_{i|F_i}, i=1,...,n) (\le)
\mu_\beta(C).
$$
\end{enumerate}
\end{theorem}
In order to appreciate the above statement, let us discuss linear
chains over the field $k$.  A \it linear $(n+1)$-chain of type
$r=(r_j, j=0,...,n)$ \rm is a tuple $V=(V_j, j=0,...,n; f_i,
i=1,...,n)$ composed of $k$-vector spaces $V_j$ with $\dim(V_j)=r_j$,
$j=0,...,n$, and linear maps $f_i\colon V_i\lra V_{i-1}$,
$i=1,...,n$.  As before, we may speak of \it subchains\rm, \it
quotient chains\rm, \it dual chains\rm, and so on. If we are given a
tuple $\alpha=(\alpha_j, j=0,...,n)$ of real numbers, we say that a
chain $V=(V_j, j=0,...,n; f_i, i=1,...,n)$ is \it
$\alpha$-(semi)stable\rm, if for every subchain $(W_j, j=0,...n)$,
the condition
$$
\frac{\alpha_0\dim(W_0)+\cdots+\alpha_n\dim(W_n)}
{\dim(W_0)+\cdots+\dim(W_n)}
(\le)
\frac{\alpha_0\dim(V_0)+\cdots+\alpha_n\dim(V_n)}
{\dim(V_0)+\cdots+\dim(V_n)}
$$
holds true. These concepts are special cases of King's general
formalism \cite{Ki}.
\begin{remark}
Any chain is $0$-semistable.
\end{remark}
Now, let $C$ be a holomorphic chain of type $t=(r_j,j=0,...,n; d_j,
j=0,...,n)$.  If $\eta$ stands for the generic point of the curve
$X$, the restriction $C_\eta$ of $C$ to the generic point is a linear
chain of type $(r_j,j=0,...,n)$ over $\CC(X)$. Condition (1) in
Theorem \ref{AlongALine} just says that $C_\eta$ is a
$\gamma$-semistable linear chain.  An immediate consequence is
\begin{corollary}
Suppose there are no $\gamma$-semistable $\CC(X)$-linear chains, then
there are no $\alpha^\la$-semistable holomorphic chains of type $t$
for $\la\gg 0$.
\end{corollary}
Therefore, the intersection of the region of parameters $\alpha$ for
which there do exist $\alpha$-semistable holomorphic chains of type
$t$ intersected with the half line $\beta+\R_{\ge 0}\cdot \gamma$ is
always bounded. Of course, one hopes that one can choose $\la_\infty$
``uniformly". We state the following.
\begin{conjecture}
\label{ParBounds}
Suppose that $0$ is the only parameter for which there exist
$\alpha$-semistable linear chains of type $(r_j,j=0,...,n)$ over
$\CC(X)$.  Then, for any type $t=(r_j, j=0,...,n; d_j, j=0,...,n)$,
there is a bounded region $R(t)\subset \R^n$, such that the existence
of an $\alpha$-semistable holomorphic chain of type $t$ implies
$\alpha\in R(t)$.
\end{conjecture}
We will prove this conjecture for $n=2$ (see Theorem
\ref{BoundedParameterRegion}).
\begin{proof}[Sketch of proof of Theorem \ref{AlongALine}]
For the proof, one has to place oneself in the more general setting
of decorated tuples of vector bundles (see \cite{Schquiv}). In our
sketch, we will use the terminology of \cite{Schquiv} without
repeating it here. Fix non-negative integers $a$, $b$, and $c$.
Then, we study tuples $(E_j, j=0,...,n;\phi)$ where the $E_j$ are
vector bundles of rank $r_j$ and degree $d_j$, $j=0,...,n$, and
$$
\phi\colon \bigl(E^{\otimes a}\bigr)^{\oplus b}\lra \det(E)^{\otimes c},\q
E:=E_0\oplus\cdots\oplus E_n,
$$
is a non-trivial homomorphism. We have a natural equivalence
relation on those objects which always identifies $\phi$ with $z\cdot
\phi$, $z\in\CC^*$.  Note that $\CC\oplus
\bigoplus_{i=1}^n\Hom(\CC^{r_i},\CC^{r_{i-1}})$ is a direct summand
of
$$
\biggl(\Bigl(\bigl(\CC^{r_0+\cdots+r_n}\bigr)^{\otimes
a}\Bigr)^{\oplus b}\biggr)\otimes
\det\bigl(\CC^{r_0+\cdots+r_n}\bigr)^{\otimes -c}
$$
for appropriate non-negative integers $a$, $b$, and $c$. Thus, the
formalism of holomorphic chains is contained in the formalism of
decorated tuples of vector bundles.  For a tuple
$\alpha=(\alpha_0,...,\alpha_n)$ of rational numbers with
$\sum_{i=0}^n \alpha_ir_i=0$ and $\la\in\Q_{>0}$, a tuple $(E_j,
j=0,...,n;\phi)$ is called \it $(\alpha,\la)$-(semi)stable\rm, if
$$
L(E^\bullet, a) -\sum_{j=0}^n r_j\cdot \biggl(\sum_{\nu=1}^s
a_\nu\cdot\Bigl(\sum_{j=0}^n \alpha_j\rk(E_j^\nu)\Bigr)\biggr)
+\la\cdot \mu\bigl(E^\bullet,a; \phi\bigr) (\ge) 0
$$
with
\begin{eqnarray*}
&&
L(E^\bullet,a):=
\\
&&\sum_{\nu=1}^s a_\nu\Biggl(\deg\Bigl(\bigoplus_{j=0}^n E_j\Bigr)\cdot
\rk\Bigl(\bigoplus_{j=0}^n E^\nu_j\Bigr)- \deg\Bigl(\bigoplus_{j=0}^n E^\nu_j\Bigr)\cdot
\rk\Bigl(\bigoplus_{j=0}^n E_j\Bigr)\Biggr)
\end{eqnarray*}
holds for every weighted filtration
\begin{eqnarray*}
(E^\bullet,a):\q \bigl(0&\subsetneq& (E_j^1, j=0,...,n)\q \subsetneq\q \cdots\q \subsetneq\q
(E_j^s, j=0,...,n)
\\
&\subsetneq& (E_j, j=0,...,n),\q a=(a_1,...,a_s)\bigr)
\end{eqnarray*}
of the ``split" vector bundle $(E_j, j=0,...,n)$.  What we want to do
is study the condition of $(\alpha^\la,\la)$-semistability for large
$\la$.  Let $V$ and $V_j$ be the fibres of $E$ and $E_j$,
respectively, over the generic point and $\sigma\in \PP((V^{\otimes
a})^{\oplus b})$ the point defined by $\phi$. Condition (1) says that
$\sigma$ must be semistable w.r.t.\ to the action of
$\SL(V)\cap(\GL(V_0)\times\cdots\times \GL(V_n))$ and its
linearization in $\OO(1)$ modified by the character corresponding to
$\gamma$.  Set
$$
\mu_{\gamma}\bigl(E^\bullet,a; \phi\bigr):=
-\sum_{j=0}^n r_j\cdot
\biggl(\sum_{\nu=1}^s a_\nu\cdot\Bigl(\sum_{j=0}^n \gamma_j\rk(E_j^\nu)\Bigr)\biggr)
+\mu\bigl(E^\bullet,a; \phi\bigr).
$$
If $\sigma$ fails to be semistable, one applies the theory of the
instability flag in order to produce a weighted filtration
$(E^\bullet, a)$ with $\mu_{\gamma}\bigl(E^\bullet,a; \phi\bigr)<0$
and $L(E^\bullet, a)$ bounded from above by a constant which depends
only on the type $t$, $a$, $b$, and $c$.  The details for this may be
easily adapted from our paper \cite{SchGlobal}.  Moreover, the
possible tuples $a$ and $(\rk(\bigoplus_{j=0}^n E_j^\nu),
\nu=1,...,s)$ belong to a finite set whence
$$
-\sum_{j=0}^n r_j\cdot\biggl(\sum_{\nu=1}^s a_\nu\cdot\Bigl(\sum_{j=0}^n \beta_j\rk(E_j^\nu)\Bigr)\biggr)
$$
may also be bounded from above by a constant which depends only on
the type $t$, $a$, $b$, and $c$.  It is now clear that we may find
$\la_0$, such that for $\la>\la_0$ an $(\alpha^\la,\la)$-semistable
tuple satisfies
\begin{enumerate}
\item $\mu_{\gamma}\bigl(E^\bullet,a; \phi\bigr)\ge 0$ for every
weighted filtration $(E^\bullet, a)$ (which is equivalent to saying
that $\sigma$ is semistable).
\item If ``=" holds, then
$$
L\bigl(E^\bullet,a; \phi\bigr)-\sum_{j=0}^n r_j\cdot
\biggl(\sum_{\nu=1}^s a_\nu\cdot\Bigl(\sum_{j=0}^n \beta_j\rk(E_j^\nu)\Bigr)\biggr)(\ge) 0.
$$
\end{enumerate}
The fact that Condition (1) and (2) are also sufficient will follow
easily once one knows that the tuples $(E_j,j=0,...,n;\phi)$ of type
$t$, satisfying (1) and (2), live in a bounded family. This is again
established along the lines of the corresponding result in
\cite{SchGlobal}.
\end{proof}
\subsection{Walls and the chamber structure}
\label{sub:chambers-walls}
In this section, we would like to subdivide $\R^n$ into locally
closed subsets, called chambers, such that the concept of
$\alpha$-(semi)stability is constant within each chamber.
\par
We fix the type $t=(r_j, j=0,...,n; d_j, j=0,...,n)$ and set
$r:=r_0+\cdots+ r_n$ and $d:=d_0+\cdots+d_n$. For a holomorphic chain
$C=(E_j, j=0,...,n;\phi_i, i=1,...,n)$, the \it total rank \rm is
given by $r(C):=\rk(E_0)+\cdots+\rk(E_n)$ and the \it total degree
\rm by $d(C):=\deg(E_0)+\cdots+\deg(E_n)$.  The idea is to first
define hyperplanes which cut out parameters for which there might
exist properly $\alpha$-semistable (i.e., semistable but not stable)
holomorphic $(n+1)$-chains of type $t$.  Suppose $C$ is such a chain
and $C^\p=(F_j, j=0,...,n)$ is a destabilizing subchain.  Then, with
$\alpha$ the stability parameter in question, we obtain the equation
\begin{eqnarray*}
&&rd(C^\p)-r(C^\p)d
\\
&=&\alpha_1\bigl(r_1r(C^\p)-\rk(F_1)r\bigr)+\cdots+
\alpha_n\bigl(r_nr(C^\p)-\rk(F_n)r\bigr).
\end{eqnarray*}
Define
\begin{eqnarray*}
{\mathcal S}:=\Bigl\{\,(s_0,...,s_n; e)&\big|& 0\le s_j\le r_j,\ j=0,...,n,\
0<s< r,
\\
&& s:=s_0+\cdots+s_n,\ e\in\Z \,\Bigr\}.
\end{eqnarray*}
For an element $\sigma\in{\mathcal S}$, let
$$
{w}_\sigma:= \Bigl\{\, \alpha\in\R^n\, \big|\,
\alpha_1(r_1s-s_1r)+\cdots+
\alpha_n(r_ns-s_nr)=re-sd\,\Bigr\}
$$
be the \it wall defined by $\sigma$\rm. Note that we may have an
empty wall ${w}_\sigma=\varnothing$ or an \it improper wall \rm
${w}_\sigma=\R^n$.  For an improper wall, we must have $r_js-rs_j=0$,
$j=1,...,n$, and $re-sd=0$.  Set $r_{\rm red}:=r/\gcd(r,s)>1$,
because $s<r$. Then, $r_{\rm red}$ divides $r_i$, $i=1,...,n,$ and
$d$.
\begin{proposition}
\label{CoPrimalityI}
If $\gcd(r_1,....,r_n,d)=1$, then there do not exist any improper
walls.
\end{proposition}
The \it $n$-dimensional chambers \rm are given as the connected
components of
$$
\R^n\setminus\bigcup_{\sigma\in{\mathcal S}: {w}_\sigma\neq \R^n} {w}_\sigma.
$$
The \it $(n-1)$-dimensional chambers \rm are given as the
connected components of
$$
\bigcup_{\sigma\in{\mathcal S}: {w}_\sigma\neq \R^n} \biggl({w}_\sigma\setminus
\bigcup_{\tau\in{\mathcal S}: {w}_\sigma\not\subseteq {w}_\tau} {w}_\tau \biggr),
$$
and so on. We label the $j$-dimensional chambers ${\mathcal
C}^j_k$, $k\in J_j$, $j=0,...,n$. Note that we have
\begin{equation}
\label{ChamberDecomposition}
\R^n=\bigsqcup_{j=0}^n \bigsqcup_{k\in J_j} {\mathcal C}^j_k.
\end{equation}
Observe also that this chamber decomposition is locally finite, i.e.,
every bounded subset $R\subset\R^n$ intersects only finitely many
chambers.  By construction, we have the following property.
\begin{proposition}
\label{ChamberDecompositionII}
{\rm i)} Let ${\mathcal C}$ be any chamber and $\alpha_1,\alpha_2\in
{\mathcal C}$.  Then, a holomorphic chain of type $t$ is
$\alpha_1$-(semi)stable, if and only if it is
$\alpha_2$-(semi)stable.
\par
{\rm ii)} Let ${\mathcal C}_1$ be any chamber and ${\mathcal C}_2$ a
chamber in the closure of ${\mathcal C}_1$. Choose $\alpha_i\in
{\mathcal C}_i$, $i=1,2$.  Then, a holomorphic chain which is
$\alpha_1$-semistable is also $\alpha_2$-semistable, and a
holomorphic chain which is $\alpha_2$-stable is also
$\alpha_1$-stable.
\end{proposition}
\begin{remark}
\label{rationality1}
By definition, any chamber contains elements of $\Q^n$, so that it
suffices to consider rational stability parameters.
\end{remark}
Finally, we note the following consequence of Proposition
\ref{CoPrimalityI}.
\begin{corollary}
\label{CoPrimalityII}
If $\gcd(r_1,....,r_n,d)=1$, then, for a stability parameter $\alpha$
which lies in an $n$-dimensio\-nal chamber, the conditions of
$\alpha$-stability and $\alpha$-semistability coincide.
\end{corollary}
\begin{proof}
The definition of the walls shows that a stability parameter $\alpha$
for which there exists a properly $\alpha$-semistable holomorphic
chain of type $t$ must lie on a wall. Since the assumption grants
that no improper wall exists, we are done.
\end{proof}
To conclude this paragraph, we remark that the chamber structure has
been obtained by rough a priori considerations. Usually, one expects
a coarser chamber structure based on a refined analysis of stability.
This seems, unfortunately, rather involved and may be carried out
only in more specialized situations.  In general, we would expect
only finitely many chambers, even if the region $R(t)$ of possible
stability parameters were not bounded. This and other phenomena will
be explained for $n=2$ in Section \ref{sec:3-chains}.
\subsection{Vortex equations and Hitchin--Kobayashi correspondence}
\label{sub:vortex-equations}
There are natural gauge-theoretic equations on a  holomorphic chain
$$
C \colon E_{n} \overset{\phi_n}{\lto} E_{n-1}
\overset{\phi_{n-1}}{\lto} \cdots \overset{\phi_{1}}{\lto} E_{0},
$$
which we describe now.
Define $\tau = (\tau_{0}, \ldots, \tau_{n})
\in \R^{n+1}$ by
\begin{equation}
  \label{eq:tau-alpha}
  \tau_{j} = \mu_{\alpha}(C) - \alpha_{j}, \quad j=0,\ldots, n,
\end{equation}
where we make the convention $\alpha_0 = 0$.  Then $\alpha$ can be
recovered from $\tau$ by
\begin{equation}
  \label{eq:alpha-tau}
  \alpha_{j} = \tau_{0} - \tau_{j}, \quad j=0,\ldots, n.
\end{equation}
The
\emph{$\tau$-vortex equations}
\begin{equation}
  \label{eq:vortex-eq}
  \sqrt{-1}\Lambda F(E_{j}) + \phi_{j+1} \phi_{j+1}^* - \phi_{j}^*\phi_{j}
  = \tau_{j} {\Id}_{E_{j}}, \quad j=0,\ldots, n,
\end{equation}
are equations for Hermitian metrics on $E_{0}, \ldots, E_{n}$.  Here,
$F(E_{j})$ is the curvature of the Hermitian connection on $E_{j}$,
$\Lambda$ is contraction with the K\"ahler form of a fixed metric on $X$ such that
$\vol(X) = 2\pi$, and $\phi_j^*$ is the adjoint of $\phi_j$.
By convention $\phi_0 = \phi_{n+1}=0$.
One has the following.
\begin{theorem}[{\cite[Theorem 3.4]{AG1}}]
  \label{thm:chain-hitchin-kobayashi}
  A holomorphic chain $C$ is $\alpha$-poly\-stable if and
  only if the $\tau$-vortex equations have a solution, where
  $\alpha$ and $\tau$ are related by \eqref{eq:tau-alpha}.
\end{theorem}
\subsection{Moduli spaces for parameters on and near the standard hyperplanes}
\label{sub:FirstEx}
A standard procedure to study moduli spaces is to start with known moduli
spaces and create new ones out of them by ``flip-type" operations.
In our setting, we might try a kind of inductive procedure, by relating
moduli spaces of holomorphic $(n+1)$-chains to moduli of ``shorter" holomorphic
chains. This is indeed possible for stability parameters in or near the
hyperplanes where the inequalities in Proposition \ref{StandardBounds} become
equalities.
\par
Let $h_i$, $i=0,...,n$, be the hyperplanes that were defined by the equations (\ref{eq:standard-hyperplanes}).
\begin{proposition}
\label{prop:braking}
Let $C=(E_j, j=0,...,n; \phi_i, i=1,...,n)$ be a holomorphic
$(n+1)$-chain.
\par
{\rm i)} Assume $C$ to be $\alpha$-semistable.  If $\alpha\in
h_{i_0+1}$, then $C$ is S-equivalent to $(E_j, j=0,...,n;
\widetilde{\phi}_i, i=1,...,n)$ with $\widetilde{\phi}_{i_0+1}=0$ and
$\widetilde{\phi}_i=\phi_i$, for $i\neq {i_0+1}$. In particular, if
$C$ is $\alpha$-polystable, then $\phi_{i_0+1}=0$.
\par
{\rm ii)} The $({i_0+1})$-chain $C'=(E_j, j=0,...,i_0; \phi_i,
i=1,..., i_0)$ is $\beta$-semistable for
$\beta=(\alpha_0,...,\alpha_{i_0})$, and the $(n-{i_0})$-chain
$C''=(E_{i_0+1+j},\allowbreak j=0,..., n-{i_0}-1; \phi_{i_0+1+i},
i=1,..., n-i_0-1)$ is $\gamma$-semistable for
$\gamma=(\alpha_{i_0+1},...,\alpha_n)$.  If, furthermore, $C$ is
$\alpha$-polystable, then $C'$ and $C''$ are $\beta$- and
$\gamma$-polystable, respectively.
\par
{\rm iii)} If, conversely, the chains $C'$ and $C''$ are $\beta$- and
$\gamma$-semistable (polystable), then $C$ is $\alpha$-semistable
(polystable).
\end{proposition}
\begin{proof}
The arguments are essentially the same as for vector bundles, working
in the Abelian category of holomorphic chains with the
$\alpha$-degree, the total rank, the $\alpha$-slope, and the notion
of $\alpha$-semistability as the semistability concept.
\par
If $\alpha\in h_{i_0+1}$, then $\widetilde{C}:= (E_0,...,E_{i_0},
0,...,0; \phi_1,...,\phi_{i_0},0,...,0)$ is a subchain with
$\mu_\alpha(\widetilde{C})=\mu_\alpha(C)$. The quotient chain
$\ol{C}:=C/\widetilde{C}$ is $(0,...,0,E_{i_0+1},...,E_{n};
\phi_1,...,\phi_{i_0},0,...,0)$, and
$\mu_\alpha(C)=\mu_\alpha(\ol{C})$.
\par
i) By definition of S-equivalence, $C$ is S-equivalent to
$\widetilde{C}\oplus \ol{C}$, and, for $\alpha$-polystable chains,
S-equivalence is the same as equivalence.
\par
ii) and iii) Standard arguments (parallel to those for semistable
vector bundles) show that $C$ is $\alpha$-semistable, if and only if
both $\widetilde{C}$ and $\ol{C}$ are $\alpha$-semistable. Now, note
that $C'$ is $\beta$-semistable, if and only if $\widetilde{C}$ is
$\alpha$-semistable and that $C''$ is $\gamma$-semistable, if and
only if $\ol{C}$ is $\alpha$-semistable. This proves the assertions
on semistability in ii) and iii). The corresponding claims about
polystability are left as an exercise to the reader.
\end{proof}
\begin{remark}
All the above observations regarding polystability may also be easily
derived from the existence of solutions to the vortex equations
(\ref{eq:vortex-eq}) on the chain $C$.
\end{remark}
\begin{corollary}
\label{cor:Start}
{\rm i)} Let $\alpha\in h_{i_0+1}$. With the same notation as in the
above proposition, we have that
$$
\mathcal{M}^{\rm red}_\alpha(t)\cong \mathcal{M}^{\rm red}_\beta(t')\times
\mathcal{M}^{\rm red}_\gamma(t''),
$$
where $t'$ and $t''$ are the types of $C'$ and $C''$,
respectively.  Here, the superscript ``{\rm red}" refers to the
induced reduced scheme structure.
\par
{\rm ii)}
Let $\{\alpha\}=\bigcap_{i=1}^{n}h_i$. Then
$$
\mathcal{M}^{\rm red}_\alpha(t)\cong
\mathcal{M}(r_1,d_1)\times...\times \mathcal{M}(r_n,d_n)
$$
with $\mathcal{M}(r,d)$ the moduli space of semistable vector
bundles of rank $r$ and degree $d$.
\end{corollary}
\begin{proposition}
Fix the type $t$, and suppose $\alpha\in h_{i_0+1}$. Let ${\mathcal
C}$ be any chamber, such that $\alpha\in \ol{\mathcal C}$, and choose
$\kappa\in{\mathcal C}$.
\par
{\rm i)} If $C$ is $\kappa$-semistable, then the chain $C'$ is
$\beta$-semistable and the chain $C''$ is $\gamma$-semistable.
\par
{\rm ii)} If $C'$ is $\beta$-stable and $C''$ is $\gamma$-stable,
then, for any $\phi_{i_0+1}\colon E_{i_0+1}\to E_{i_0}$ different
from zero, the resulting chain $C$ is $\kappa$-stable.
\end{proposition}
\begin{proof}
i) This is a trivial continuity statement, observing the discussions
in the proof of Proposition \ref{prop:braking}.
\par
 ii) We will demonstrate the following property:
\par\it
There is an open subset ${\mathcal U}\subset \R^n$, containing
$\alpha$, such that, for $\kappa\in {\mathcal U}\cap
(H_{i_0+1}\setminus h_{i_0+1})$, we have the following. If $C'$ is
$\beta$-stable and $C''$ is $\gamma$-stable, then, for any
$\phi_{i_0+1}\colon E_{i_0+1}\to E_{i_0}$ different from zero, the
resulting chain $C$ is $\kappa$-stable.
\par
\rm In view of the general properties of the chamber decomposition
(Pro\-position \ref{ChamberDecompositionII}), this will imply the

assertion of the proposition.
\par
Define
${\mathcal S}$ as in Section \ref{sub:chambers-walls}, and
${\mathcal S}^{\rm real}$ as the set of elements of ${\mathcal S}$
which come from a subchain of a holomorphic chain $C$ of type $t$,
such that $C'$ is $\beta$-stable and $C''$ is $\gamma$-stable.
Declare the finite set
\begin{eqnarray*}
{\mathcal R}:=\Bigl\{\, (s_0,...,s_n)&\big|& 0\le s_j\le r_j,\ j=0,...,n,
\\
&&\bigl(0< s_0+\cdots+s_{i_0} < r_0+\cdots+ r_{i_0}\bigr)\vee
\\
&& \vee \bigl(0< s_{i_0+1}+\cdots+ s_n < r_{i_0+1}+\cdots+ r_n\bigr) \,\Bigr\}.
\end{eqnarray*}
For $s\in {\mathcal R}$, we define ${\mathcal S}^{\rm real}_s$ as the set of
elements $(s_0,...,s_n, e)\in {\mathcal S}^{\rm real}$ with $(s_0,...,s_n)=s$, and set
$$
e_0(s):=\max\bigl\{\, e\,|\, (s,e)\in {\mathcal S}^{\rm real}_s\,\bigr\}.
$$
For each $s\in {\mathcal R}$, there is the function
\begin{eqnarray*}
\chi_{s}\colon \R^n &\lra& \R^n
\\
(\kappa_1,...,\kappa_n) &\lma& \sum_{j=0}^n s_j\cdot \sum_{j=0}^{n}
\kappa_j r_j
-\sum_{j=0}^n r_j\cdot \sum_{j=0}^{n} \kappa_i s_i.
\end{eqnarray*}
Obviously, we can define open subsets ${\mathcal U}_s$ by the condition
$\chi_{s}(\kappa)>e_0(s)(r_0+\cdots+r_n)-d(s_0+\cdots+s_n)$
for any $\kappa\in {\mathcal U}_s$, $s\in {\mathcal R}$.
We set
$$
{\mathcal U}:=\bigcap_{s\in {\mathcal R}} {\mathcal U}_s.
$$
\par
Let $\widetilde{C}$ be any proper subchain of $C$. We view $C'$ as a
subchain and $C''$ as a quotient chain of $C$.  Suppose that both
$\widehat{C}:=\widetilde{C}\cap C'$ and the projection $\ol{\ol{C}}$
of $\widetilde{C}$ to $C''$ are trivial subchains of $C'$ and $C''$,
respectively.  Because $\phi_{i_0+1}$ is non-trivial, this can happen
only for $\widetilde{C}=0$, $\widetilde{C}=C$, and
$\widetilde{C}=C'$. In the third case, $\mu_\kappa(\widetilde{C})<
\mu_\kappa(C)$ is equivalent to the assumption $\kappa\in
H_{i_0+1}\setminus h_{i_0+1}$.  Thus, we may assume that
$\widehat{C}$ is a non-trivial subchain of $C'$ or that $\ol{\ol{C}}$
is a non-trivial subchain of $C''$.  But then,
$\mu_\kappa(\widetilde{C})<\mu_\kappa(C)$, for all $\kappa\in
{\mathcal U}$, follows immediately from the definition of ${\mathcal
U}$, because $(\rk(F_0),...,\rk(F_n))\in {\mathcal R}$.
\par
We still have to show that ${\mathcal U}$ contains $\alpha$.  For
given $s\in {\mathcal R}$, we may choose a chain $C$ of type $t$,
such that $C'$ is $\beta$-stable and $C''$ is $\gamma$-stable and a
subchain $\widetilde{C}$ of $C$ with $(\rk(F_0),...,\rk(F_n))=s$ and
$\deg(F_0)+\cdots+\deg(F_n)=e_0(s)$.  Then,
\begin{eqnarray*}
&&
\mu_\alpha(\widetilde{C})
\\
&=& \frac{\sum_{j=0}^{n} \deg(F_j)+\alpha_1\rk(F_1)+\cdots+
\alpha_n\rk(F_n)}{\rk(F_0)+\cdots+\rk(F_n)}
\\
&=&\frac{\sum_{j=0}^{i_0} \rk(F_j)\cdot \mu_{(0,\alpha_1,...,\alpha_{i_0})}(\widehat{C})
+\sum_{j=i_0+1}^{n} \rk(F_j)\cdot\mu_{(\alpha_{i_0+1},...,\alpha_n)}(\ol{\ol{C}})}
{\rk(F_0)+\cdots+\rk(F_n)}
\\
&<& \frac{\sum_{j=0}^{i_0} \rk(F_j)\cdot\mu_{(0,\alpha_1,...,\alpha_{i_0})}(C')
+\sum_{j=i_0+1}^{n} \rk(F_j)\cdot \mu_{(\alpha_{i_0+1},...,\alpha_n)}(C'')}
{\rk(F_0)+\cdots+\rk(F_n)}
\\
&=& \frac{\sum_{j=0}^{i_0} \rk(F_j)\cdot\mu_{\alpha}(C)
+\sum_{j=i_0+1}^{n} \rk(F_j)\cdot\mu_{\alpha}(C)}
{\rk(F_0)+\cdots+\rk(F_n)}
\\
&=&\mu_\alpha(C).
\end{eqnarray*}
This inequality is equivalent to $\alpha\in {\mathcal U}_s$.
\end{proof}
\begin{remark}
The same game can be played with a parameter $\alpha$ lying on
several of the hyperplanes $h_i$.  In particular, we can apply it to
the parameter $\alpha$ defined in Corollary \ref{cor:Start}, ii).
\end{remark}

\section{Extensions and deformations of chains}
\label{sec:extensions}

As a first step to study the variation of the moduli spaces of
$\alpha$-semistable holomorphic chains as the parameter $\alpha$
changes, we study in this section the deformation theory of
holomorphic chains.  As in the case of holomorphic triples, which is
treated in \cite[\S 3]{BGG}, the infinitesimal deformations of
holomorphic chains are given by the hypercohomology groups of certain
sheaf complexes associated to the holomorphic chains.

Througout this section we fix a stability parameter $\alpha = (\alpha_j,
j=0,...,n)$ and two holomorphic chains $C'$ and $C''$, of types
$t'=(r_j',j=0,...,n; d_j', j=0,...,n)$ and $t''=(r_j'',j=0,...,n; d_j'',
j=0,...,n)$, respectively, given by
\begin{equation}
\label{eq:C'}
C' \colon E_n' \overset{\phi_n'}{\lto} E_{n-1}'
  \overset{\phi_{n-1}'}{\lto} \cdots
  \overset{\phi_{1}'}{\lto} E_{0}' ,
\end{equation}
\begin{equation}
\label{eq:C''}
C'' \colon E_n'' \overset{\phi_n''}{\lto} E_{n-1}''
  \overset{\phi_{n-1}''}{\lto} \cdots
  \overset{\phi_{1}''}{\lto} E_{0}''.
\end{equation}
Given two vector bundles $E$ and $F$ over $X$, $\Hom(E,F)$ and $\End
(E)$ denote the vector bundles of homomorphisms from $E$ to $F$ and
of endomorphisms of $E$, respectively. The corresponding spaces of
global sections are denoted $\Hom_X(E,F)$ and $\End_X (E)$,
respectively.

\subsection{Hypercohomology}
\label{sub:hypercohomology}

In this subsection we analyze $\Ext^1(C'',C')$ using the
hypercohomology groups $\HH^i(F^\bullet(C'',C'))$ of a 2-step complex
of vector bundles
\begin{equation}
\label{eq:extension-complex}
  F^\bullet(C'',C')\colon F^0 \overset{b}{\lto} F^1.
\end{equation}
This complex has terms
\[
  F^0=\bigoplus_{i=0}^n \Hom(E_i'',E_i'), \qquad
  F^1=\bigoplus_{i=1}^n \Hom(E_i'',E_{i-1}'),
\]
and differential
\[
  b(\psi_0,\ldots,\psi_n)=\sum_{i=1}^n b_i(\psi_{i-1},\psi_i),
  \qquad \textnormal{ for } \psi_i\in\Hom(E_i'',E_i'),
\]
where
\[
  b_i\colon \Hom(E_{i-1}'',E_{i-1}')\oplus
\Hom(E_i'',E_i')\to \Hom(E_i'',E_{i-1}')\hra F^1
\]
is given by
\[
  b_i(\psi_{i-1},\psi_i)=\psi_{i-1}\circ\phi_i''-\phi_i'\circ\psi_i.
\]
Applying the cohomology functor to this complex of vector bundles, we
obtain maps of vector spaces
\begin{equation}
\label{eq:H^*}
 d=H^p(b)\colon  H^p(F^0) \lto H^p(F^1),
\end{equation}
for $p=0,1$, where
\[
  H^p(F^0)=\bigoplus_{i=0}^n \Ext^p_X(E_i'',E_i'), \qquad
  H^p(F^1)=\bigoplus_{i=1}^n \Ext^p_X(E_i'',E_{i-1}'),
\]
and
\[
  d(\psi_0,\ldots,\psi_n)=\sum_{i=1}^n d_i(\psi_{i-1},\psi_i),
  \qquad \text{for } \psi_i\in\Ext^i_X(E_i'',E_i').
\]
Here,
\begin{equation}
\label{eq:d_i}
  d_i\colon
  \begin{array}{c}
  \Ext^p_X(E_{i-1}'',E_{i-1}')\\
  \oplus\\
  \Ext^p_X(E_i'',E_i')
  \end{array}
  \to \Ext^p_X(E_i'',E_{i-1}')
  \hra H^p(F^1)
\end{equation}
is given by
\[
  d_i(\psi_{i-1},\psi_i)=\psi_{i-1}\circ\phi_i''-\phi_i'\circ\psi_i,
\]
where $\circ$ is composition of maps for $p=0$, and the Yoneda
product for $p=1$.

The following result generalizes \cite[Proposition~3.1]{BGG} from
holomorphic triples to holomorphic chains of arbitrary length.

\begin{proposition}
\label{prop:hypercohology}
  There are natural isomorphisms
\begin{align*}
  \Hom(C'',C') &\cong \HH^{0}(F^{\bullet}(C'',C')), \\
  \Ext^{1}(C'',C') &\cong \HH^{1}(F^{\bullet}(C'',C')),
\end{align*}
  and an exact sequence
\begin{eqnarray*}
0 & \lto  & \HH^0(F^{\bullet}(C'',C')) \lto H^0(F^0) \overset{d}{\lto} H^0(F^1) \lto
\HH^1(F^{\bullet}(C'',C')) \\ & \lto & H^1(F^0) \overset{d}{\lto} H^1(F^1) \lto
\HH^2 (F^{\bullet}(C'',C')) \lto 0.
\end{eqnarray*}
\end{proposition}

\begin{proof}
  This follows from \cite[Theorem~4.1 and Theorem~5.1]{GK}, since a
  holomorphic chain is a holomorphic quiver bundle, for the quiver
  $n \lto n-1 \lto \cdots \lto 1 \lto 0$.
\end{proof}

Given two sheaves or vector bundles $E$ and $F$, we define
$h^i(E,F)=\dim (\Ext^i_X(E,F))$ and $\chi(E,F)=h^0(E,F)-h^1(E,F)$.
Similarly, for any pair of chains $C'$ and $C''$ as before, we define
\begin{align*}&
  h^i(C'',C')=\dim (\HH^i(C'',C')), \\ &
  \chi(C'',C')=h^0(C'',C')-h^1(C'',C')+h^2(C'',C').
\end{align*}
Recall that
\[
  r_i'=\rk (E_i'), \quad d_i'=\deg (E_i'), \quad
  r_i''=\rk (E_i''), \quad d_i''=\deg (E_i'').
\]

\begin{proposition}
\label{prop:chi(C'',C')}
\begin{align*}
  \chi(C'',C')
  & = \sum_{i=0}^n \chi(E_i'',E_i') -  \sum_{i=1}^n \chi(E_i'',E_{i-1}') \\
  & = (1-g) \left( \sum_{i=0}^n r_i'' r_i' - \sum_{i=1}^n r_i'' r_{i-1}' \right) \\
  & \qquad  + \sum_{i=0}^n (r_i'' d_i' -r_i' d_i'') - \sum_{i=1}^n (r_i'' d_{i-1}' - r_{i-1}' d_i'').
\end{align*}
\end{proposition}

\begin{proof}
 The first equality is a consequence of the exact sequence in
 Proposition~\ref{prop:hypercohology}. The Riemann--Roch formula 
\[
  \chi(E,F)=(1-g) \rk (E) \rk (F) + \rk (E) \deg (F) - \rk (F) \deg (E),
\]
 for vector bundles $E$ and $F$, implies now the second equality.
\end{proof}

The previous proposition shows that $\chi(C'',C')$ only depends on
the types $t'$ and $t''$ of $C'$ and $C''$, respectively, so we may
use the notation
\[
  \chi(t'',t'):=\chi(C'',C').
\]

\begin{corollary}
\label{cor:chi-relation}
  For any extension $0\to C'\to C \to C''\to 0$ of holomorphic chains,
\[
  \chi(C,C)=\chi(C',C')+ \chi(C'',C'')+ \chi(C'',C')+ \chi(C',C'').
\]
\end{corollary}

\subsection{Vanishing of $\HH^0$ and $\HH^2$}
\label{sub:vanishing-hypercohology}

The following result is proved as in the case of semistable vector
bundles, given the identification of $\HH^0(F^\bullet(C'',C'))$ with
$\Hom(C'',C')$ of Proposition~\ref{prop:hypercohology}.

\begin{proposition}
\label{prop:h0-vanishing}
  Suppose that $C'$ and $C''$ are $\alpha$-semistable.
\par
{\rm i)} If $\mu_\alpha(C')<\mu_\alpha(C'')$, then $\HH^0(F^\bullet(C'',C'))=0$.
\par
{\rm ii)} If $\mu_\alpha(C')=\mu_\alpha(C'')$ and $C''$ is
$\alpha$-stable, then
 \[
     \HH^0(F^\bullet(C'',C')) \cong
     \begin{cases}
       \C \quad &\text{if $C' \cong C''$,} \\
       0 \quad &\text{if $C' \not\cong C''$}.
     \end{cases}
\]
  \qed
\end{proposition}

In the following result, $\{ u_0,\ldots,u_n\}$ is the standard basis
of $\R^{n+1}$.

\begin{proposition}
\label{prop:h2-vanishing}
Let $C'$ and $C''$ be two holomorphic chains.
\par
{\rm i)}
 Let $\cD\subset\{ 1,\ldots,n\}$ and for each $i\in\cD$, let
 $\epsilon_i\geq 0$.
 Suppose that the following three conditions hold:
\begin{itemize}
\item
 $C',C''$ are $\alpha$-semistable with $\mu_\alpha(C')=\mu_\alpha(C'')$.
\item
 For all $i\in\{ 1,\ldots,n\}\setminus\cD$, $\alpha_i-\alpha_{i-1}>
 2g-2$.
\item
 For all $i\in\cD$, one of $C',C''$ is $(\alpha+\epsilon_i
 u_i)$-stable and $\alpha_i-\alpha_{i-1}\geq 2g-2$.
\end{itemize}
 Then $\HH^2(F^\bullet(C'',C'))=0$.
\par
{\rm ii)}
 If $C'$ and $C''$ are $\alpha$-semistable with the same
 $\alpha$-slope, and $\alpha_i-\alpha_{i-1}>2g-2$ for all
 $i=1,\ldots,n$, then $\HH^2(F^\bullet(C'',C'))=0$.
\par
{\rm iii)}
 If one of $C', C''$ is $\alpha$-stable and the other one is
 $\alpha$-semistable with the same $\alpha$-slope, and $\alpha_i -
 \alpha_{i-1} \geq 2g-2$ for all $i=1,\ldots,n$, then
 $\HH^2(F^\bullet(C'',C'))=0$.
\par
{\rm iv)}
 If for all $i=1,\ldots,n$, $\phi_i''$ is injective or $\phi_i'$ is
 generically surjective,
 then $\HH^2(F^\bullet(C'',C'))=0$.
\end{proposition}

\begin{proof}
 From the exact sequence of Proposition~\ref{prop:hypercohology},
 $\HH^2(F^\bullet(C'',C'))=0$ if and only if the map $d\colon
 H^1(F^0)\to H^1(F^1)$ (defined as in~\eqref{eq:H^*} for $p=1$) is
 surjective, that is, the maps $d_i$ in~\eqref{eq:d_i} are surjective
 for all $i=1,\ldots,n$. Using Serre duality, $d_i$ is surjective if
 and only if the map
\[
  P_i\colon \! \Hom_X(E_{i-1}',E_{i}''\otimes K) \!\to\! 
\Hom_X(E_{i-1}'',E_{i-1}'\otimes K)\oplus \Hom_X(E_i',E_i''\otimes K)
\]
given by
\[
  P_i(\xi_i)=((\phi_i''\circ{\Id}_K)\circ\xi_i,\xi_i\circ\phi_i'),
  \quad \text{ for } \xi_i\in\Hom_X(E_{i-1}',E_{i}''\otimes K),
\]
is injective. Let $i\in\{1,\ldots,n\}$ and $\xi_i\colon E_{i-1}'\to
E_i''\otimes K$ be a map in $\ker (P_i)$. Let
\[
 I_i =\im (\xi_i)\subset E_i''\otimes K \qquad\text{and}\qquad
 N_{i-1} =\ker (\xi_i) \subset E_{i-1}'.
\]
Then the fact that $\xi_i\in\ker (P_i)$ is equivalent to the fact that
the maps
\begin{equation}
\label{eq:zero-maps}
E_{i-1}'\overset{\xi_i}{\lto} E_i''\otimes K
\overset{\phi_i''\otimes\Id}{\lto} E_{i-1}''\otimes K
\text{ and }
E_{i}'\overset{\phi_i'}{\lto} E_{i-1}' \overset{\xi_i}{\lto} E_i''\otimes K
\end{equation}
are both zero. The first map is zero if and only if
\begin{equation}
\label{eq:I-ker}
 I_i\subset\ker(\phi_i''\otimes{\Id}),
\end{equation}
i.e., $ I_i\otimes K^*\subset\ker (\phi_i'')$, so the diagram
\[
  \xymatrix{
  I_i\otimes K^* \ar[r] \ar @{^{(}->}[d]
  & 0 \ar @{^{(}->}[d]
  \\
  E_i'' \ar[r]^-{\phi_i''} & E_{i-1}''}
\]
commutes. In other words, there is a subchain $C_{I_i}''\hra C''$
given by
\[
\xymatrix@C=20pt{
  **[r] C_{I_i}''\colon \ar@{^{(}->}[d] & 0 \ar[r]\ar@{^{(}->}[d] & \cdots
  \ar[r] & 0 \ar[r] \ar@{^{(}->}[d] &  I_i\otimes K^* \ar[r]
  \ar[d]\ar@{^{(}->}[d]  & 0 \ar[d]\ar@{^{(}->}[d]\ar[r] & \cdots
  \ar[r] & 0 \ar@{^{(}->}[d] \\
  **[r] C''\colon & E_n'' \ar[r]^-{\phi_n''}  & \cdots \ar[r]^-{\phi_{i+2}''}  & E_{i+1}'' \ar[r]^-{\phi_{i+1}''} & E_i'' \ar[r]^-{\phi_{i}''}& E_{i-1}'' \ar[r]^-{\phi_{i-1}''} & \cdots \ar[r]^-{\phi_{1}''} & E_0''.}\]
Similarly, the second map in \eqref{eq:zero-maps} is zero if and only if
\begin{equation}
\label{eq:im-N}
  \im (\phi_i')\subset N_{i-1},
\end{equation}
i.e., the diagram
\[
  \xymatrix{
  E_i' \ar[r]^-{\phi_i'} \ar @{=}[d]
  & N_{i-1} \ar @{^{(}->}[d]
  \\
  E_i' \ar[r]^-{\phi_i'} & E_{i-1}'}
\]
commutes, so we can define a subchain $C_{N_{i-1}}'\hra C'$ by
\[
\xymatrix@C=19pt{
  **[r] C_{N_{i-1}}'\colon \ar@{^{(}->}[d]  & E_n' \ar[r]^-{\phi_n'} \ar@{=}[d] & \cdots
  \ar[r]^-{\phi_{i+1}'} & E_{i}' \ar[r]^-{\phi_{i}'} \ar@{=}[d] &  N_{i-1} \ar[rr]^-{\phi_{i-1}'|_{N_{i-1}}}
  \ar@{^{(}->}[d]  && E_{i-2}' \ar@{=}[d]\ar[r]^-{\phi_{i-2}'} & \cdots
  \ar[r]^-{\phi_{1}'} & E_0' \ar@{=}[d] \\
  **[r] C'\colon & E_n' \ar[r]^-{\phi_n'}  & \cdots \ar[r]^-{\phi_{i+1}'}  & E_{i}' \ar[r]^-{\phi_{i}'} & E_{i-1}' \ar[rr]^-{\phi_{i-1}'}&& E_{i-2}' \ar[r]^-{\phi_{i-2}'} & \cdots \ar[r]^-{\phi_{1}'} & E_0'.}\]
Let $k_{i-1}=\rk (N_{i-1})$, $l_{i-1}=\deg (N_{i-1})$, so
\[
  \mu_\alpha(C_{N_{i-1}}')=\frac{\left(l_{i-1}+\sum_{j\neq i-1} d_j'\right)+
  \left(\alpha_{i-1}k_{i-1}+\sum_{j\neq i-1} \alpha_jr_j'\right)}
  {k_{i-1} +\sum_{j\neq i-1} r_j'}.
\]
From the short exact exact sequence
\[
0\lto N_{i-1}\lto E_{i-1}'\lto I_i\lto 0,
\]
we see that
\begin{align*}
  \rk(I_i\otimes K^*)   & =\rk (I_i) =r_{i-1}'-k_{i-1}, \\
  \deg(I_i\otimes K^*)  & =\deg (I_i) +\deg (K^*) \rk (I_i)\\
                        & = d_{i-1}'-l_{i-1}+(2-2g)(r_{i-1}'-k_{i-1}).
\end{align*}
Hence
\begin{align*}
  \mu_\alpha(C_{I_i}'') & =\frac{\deg (I_i\otimes K^*)
  +\alpha_i\rk (I_i\otimes K^*)}{\rk (I_i\otimes K^*)} \\
  & = \frac{d_{i-1}'-l_{i-1} - l_{i-1}}{r_{i-1}' - k_{i-1}} + 2-2g
  +\alpha_i.
\end{align*}
Using these formulae for $\mu_\alpha(C_{N_{i-1}}')$ and
$\mu_\alpha(C_{I_i}'')$, we obtain
\begin{multline}
\label{eq:weighted-slopes}
\left( k_{i-1} +\sum_{j\neq i-1} r_j'\right) \mu_\alpha(C_{N_{i-1}}') +
(r_{i-1}' - k_{i-1}) \mu_\alpha(C_{I_i}'')=\\
  \sum_{j=0}^n d_j' + (r_{i-1}' - k_{i-1}) (2-2g) +
  (\alpha_{i-1}-\alpha_i) k_{i-1} +
  \left( \sum_{j\neq i-1} \alpha_jr_j' + \alpha_i r_{i-1}' \right).
\end{multline}
To prove Part~i), suppose first that $i\in\{ 1, \ldots, n\}
\setminus \cD$, so $\alpha_i-\alpha_{i-1}> 2g-2$. Since $C'$ and
$C''$ are $\alpha$-semistable with the same $\alpha$-slope,
$\mu_\alpha(C_{N_{i-1}}') \leq \mu_\alpha(C')$,
$\mu_\alpha(C_{I_i}'') \leq \mu_\alpha(C'')$ and $\mu_\alpha(C') =
\mu_\alpha(C'')$, so \eqref{eq:weighted-slopes} is smaller than or
equal to
\begin{align*}
  \left( k_{i-1} +\sum_{j\neq i-1} r_j'\right) \mu_\alpha(C') +
  (r_{i-1}' - k_{i-1}) \mu_\alpha(C'') & =  \left(\sum_{j=0}^n
  r_j'\right) \mu_\alpha(C') \\ & =\sum_{j=0}^n d_j' +
  \sum_{j=0}^n \alpha_j r_j'.
\end{align*}
This is equivalent to the inequality
\[
  (\alpha_i-\alpha_{i-1})(r'_{i-1}-k_{i-1})\leq (r_{i-1}'-k_{i-1})
  (2g-2).
\]
If $\xi_i\neq 0$, then $N_i\neq 0$, so $r_{i-1}'-k_{i-1}>0$ and we
see that
\begin{equation}
\label{eq:inequality-alpha}
  \alpha_i-\alpha_{i-1}\leq 2g-2,
\end{equation}
which contracticts $\alpha_i-\alpha_{i-1}> 2g-2$. Therefore,
$\xi_i=0$.  Thus, $\ker (P_i)=0$, i.e., $P_i$ is injective for all
$i\in\{ 1, \ldots, n\} \setminus \cD$.  If $i\in\cD$, the fact that
$C''$ is $\alpha_i$-stable, where $\alpha_i:=\alpha+\epsilon_i u_i$,
implies the strict inequality $\mu_{\alpha_i}(C_{I_i}'') <
\mu_{\alpha_i}(C'')$ (note that $C_{I_i}''$ is a proper subchain of
$C''$, since $r_j''\neq 0$ for all $j$, so $C_{I_i}''\neq C''$, and
$\xi\neq 0$, so $C_{I_i}''\neq 0$). Defining $\lambda_i(C'') = r_i''/
\sum_{j=0}^n r_j''$, this strict inequality can be written as
\[
  \mu_\alpha(C_{I_i}'')-\mu_\alpha(C'') <
  \epsilon_i\lambda_i(C'')-\epsilon_i=
  \epsilon_i\left(\lambda_i(C'')-1\right),
\]
so $\mu_\alpha(C_{I_i}'')<\mu_\alpha(C'')$, as $\epsilon_i\geq 0$.
Hence, replacing the inequality $\mu_\alpha(C_{I_i}'') \leq
\mu_\alpha(C'')$ by the strict inequality, we obtain a strict
inequality in \eqref{eq:inequality-alpha}, which contradicts the
hypothesis $\alpha_i-\alpha_{i-1}\geq 2g-2$. Therefore, $\xi_i$ must
be zero, i.e., $P_i$ is also injective when $i\in\cD$.  Thus, we
conclude that $\HH^2(F^\bullet(C'',C'))=0$. This completes the proof
of Part~i).

Defining $\cD$ to be the empty set, we see that Part~i) implies
Part~ii), and defining $\cD=\{1,\ldots,n\}$ and $\epsilon_i= 0$ for
all $i=1,\ldots,n$, we see that Part~i) implies Part~iii).

To prove~iv), note that \eqref{eq:I-ker} implies that, if $\phi_i''$
is injective, then $I_i=0$, whereas \eqref{eq:im-N} implies that, if
$\phi_i'$ is generically surjective, then $N_{i-1}=E_{i-1}$. In both
cases, we deduce that $\xi_i=0$. Hence $\ker (P_i) =0$, i.e., $P_i$
is injective, for all $i=1,\ldots,n$. As explained above, this is
equivalent to $\HH^2(F^\bullet(C'',C'))=0$.
\end{proof}

\begin{remark}
Note that Part~iv) of Proposition~\ref{prop:h2-vanishing} generalizes
and gives a direct proof of \cite[Proposition~6.3]{BG} without
dimensional reduction techniques.
\end{remark}

The following is an immediate consequence of
Proposition~\ref{prop:h2-vanishing}. 

\begin{corollary}
\label{dimension-ext1}
 Suppose that $C',C''$ are $\alpha$-semistable chains such that
 $\mu_\alpha(C')=\mu_\alpha(C'')$, where $\alpha_i-\alpha_{i-1}>2g-2$
 for all $i=1,\ldots,n$. Then
\[
 \dim \Ext^1(C'',C')=h^0(C'',C') - \chi(C'',C').
\]
The same holds if the conditions of
\textnormal{Proposition~\ref{prop:h2-vanishing}, Part~i)}, are
satisfied for some subset $\cD\subset\{ 1,\ldots,n\}$ and
$\epsilon_i\geq 0$ ($i\in\cD$).
\end{corollary}

\subsection{Deformation theory of chains}
\label{sub:deformation-theory}

Let $\mathcal{M}^s_\alpha(t)$ be the moduli space of $\alpha$-stable
chains of type $t=(r_j,j=0,\ldots,n; d_j, j=0,\ldots,n)$.

\begin{theorem}
\label{thm:smoothdim}
 Let $C$ be an $\alpha$-stable chain of type $t$.
\par {\rm i)}
 The Zariski tangent space at the point defined by $C$ in the moduli
 space $\mathcal{M}^s_\alpha(t)$ is isomorphic to $\HH^1(F^\bullet(C,C))$.
\par {\rm ii)}
 If $\HH^2(F^\bullet(C,C))=0$, then the moduli space
 $\mathcal{M}^s_\alpha(t)$ is smooth in a neighbourhood of the point
 defined by $C$.
\par {\rm iii)}
 $\HH^2(F^\bullet(C,C))=0$ if and only if the homomorphism in the
 corresponding exact sequence of \textnormal{Proposition~\ref{prop:hypercohology}},
\[
  d\colon
  \bigoplus_{i=0}^n \Ext^1_X(E_i, E_i) \lto \bigoplus_{i=1}^n \Ext^1_X(E_i,E_{i-1}),
\]
 is surjective.
\par {\rm iv)}
 At a smooth point $C\in\mathcal{M}^s_\alpha(t)$, the dimension of
 the moduli space $\mathcal{M}^s_\alpha(t)$ is
 \begin{align*}
   \dim \mathcal{M}^s_\alpha(t) & = h^1(C,C)=1-\chi(C,C) \\
                              & =
   (g-1) \left( \sum_{i=0}^n r_i^2 - \sum_{i=1}^n r_i r_{i-1} \right)\\
   & \qquad  + \sum_{i=1}^n (r_i d_{i-1} - r_{i-1} d_i) +1.
 \end{align*}
\par {\rm v)}
 If for each $i=1,\ldots,n$, $\phi_i\colon E_i\to E_{i-1}$ is
 injective or generically surjective, then $\HH^2(F^\bullet(C,C))=0$
 and therefore $C$ defines a smooth point of
 $\mathcal{M}^s_\alpha(t)$.
\par {\rm vi)}
  If $\alpha_i-\alpha_{i-1}\geq 2g-2$ for all $i=1,\ldots,n$, then
  $C$ defines a smooth point in the moduli space, and hence
  $\mathcal{M}^s_\alpha(t)$ is smooth.
\end{theorem}

\begin{proof}
 All the results of this theorem follow immediately
 from the results of this section, as in the proof of
 \cite[Theorem~3.8]{BGG}, noting that the space of infinitesimal
 deformations of $C$ is $\HH^1(F^\bullet(C,C))$.
\end{proof}

Part iv) of Proposition~\ref{thm:smoothdim} highlights the importance
of the following.

\begin{definition}
 \label{def:R_g}
  The region $R_{2g-2}\subset\R^{n+1}$ is the set of points $\alpha$
  such that $\alpha_i-\alpha_{i-1}\geq 2g-2$ for all $i=1,\ldots,n$.
\end{definition}

\section{Wall-crossing}        
\label{sec:wall-crossing}

\subsection{Flip loci}
\label{sub:flip-loci}

In this section we study the variations in the moduli spaces
$\mathcal{M}^s_\alpha(t)$, for fixed type $t$ and different values of
$\alpha$. We begin with a set-theoretic description of the
differences between two spaces $\mathcal{M}^s_{\alpha}(t)$ and
$\mathcal{M}^s_{\beta}(t)$ when $\alpha$ and $\beta$ are separated by
a wall $w$ (as defined in Section \ref{sub:chambers-walls}). For the
rest of this section we adopt the following notation: Let $\alpha_w$
be the point in the parameter space obtained by intersecting the line
determined by $\alpha$ and $\beta$ with the wall $w$.  Set
\begin{equation}\nonumber
\alpha_w^+ = \alpha_w + \epsilon,\quad
\alpha_w^- = \alpha_w - \epsilon,
\end{equation}
\noindent where $|\epsilon| > 0$ is small enough
so that $\alpha_w^-$ is in the same chamber  as $\alpha$ and
$\alpha_w^+$ is in the same chamber  as $\beta$.

We define {\it flip loci}
$\mathcal{S}_{\alpha_w^{\pm}}\subset\mathcal{M}^s_{\alpha_w^{\pm}}$
by the conditions that the points in $\mathcal{S}_{\alpha_w^+}$ represent
chains which are $\alpha_w^+$-stable  but $\alpha_w^-$-unstable,
while the points in $\mathcal{S}_{\alpha_w^-}$ represent chains
which are $\alpha_w^-$-stable  but $\alpha_w^+$-unstable.
The following is immediate.

\begin{lemma}\label{lemma:fliploci}In the above notation:
\begin{equation}\nonumber
\mathcal{M}^s_{\alpha_w^+}-\mathcal{S}_{\alpha_w^+}=
\mathcal{M}^s_{\alpha_w}=
\mathcal{M}^s_{\alpha_w^-}-\mathcal{S}_{\alpha_w^-}.
\end{equation}
\end{lemma}

\begin{remark} If the wall $w$ is included in one of the
bounding hyperplanes we have that either
$\mathcal{S}_{\alpha_w^+}=\mathcal{M}^s_{\alpha_w^+}$, or
$\mathcal{S}_{\alpha_w^-}=
\mathcal{M}^s_{\alpha_w^-}$. The only interesting cases are thus those
for which  $w$ is not a bounding wall.
\end{remark}

A chain $C\in \mathcal{S}_{\alpha_w^\pm}$ must be strictly
$\alpha_w$-semistable. Hence in order to compute the codimension of
$\mathcal{S}_{\alpha_w^\pm}$ in $\mathcal{M}^s_{\alpha_w^\pm}$ we
have to ``count'' Jordan-H\"older filtrations of $C$. Let $C$ be a
strictly $\alpha_w$-semistable chain. As we have seen in Section
\ref{sec:definitions}, there is a filtration of chains given by
\begin{equation}\label{j-h-filtration}
0=C_0\subset C_1\subset \cdots  \subset C_m=C,
\end{equation}
with $G_i=C_i/C_{i-1}$ $\alpha_w$-stable and
$\mu_{\alpha_w}(G_i)=\mu_{\alpha_w}(C)$ for $1\leq i \leq m$.
Let $\gr(C)=\bigoplus_i G_i$ be the graduation of $C$.

\begin{proposition} Let $w$ be a wall contained in the region
$R_{2g-2}$.  Let ${\mathcal{S}}$ be a family of $\alpha_w$-semistable
chains $C$ of type $t$, all of them pairwise non-isomorphic, and
whose Jordan-H\"older filtration (\ref{j-h-filtration}) has
graduation $\gr(C)=\oplus_{i=1}^m G_i$, with $G_i$ of type $t_i$.
Then
\begin{equation}\label{dim-family}
\dim {\mathcal{S}}\leq
- \sum_{i\leq j}\chi(t_j,t_i)  -\frac
{m(m-3)}{2}.
\end{equation}
\end{proposition}
\begin{proof}

It is clear that
$$
\dim \mathcal{S}\leq \sum_{1\leq i \leq m} \dim \mathcal{M}^s_{\alpha_w}(t_i) +
\sum_{1\leq i<j\leq m} \dim \PP(\Ext^1(G_j,G_i)),
$$
since $\PP(\Ext^1(G_j,G_i))$ parametrizes equivalence classes of extensions
$$
0\lra G_i\lra G\lra G_j \lra 0.
$$
By Theorem \ref{thm:smoothdim}, $\mathcal{M}^s_{\alpha_w}(t_i)$ is
smooth and $ \dim \mathcal{M}^s_{\alpha_w}(t_i) = 1-\chi(t_i,t_i)$.
From Corollary \ref{dimension-ext1}, we have that the dimension of
$\Ext^1(G_j,G_i)$ is given by $h^1(G_j, G_i)= -\chi(t_j,t_i)$.  Here
we are using the vanishing of $h^2(G_j,G_i)$ given by Proposition
\ref{prop:h2-vanishing}, and we have assumed that $G_i$ and $G_j$ are
not isomorphic, and hence by Proposition \ref{prop:h0-vanishing}
$h^0(G_j,G_i)=0$, since otherwise we would have a subfamily of
positive codimension in ${\mathcal{S}}$). The result follows now by
adding up these dimensions.
\end{proof}

In order to show that the flip loci
$\mathcal{S}_{\alpha_w^{\pm}}\subset\mathcal{M}^s_{\alpha_w^{\pm}}$
have positive codimension we need to bound the values of
$\chi(t_j,t_i)$ in (\ref{dim-family}). This is what we do next.

\subsection{Bounds for $\chi$}
\label{sub:chi-bounds}

\begin{proposition}
\label{prop:stability-for-H_*}
Let $C', C''$ be two holomorphic chains, as in \eqref{eq:C'},
\eqref{eq:C''}, and $b$ the differential of the complex
$F^\bullet(C'',C')$, as in \eqref{eq:extension-complex}. If $C'$ and
$C''$ are $\alpha$-polystable and $\alpha_i-\alpha_{i-1}\geq 2g-2$
for all $i=1,\ldots,n$, then
 \begin{align}
   \label{eq:ker-b}
   \mu(\ker (b)) & \leq\mu_\alpha(C')-\mu_\alpha(C''),\\
   \label{eq:coker-b}
   \mu(\coker (b)) & \geq \mu_\alpha(C')-\mu_\alpha(C'') + 2g-2.
 \end{align}
\end{proposition}

\begin{proof}
We start constructing a holomorphic chain
\[
  \widetilde{C}(C'',C')\colon F^{-1} \overset{c}{\lto} F^0
  \overset{b}{\lto} F^1 \overset{a}{\lto} F^2,
\]
where $F^0, F^1$ and $b$ are defined as in \eqref{eq:extension-complex},
\[
  F^{-1}=\bigoplus_{i=1}^n \Hom(E_{i-1}'',E_i'), \qquad
  F^2=\bigoplus_{i=2}^n \Hom(E_i'',E_{i-2}'),
\]
and
\[
   a(\theta_1,\ldots,\theta_n)=\sum_{i=2}^n a_i(\theta_{i-1},\theta_i),
  \qquad
   c(\omega_1,\ldots,\omega_n)=\sum_{i=0}^n c_i(\omega_i,\omega_{i+1}),
\]
where $\theta_i\in\Hom(E_i'',E_{i-1}')$,
$\omega_i\in\Hom(E_{i-1}'',E_i')$, and
\begin{align*}
 a_i\colon & \Hom(E_{i-1}'',E_{i-2}') \oplus \Hom(E_i'',E_{i-1}')\to
 \Hom(E_i'',E_{i-2}') \hra F^2, \\
  c_i\colon & \Hom(E_{i-1}'',E_i')\oplus
\Hom(E_i'',E_{i+1}')\to \Hom(E_i'',E_i')\hra F^0,
\end{align*}
are given by
\[
  a_i(\theta_{i-1},\theta_i)=\theta_{i-1}\circ\phi_i''-\phi_{i-1}'\circ\theta_i,
  \quad
  c_i(\omega_i,\omega_{i+1})=\omega_i\circ\phi_i''-\phi_{i+1}'\circ\omega_{i+1},
\]
with $E_{n+1}'=0=E_{-1}''$ by convention. Note that the holomorphic
chain $\widetilde{C}(C'',C')$ is not in general a complex.

Suppose that $C',C''$ are $\alpha$-polystable. Then by
Theorem~\ref{thm:chain-hitchin-kobayashi}, there are Hermitian
metrics on the vector bundles $E'_i$ and $E''_i$ satisfying the
$\tau'$- and $\tau''$-vortex equations
\begin{equation}
\label{eq:chain-vortex}
\begin{split}
  \sqrt{-1}\Lambda F(E'_{i}) + \phi'_{i+1} \phi_{i+1}^{\prime *} -
  \phi_{i}^{\prime *}\phi_{i}
  &= \tau'_{i} {\Id}_{E'_{i}}, \\
  \sqrt{-1}\Lambda F(E''_{i}) + \phi''_{i+1} \phi_{i+1}^{\prime\prime *} - \phi_{i}^{\prime\prime *}\phi''_{i}
  &= \tau''_{i} {\Id}_{E''_{i}},
\end{split}
\end{equation}
for $i=0,\ldots,n$, where $\tau', \tau''\in\R^{n+1}$ are given by
\begin{equation}
\label{eq:tau'-tau''}
\tau'_i=\mu_\alpha(C') - \alpha_i
\qquad \text{and} \qquad
\tau''_i = \mu_\alpha(C'') - \alpha_i.
\end{equation}
Using these equations, we now show that the induced metrics on the
bundles $F^i$, for $i=-1,0,1,2$, which are the terms of the
holomorphic chain $\widetilde{C}(C',C'')$, satisfy the equations
\begin{equation}
\label{eq:chain-vortex-F}
\begin{split}
  \sqrt{-1}\Lambda F(F^0) + c\circ c^* - b^*\circ b
  = & \left(\mu_\alpha(C') - \mu_\alpha(C'')\right){\Id}, \\
  \sqrt{-1}\Lambda F(F^1) + b\circ b^* - a^*\circ a
  = & \left(\mu_\alpha(C') - \mu_\alpha(C'')\right){\Id}  \\
    & + \sum_{i=1}^n (\alpha_i-\alpha_{i-1})\pi^1_{i},
\end{split}
\end{equation}
where $\pi^1_{i}\colon F^1\to \Hom(E_i'',E_{i-1}')$ is the
canonical projection. To prove this, let $\psi_i \in
\Hom(E_i'',E_i')$ and $\zeta_i \in \Hom (E_i'',E_{i-1}')$. The
curvature $F(F^p)$ of the induced connection on $F^p$, for $p=0,1$,
is the $\left(\End (F^p)\right)$-valued $2$-form given by
\begin{align*}
  F(F^0)(\psi_i)  & = F(E_i')\wedge\psi_i - \psi_i\wedge F(E_i''), \\
  F(F^1)(\zeta_i) & = F(E_{i-1}')\wedge\zeta_i - \zeta_i\wedge F(E_i''),
\end{align*}
so the first terms in the left-hand sides of
\eqref{eq:chain-vortex-F} are given by
\begin{align}\label{eq:chain-vortex-F0-1}
  \sqrt{-1}\Lambda F(F^0)(\psi_i) &=\sqrt{-1}\Lambda F(E_i')\circ\psi_i - \psi_i\circ \sqrt{-1}\Lambda F(E_i''), \\
\label{eq:chain-vortex-F1-1}
  \sqrt{-1}\Lambda F(F^1)(\zeta_i) &=\sqrt{-1}\Lambda
  F(E_{i-1}')\circ\zeta_i - \zeta_i\circ \sqrt{-1}\Lambda F(E_i'').
\end{align}
The remaining terms in the left-hand side of the first equation in
\eqref{eq:chain-vortex-F} are
\begin{align}
\begin{split}\label{eq:chain-vortex-F0-2}
 c\circ c^*(\psi_i) = & \psi_i\circ  \phi_{i}^{\prime\prime *}\phi''_{i} +
 \phi'_{i+1} \phi_{i+1}^{\prime *} \circ \psi_i \\
 & - \phi^{\prime *}_{i+1} \circ \psi_i \circ \phi''_{i+1} - \phi'_{i} \circ \psi_i
 \circ \phi^{\prime\prime *}_{i}, \\
 b^*\circ b(\psi_i) = & \psi_i\circ \phi''_{i+1} \phi_{i+1}^{\prime\prime *} +
 \phi_{i}^{\prime *} \phi'_{i} \circ \psi_i \\
 & - \phi'_{i} \circ \psi_i \circ \phi^{\prime\prime *}_{i} - \phi^{\prime *}_{i+1}
\circ \psi_i \circ \phi''_{i+1},
\end{split}
\end{align}
whereas the remaining terms in the left-hand side of the second equation in
\eqref{eq:chain-vortex-F} are
\begin{align}
\begin{split}\label{eq:chain-vortex-F1-2}
 b\circ b^*(\zeta_i) = & \zeta_i\circ  \phi_{i}^{\prime\prime *}\phi''_{i} +
 \phi'_{i} \phi_{i}^{\prime *} \circ \zeta_i \\
 & - \phi^{\prime *}_{i} \circ \zeta_i \circ \phi''_{i+1} - \phi'_{i-1} \circ \zeta_i
 \circ \phi^{\prime\prime *}_{i}, \\
 a^*\circ a(\zeta_i) = & \zeta_i\circ \phi''_{i+1} \phi_{i+1}^{\prime\prime *} +
 \phi_{i-1}^{\prime *} \phi'_{i-1} \circ \zeta_i \\
 & - \phi'_{i-1} \circ \zeta_i \circ \phi^{\prime\prime *}_{i} -
 \phi^{\prime *}_{i} \circ \zeta_i \circ \phi''_{i+1}.
\end{split}
\end{align}
Using \eqref{eq:chain-vortex-F0-1} and \eqref{eq:chain-vortex-F0-2},
and \eqref{eq:chain-vortex-F1-1} and \eqref{eq:chain-vortex-F1-2},
together with \eqref{eq:chain-vortex}, it follows immediately that
the left-hand sides of \eqref{eq:chain-vortex-F} are
\begin{align}
\label{eq:chain-vortex-F0-3}
 \sqrt{-1}\Lambda F(F^0) + c\circ c^* - b^*\circ b &
 = \sum_{i=0}^n (\tau_i' - \tau_i'')\pi_i^0,
\\
\label{eq:chain-vortex-F1-3}
 \sqrt{-1}\Lambda F(F^1) + b\circ b^* - a^*\circ a &
 = \sum_{i=1}^n (\tau_{i-1}' - \tau_i'')\pi_i^1,
\end{align}
where $\pi^0_{i}\colon F^0\to \Hom(E_i'',E_i')$ is the canonical
projection. Now, it follows from \eqref{eq:tau'-tau''} that the
right-hand sides of \eqref{eq:chain-vortex-F0-3} and
\eqref{eq:chain-vortex-F1-3} equal the right-hand sides of the first
and the second equation in \eqref{eq:chain-vortex-F}, respectively,
so equations \eqref{eq:chain-vortex-F} are satisfied.

We can now use the equations \eqref{eq:chain-vortex-F} to obtain the
inequalities \eqref{eq:ker-b} and \eqref{eq:coker-b}. Let $G\subset
\coker (b)$ be the maximal vector subbundle of $\coker (b)$. Note
that
\begin{equation}
\label{eq:rk-deg-vbcok}
 \rk(\coker (b)) =\rk (G), \qquad \deg(\coker (b))\geq \deg (G),
\end{equation}
and, by standard results (see e.g. \cite{UY,S,B}),
\begin{align}
\label{eq:deg-ker}
 \deg(\ker (b)) & =\frac{1}{2\pi}\left(\int_X\tr\left(\pi_0\sqrt{-1}\Lambda F(F^0)\right) - \|\beta_0\|^2_{L^2}\right), \\
\label{eq:deg-coker}
 \deg (G) & =\frac{1}{2\pi}\left(\int_X\tr\left(\pi_{1}\sqrt{-1}\Lambda F(F^1))\right) + \|\beta_{1}\|^2_{L^2}\right),
\end{align}
where
\[
 \pi_0\colon F_0\lto F_0
 \qquad \text{and} \qquad \pi_1\colon F^1\lto F^1
\]
are the orthogonal projection operators onto $\ker (b)$ and $G$, with
respect to the induced Hermitian metrics on $F^0$ and $F^{1}$ (using
the metric on $F^1$, $G$ is regarded as a (smooth) subbundle of
$F^1$), and $\beta_0$ and $\beta_1$ are the corresponding second
fundamental forms, i.e., the
$\left(\End\left(F^0\right)\right)$-valued and
$\left(\End\left(F^1\right)\right)$-valued $(0,1)$-forms
\[
 \beta_0=\dbar_{F^0}(\pi_0), \qquad
 \beta_{1}=\dbar_{F^1}(\pi_1).
\]
Let $a_\perp=a\circ\pi_1\colon F^1\to F^2$ and
$c_\perp=\pi_0\circ c\colon F^{-1}\to F^0$. Since $a_\perp^*=\pi_1
\circ a^*$ and $c_\perp^*=c^*\circ \pi_0$,
\begin{equation}
\label{eq:tr-a-b-perp}
\begin{split}
  & \tr(\pi_1\circ a^*a) = \tr((\pi_1 a^*)(a\pi_1))=\tr(a_\perp^* a_\perp)=|
  a_\perp |^2, \\
  & \tr(cc^*\circ \pi_0)  = \tr((\pi_0 c)(c^*\pi_0))=\tr(c_\perp c_\perp^*)=|
  c_\perp |^2,
\end{split}
\end{equation}
(where $|\cdot|$ are the induced norms).
By the definitions of $\pi_0$ and $\pi_1$,
\begin{equation}
\label{eq:b-c-vanishing}
b\circ\pi_0=0 \qquad\text{and}\qquad \pi_1\circ b=0.
\end{equation}
Applying $\tr(-\circ\pi_0)$ and $\tr(\pi_1\circ-)$ to
\eqref{eq:chain-vortex-F}, and using
\eqref{eq:tr-a-b-perp} and \eqref{eq:b-c-vanishing}, we obtain
\begin{equation}
\label{eq:ineq-ker-b-1}
\tr\left(\sqrt{-1}\Lambda F(F^0)\pi_0\right) + |c_\perp|^2 =
\left(\mu_\alpha(C') - \mu_\alpha(C'')\right)\rk(\ker (b)),
\end{equation}
\begin{equation}
\begin{split}
\label{eq:ineq-coker-b-1}
\tr\left(\sqrt{-1}\Lambda F(F^1)\pi_1\right) - |a_\perp|^2 =
      \left(\mu_\alpha(C') - \mu_\alpha(C'')\right)\rk(\coker (b))\\
      + \sum_{i=1}^n (\alpha_i-\alpha_{i-1})\tr(\pi^1_{i}\pi_1),
\end{split}
\end{equation}
respectively (where $\tr(\pi_1)=\rk(\coker (b))$ by
\eqref{eq:rk-deg-vbcok}). Since $\tr(\pi^1_{i}\pi_1)\geq 0$,
$\alpha_i-\alpha_{i-1}\geq 2g-2$ for all $i$, and $\sum_{i=1}^n
\pi^1_{i}= {\Id}_{F^1}$, the last term in the right-hand side of
\eqref{eq:ineq-coker-b-1} satisfies
\[
  \sum_{i=1}^n (\alpha_i-\alpha_{i-1})\tr(\pi^1_{i}\pi_1)\geq
 (2g-2) \sum_{i=1}^n \tr(\pi^1_{i}\pi_1) = (2g-2) \rk(\coker (b)),
\]
so \eqref{eq:ineq-coker-b-1} implies
\begin{multline}
\label{eq:ineq-coker-b-2}
\tr\left(\sqrt{-1}\Lambda F(F^1)\pi_1\right) - |a_\perp|^2 \\
\geq  \left(\mu_\alpha(C') - \mu_\alpha(C'') + 2g-2\right)\rk(\coker (b)).
\end{multline}
Integrating \eqref{eq:ineq-ker-b-1} and \eqref{eq:ineq-coker-b-2}
over $X$, using \eqref{eq:rk-deg-vbcok}, \eqref{eq:deg-ker} and
\eqref{eq:deg-coker}, and dividing by $\vol(X)=2\pi$, we obtain
\[
 \deg(\ker (c))+\frac{1}{2\pi}\left(\|\beta_0\|^2_{L^2} + \|c_\perp\|^2_{L^2} \right)
   = \left(\mu_\alpha(C')-\mu_\alpha(C'')\right)\rk(\ker (b)),
\]
\begin{multline*}
 \deg(\coker (c))-\frac{1}{2\pi}\left(\|\beta_1\|^2_{L^2} + \|a_\perp\|^2_{L^2} \right)
\\
\geq \left(\mu_\alpha(C')-\mu_\alpha(C'') + 2g-2\right)\rk(\coker (b)),
\end{multline*}
respectively, which imply \eqref{eq:ker-b} and \eqref{eq:coker-b}.
\end{proof}

\begin{proposition}
\label{prop:chi-bounds}
 Let $C'$ and $C''$ be non-zero holomorphic chains of types $t'$ and
 $t''$, respectively, and let $\alpha\in \R^{n+1}$. Suppose that the
 following conditions hold:
 \begin{itemize}
   \item $C'$ and $C''$ are $\alpha$-polystable with $\mu_\alpha(C')=\mu_\alpha(C'')$,
   \item $\alpha_i-\alpha_{i-1}\geq 2g-2$ for all $i=1,\ldots,n$,
   \item the map $b\colon F^0\lto F^1$ of
   \eqref{eq:extension-complex} is not generically an isomorphism.
 \end{itemize}
 Then $\chi(C'',C')\leq 1-g$.
 In particular, if $g\geq 2$, then $\chi(C'',C')<0$.
\end{proposition}

\begin{proof}
 Let $F^\bullet(C'',C')\colon F^0 \overset{b}{\lto} F^1$ be the
 complex~\eqref{eq:extension-complex}. By Proposition~\ref{prop:chi(C'',C')},
\[
 \chi(C'',C')= (1-g) \left(\rk (F^0) - \rk (F^1)\right) + \deg (F^0) - \deg (F^1).
\]
 Using
 \begin{align*}
 & \deg (F^0)=\deg(\ker (b)) + \deg(\im (b)), \\  
 & \deg (F^1)  =\deg(\im (b))  + \deg(\coker (b)), \\
 & \rk (F^1) =\rk(\im (b))   + \rk(\coker (b)),
 \end{align*}
 and the inequalities \eqref{eq:ker-b} and \eqref{eq:coker-b} with
 $\mu_\alpha(C')=\mu_\alpha(C'')$, we see that
 \begin{align*}
  \deg (F^0) -\deg (F^1) & \leq 2(1-g) \rk (\coker (b))\\
  & =2(1-g) \left(\rk (F^1) -\rk (\im (b))\right),
 \end{align*}
so
\begin{equation}
\label{eq:chi-first-bound}
 \chi(C'',C')\leq (1-g) \left(\rk (F^0) + \rk (F^1) - 2 \rk(\im (b))\right).
\end{equation}
Note that
\begin{equation}
\label{eq:ineq-rk-im-b}
  \rk (F^0) + \rk (F^1) - 2 \rk (\im (b))\geq 0,
\end{equation}
with equality in \eqref{eq:ineq-rk-im-b} if and only if $\rk (F^0) =
\rk (\im (b))=\rk (F^1)$. But $b$ is not generically an isomorphism,
so the equality in \eqref{eq:ineq-rk-im-b} does not hold, i.e., $\rk
(F^0) + \rk (F^1) - 2 \rk (\im (b)) \geq 1$.  Therefore,
\eqref{eq:chi-first-bound} implies $\chi(C'',C')\leq 1-g$.
\end{proof}

\subsection{The birationality region}
\label{sub:birationality-region}

Proposition~\ref{prop:chi-bounds} motivates the definition of a
region $\wt{R}(t)\subset\R^{n+1}$ as follows. First, we recall
from~\secref{sub:parameter-region} that a \emph{linear chain} $V$
(over $\CC$) is a chain in the category of complex vector spaces,
i.e., a diagram of complex vector spaces $V_i$ and linear maps $f_i$
which compose as follows.
\begin{equation}
\label{eq:linear-chain}
  V \colon V_n \overset{f_n}{\lto} V_{n-1}
  \overset{f_{n-1}}{\lto} \cdots
  \overset{f_{1}}{\lto} V_{0} .
\end{equation}
Note that this is simply a holomorphic chain when $X$ is a point. The
\emph{dimension vector} of $V$ is the $(n+1)$-tuple of integers
$r=(r_j,j=0,\ldots ,n)$, with $r_j=\dim V_j$. Given two linear chains
$V', V''$, we define a 2-step complex of vector spaces over $\CC$,
\begin{equation}
\label{eq:extension-vector-complex}
  F^\bullet(V'',V')\colon F^0 \overset{b}{\lto} F^1,
\end{equation}
exactly as in~\eqref{eq:extension-complex}, where now $X$ is a point,
so $b$ is simply a linear map.
Note also that the $\alpha$-slope $\mu_\alpha(C)$ of a holomorphic
chain $C$ only depends on the type $t=(r_j,j=0,\ldots,n; d_j,
j=0,\ldots,n)$ of $C$, so we may use the notation
\[
  \mu_\alpha(t) := \frac{\sum_{i=0}^n (d_i+\alpha_i r_i)}{\sum_{i=0}^n r_i}.
\]

\begin{definition}
\label{def:Rtilde}
  Fix a type $t=(r_j,j=0,\ldots,n; d_j, j=0,\ldots,n)$.
  The region $\wt{R}(t)\subset\R^{n+1}$ is the set of points $\alpha$
  such that for all types $t'=(r_j',j=0,\ldots,n; d_j', j=0,\ldots,n)$ and
  $t''=(r_j'',j=0,\ldots,n; d_j'', j=0,\ldots,n)$, with $t'+t''=t$ and
  $\mu_\alpha(t')=\mu_\alpha(t'')$, and for all linear chains $V'$
  and $V''$ with dimension vectors $r'=(r_j',j=0,\ldots,n)$ and
  $r''=(r_j'',j=0,\ldots,n)$, respectively, the map $b$
  of~\eqref{eq:extension-vector-complex} is not an isomorphism.
\end{definition}

Note that the region $\wt{R}(t)$ is generally non-connected.
Furthermore, the definition of $\wt{R}(t)$ does not involve the
geometry of $X$ but only linear algebra.

The following is an immediate consequence of Proposition~\ref{prop:chi-bounds}.

\begin{theorem}
\label{thm:chi-bounds}
 Let $C'$ and $C''$ be non-zero holomorphic chains of types $t'$ and
 $t''$, respectively, and let $\alpha\in \R^{n+1}$. Suppose that the
 following conditions hold:
 \begin{itemize}
   \item $C'$ and $C''$ are $\alpha$-polystable with $\mu_\alpha(C')=\mu_\alpha(C'')$,
   \item $\alpha\in \wt{R}(t)\cap R_{2g-2}$.
 \end{itemize}
 Then $\chi(C'',C')\leq 1-g$.
 In particular, if $g\geq 2$, then $\chi(C'',C')<0$.
\end{theorem}

Because of this theorem, it becomes an important problem to
characterize the \emph{birationality region} $\wt{R}(t)$ of
Definition~\ref{def:Rtilde}.  
The rest of this subsection is devoted to the determination of
$\wt{R}(t)$. 

\begin{definition}
\label{def:region-V}
\begin{enumerate}
  \item[(1)] 
  Fix a dimension vector $r=(r_j, j=0, \ldots, n)$. Let $\cV(r)$ be
  the set of pairs $(r',r'')$ of dimension vectors
  $r'=(r_j',j=0,\ldots,n)$ and $r''=(r_j'',j=0,\ldots,n)$ with
  $r'+r''=r$, such that there exist linear chains $V'$ and $V''$ of
  dimension vectors $r'$ and $r''$, respectively, for which the
  differential $b$ of the
  complex~\eqref{eq:extension-vector-complex}, corresponding to $V'$
  and $V''$, is an isomorphism.
\item[(2)] 
  Given a type $t=(r,d)$, let $\cT(t)$ be the set of pairs $(t',t'')$
  of types $t'=(r_j',j=0,\ldots,n; d_j', j=0,\ldots,n)$ and 
  $t''=(r_j'',j=0,\ldots,n; d_j'', j=0,\ldots,n)$, with $t'+t''=t$, such
  that $(r',r'')\in\cV(r)$, where $r'=(r_j',j=0,\ldots,n)$ and
  $r''=(r_j'',j=0,\ldots,n)$.
\item[(3)] 
  Given $(t',t'')\in\cT(t)$, let $\cB(t', t'')\subset \R^{n+1}$ be
  the set of points $\alpha$ on the hyperplane
\[
  \mu_\alpha(t') =  \mu_\alpha(t'').
\]
\end{enumerate}
\end{definition}

\begin{proposition}
\label{prop:Rtilda}
  $\wt{R}(t)=\R^{n+1}\setminus \cB(t)$, where the `boundary' is
\[
  \cB(t) =  \bigcup_{(t,',t´)\in \cT(t)} \cB(t', t'').
\]
\end{proposition}

\begin{proof}
  This follows automatically from Definitions~\ref{def:Rtilde}
  and~\ref{def:region-V}. 
\end{proof}

Thus, to determine the region $\wt{R}(t)$ we need to find the set
$\cV(r)$. This can be done by using the following results. Note first
that the Euler characteristic of two linear chains $V'$ and $V''$ is
\[
  \chi(V'',V'):= \dim\Hom(V'',V') - \dim\Ext^1(V'',V').
\]
The following lemma can be compared with Propositions~\ref{prop:hypercohology}
and~\ref{prop:chi(C'',C')}. 

\begin{lemma}
\label{lem:cohomology-linear-chains}
Let $V'$ and $V''$ be two linear chains with dimension vectors
$r'=(r'_j,j=0,\ldots,n)$ and $r''=(r''_j,j=0,\ldots,n)$, respectively. Let
$b$ be the map in~\eqref{eq:extension-vector-complex}. Then there is
a canonical exact sequence
\begin{equation}
\label{eq:canonical-sequence}
  0 \to \Hom(V'',V') \lto F^0 \overset{b}{\lto} F^1
  \lto \Ext^1(V'',V') \lto 0.
\end{equation}
Hence, the Euler characteristic is given by
\begin{equation}
\label{eq:chi(V'',V')}
  \chi(V'',V')=\dim(F^0) - \dim(F^1) = \sum_{i=0}^n r''_i r'_i -
  \sum_{i=1}^n r''_i r'_{i-1}.
\end{equation}
\end{lemma}

\begin{proof}
The exact sequence~\eqref{eq:canonical-sequence} is obtained as in
Proposition~\ref{prop:hypercohology} when $X$ is point, whereas 
equation~\eqref{eq:chi(V'',V')} follows immediately from this.
\end{proof}

By Lemma~\ref{lem:cohomology-linear-chains}, $\chi(V'',V')$ only
depends on the dimension vectors $r'$ and $r''$ of $V'$ and $V''$,
respectively, so we may use the notation
\[
  \chi(r'',r'):=\chi(V'',V').
\]

\begin{proposition}
\label{prop:isom-indec}
 Let $V', V''$ be linear chains, given by 
\[
  V'=\bigoplus_{j=1}^{s'} V'_j, 
  \qquad
  V''=\bigoplus_{k=1}^{s''} V''_k, 
\]
 where $V'_j, V''_k$ are linear chains. Let
\[
  F^\bullet(V'',V')\colon F^0 \overset{b}{\lto} F^1
\]
be the 2-step complex corresponding to $V'$ and $V''$, as
in~\eqref{eq:extension-vector-complex}, and
\[
  F^\bullet(V''_k,V'_j)\colon F^0_{jk} \overset{b_{jk}}{\lto}
  F^1_{jk},  \qquad \text{ for $1\leq j\leq s'$, $1\leq k\leq s''$}, 
\]
the 2-step complexes corresponding to the subchains $V''_k$ and
$V'_j$. The following conditions are equivalent:
\par {\rm i)}
The map $b$ is an isomorphism.
\par {\rm ii)}
$\Hom(V'',V')=0=\Ext^1(V'',V')$. 
\par {\rm iii)}
The maps $b_{jk}$ are isomorphisms, for all $1\leq j\leq s'$, $1\leq k\leq s''$. 
\par {\rm iv)} $\Hom(V''_k,V'_j)=0=\Ext^1(V''_k,V'_j)$, for all $1\leq j\leq
  s'$, $1\leq k\leq s''$.
\end{proposition}

\begin{proof}
The equivalence ii)$\Leftrightarrow$iv) comes from the obvious
isomorphisms
\[
  \Hom(V'',V')=\!\!\bigoplus_{\substack{1\leq j\leq s'\\1\leq k \leq s''}}\!\!\Hom(V''_k,V'_j),
  \,\, 
  \Ext^1(V'',V')=\!\!\bigoplus_{\substack{1\leq j\leq s'\\1\leq k \leq s''}}\!\!\Ext^1(V''_k,V'_j), 
\]
whereas i)$\Leftrightarrow$ii) and iii)$\Leftrightarrow$iv) follow
from Lemma~\ref{lem:cohomology-linear-chains}.
\end{proof}

\begin{definition}
\label{def:indec-rep}
Given integers $0\leq p\leq q\leq n$, the linear chain
$\delta_{[p,q]}$ is defined by the diagram
\begin{equation*}
\xymatrix@C=12pt@R=0pt{
\delta_{[p,q]}\colon & 0 \ar[r] & \cdots \ar[r] &
0\ar[r] & \CC \ar@{=}[r] & \CC \ar@{=}[r] & \cdots \ar@{=}[r] & \CC
\ar@{=}[r] & \CC \ar[r] & 0 \ar[r] & \cdots \ar[r] & 0 \\
& & & & {\scriptstyle q}  & & & & {\scriptstyle p} & & &   
}\end{equation*}
Thus, $\delta_{[p,q]}$ is given by~\eqref{eq:linear-chain},
where $V_i=\CC$ if $p\leq i\leq q$ and $V_i=0$ otherwise, whereas
$f_i=\id\colon \CC\to\CC$ if $p<i\leq q$ and $f_i=0$ otherwise.
The dimension vector of $\delta_{[p,q]}$ is denoted by $r_{[p,q]}$. 
\end{definition}

\begin{proposition}
\label{prop:indec-rep}
\par {\rm i)}
All linear chains are direct sums of indecomposable ones.
\par {\rm ii)}
The linear chains $\delta_{[p,q]}$ are indecomposable and any
indecomposable linear chain is isomorphic to such a $\delta_{[p,q]}$.
\end{proposition}

\begin{proof}
Part i) is well-known. Part ii) can be found, e.g., in~\cite{Ga}.
\end{proof}

Thus, the linear chains $V'$ and $V''$ can be written as 
\begin{equation}
\label{eq:decomposition-rep}
  V' \cong \bigoplus_{j=1}^{s'} \delta_{[p'_j,q'_j]}, 
  \qquad
  V''\cong \bigoplus_{k=1}^{s''} \delta_{[p''_k,q''_k]}, 
\end{equation}
for sets of pairs of integers $(p'_j,q'_j)$ and $(p''_k,q''_k)$, for
$1\leq j\leq s'$ and $1\leq k\leq s''$, satisfying $0\leq p'_j\leq
q'_j\leq n$ and $0\leq p''_k\leq q''_k\leq n$, for all $j$ and $k$.
Now, Proposition~\ref{prop:Rtilda} reduces the problem of finding
$\wt{R}(t)$ to the problem of finding $\cT(t)$ or equivalently (by
Definition~\ref{def:region-V}), $\cV(r)$, whereas
Propositions~\ref{prop:isom-indec} and~\ref{prop:indec-rep} reduce
this problem to finding when there are no homomorphisms and
extensions between the indecomposable linear chains
$\delta_{[p',q']}$ and $\delta_{[p'',q'']}$. This is in
Proposition~\ref{prop:coho-indec-rep} below.

Given two integers $p$ and $q$, let $[p,q]=\{ p, p+1, \ldots, q-1,
q\}$ if $p\leq q$, and $[p,q]=\varnothing$ otherwise. Note that,
by~\eqref{eq:chi(V'',V')}, for all pairs of integers $(p',q')$ and
$(p'',q'')$, with $0\leq p'\leq q'\leq n$ and $0\leq p''\leq q''\leq
n$, we have
\[
  \chi \left(\delta_{[p'',q'']}, \delta_{[p',q']}\right) 
  = \# \left([p',q']\cap [p'',q''] \right)
  - \# \left([p'+1,q'+1]\cap [p'',q''] \right), 
\]
where, given any set $S$, $\# S$ denotes its cardinal.

\begin{proposition}
\label{prop:coho-indec-rep}
Given pairs of integers $(p',q')$ and $(p'',q'')$, with $0\leq p'\leq
q'\leq n$ and $0\leq p''\leq q''\leq n$, the following conditions are
equivalent:
\par {\rm i)}
 $\Hom\left(\delta_{[p'',q'']},\delta_{[p',q']}\right)=0
 =\Ext^1\left(\delta_{[p'',q'']},\delta_{[p',q']}\right)$.
\par {\rm ii)} 
 At least one of the following inequalities is satisfied
 \begin{equation}
 \label{eq:Hom=0-indec-rep}
   p''>p',\qquad q''>q',\qquad p'>q'',\qquad p''>q', 
 \end{equation}
 and, furthermore, 
\begin{equation}
\label{eq:chi=0}
  \# \left([p',q']\cap [p'',q''] \right)
  = \# \left([p'+1,q'+1]\cap [p'',q''] \right). 
\end{equation}
\end{proposition}

\begin{proof}
We first show that $\Hom\left(\delta_{[p'',q'']},
  \delta_{[p',q']}\right)=0$ if and only if one of the
inequalities~\eqref{eq:Hom=0-indec-rep} is satisfied. Note that the
only way that the following two diagrams can commute is that the maps
$f$ and $g$ are zero.
\[
  \xymatrix@C=20pt@R=20pt{
  \CC \ar[r] \ar[d]_-{f} & 0 \ar [d]
\\
  \CC \ar@{=}[r] & \CC}
\qquad\qquad
  \xymatrix@C=20pt@R=20pt{
  \CC \ar@{=}[r] \ar[d] & \CC \ar [d]^-{g} 
  \\
  0 \ar[r] & \CC.}
\]
Hence, all the maps $f_i$ and $g_i$ in the following diagrams are
zero, provided they commute, so that they define morphisms
$\delta_{[p'',q'']} \to \delta_{[p',q']}$.
\[
\text{(Case $p''>p'$)}
\qquad\qquad
\vcenter{\xymatrix@C=20pt@R=20pt{
  \cdots \ar@{=}[r] & \CC \ar@{=}[r] \ar[d]_-{f_{p''+2}}  & \CC \ar@{=}[r] \ar[d]_-{f_{p''+1}}  & \CC \ar[r] \ar[d]_-{f_{p''}} & 0 \ar[d] \ar[r] & 0 \ar[d] \ar[r] & \cdots
  \\ 
  \cdots \ar@{=}[r] & \CC \ar@{=}[r] & \CC \ar@{=}[r] & \CC \ar@{=}[r]  & \CC \ar@{=}[r] & \CC \ar@{=}[r] & \cdots
}}
\]
\[
\text{(Case $q''>q'$)}
\qquad\qquad
\vcenter{\xymatrix@C=20pt@R=20pt{
  \cdots \ar@{=}[r] & \CC \ar@{=}[r] \ar[d]  & \CC \ar@{=}[r] \ar[d]_-{g_{q'}}  & \CC \ar@{=}[r] \ar[d]_-{g_{q'-1}} & \CC \ar[d]_-{g_{q'-2}} \ar@{=}[r] & \CC \ar[d] \ar@{=}[r] & \cdots
  \\
  \cdots \ar[r] & 0 \ar[r] & \CC \ar@{=}[r] & \CC \ar@{=}[r]  & \CC \ar@{=}[r] & \CC \ar@{=}[r] & \cdots
}}
\]
Thus, if $p''>p'$ or $q''>q'$, then $\Hom\left(\delta_{[p'',q'']},
  \delta_{[p',q']}\right)=0$.
If $p'>q''$ or $p''>q'$, then $\Hom\left(\delta_{[p'',q'']},
  \delta_{[p',q']}\right)=0$ as well, because the set $[p',q']\cap
[p'',q'']$ is empty.
Conversely, if none of the inequalities~\eqref{eq:Hom=0-indec-rep}
holds, i.e., if $q'\geq q''\geq p'\geq p''$ (so $[p',q']\cap
[p'',q'']\neq\varnothing$), then the following commutative diagram
shows that $\Hom\left(\delta_{[p'',q'']}, \delta_{[p',q']}\right)=\CC$.
\[
\xymatrix@C=20pt@R=20pt{
  \cdots \ar[r] & 0 \ar[r] \ar[d]  & 0 \ar[r]
  \ar[d]  & \CC \ar@{=}[r] \ar@{=}[d] & \CC
  \ar@{=}[d] \ar@{=}[r] & \cdots \ar@{=}[r] & \CC \ar@{=}[d] \ar@{=}[r] & \CC \ar[d] \ar@{=}[r] & \CC \ar[d] \ar@{=}[r] &\cdots
  \\
  \cdots \ar@{=}[r] & \CC \ar@{=}[r] & \CC \ar@{=}[r] & \CC \ar@{=}[r]  & \CC
  \ar@{=}[r] & \cdots \ar@{=}[r] &\CC \ar[r] & 0 \ar[r] & 0 \ar[r] & \cdots
}
\]
Finally, if $\Hom\left(\delta_{[p'',q'']},
  \delta_{[p',q']}\right)=0$, then $\Ext^1\left(\delta_{[p'',q'']},
  \delta_{[p',q']}\right)=0$ if and only if $\chi
\left(\delta_{[p'',q'']}, \delta_{[p',q']}\right)=0$, which is equivalent
to \eqref{eq:chi=0}, by the observation before this proposition.
\end{proof}

\begin{remark}
As observed in the proof of Proposition~\ref{prop:coho-indec-rep},
given pairs of integers $(p',q')$ and $(p'',q'')$ with $1\leq p'\leq
q'\leq n$ and $1\leq p''\leq q''\leq n$, the fact that at least one
of the inequalities~\eqref{eq:Hom=0-indec-rep} is satisfied is
equivalent to the fact that the following does not hold:
\begin{equation}\label{eq:Hom-neq-0-indec-rep}
  q'\geq q''\geq p'\geq p''. 
\end{equation} 
\end{remark}

The regions $\cV(r)$ and $\wt{R}(t)$ can now be obtained by applying
Propositions~\ref{prop:isom-indec} and~\ref{prop:coho-indec-rep} to all
possible direct sums~\eqref{eq:decomposition-rep} with dimension vectors
$r'$ and $r''$ such that $r'+r'' = r$. Thus, we have proved the following.

\begin{theorem}
\label{thm:BirRegionFund}
 Fix a dimension vector $r=(r_j,j=0,\ldots,n)$. A pair of dimension
 vectors $(r',r'')$, with $r'+r''=r$, belongs to $\cV(r)$ if and only
 if there are decompositions
 \[
 r'=r_{[p'_1,q'_1]} + \cdots + r_{[p'_{s'},q'_{s'}]} 
 \qquad \text{ and } \qquad 
 r''=r_{[p''_1,q''_1]} + \cdots + r_{[p''_{s''},q''_{s''}]},
 \]
 for two sequences of pairs of integers 
 \[
 ((p'_i,q'_i); i=1, \ldots, s') 
 \qquad \text{ and } \qquad 
 ((p''_j,q''_j); j=1,\ldots, s''), 
 \]
 with $0\leq p'_i\leq q'_i\leq n$ and $0\leq p''_j\leq q''_j\leq n$,
 such that the following conditions hold for all $i=1, \ldots, s'$
 and $j=1,\ldots, s''$. At least one of the following inequalities is
 satisfied
\begin{equation}
\label{eq:region-ineq}
   p''_j>p'_i,\qquad q''_j>q'_i,\qquad p'_i>q''_j,\qquad p''_j>q'_i, 
\end{equation}
 and, furthermore,
\begin{equation}
\label{eq:region-eq}
  \# \left([p'_i,q'_i]\cap [p''_j,q''_j] \right)
  = \# \left([p'_i+1,q'_i+1]\cap [p''_j,q''_j] \right). 
\end{equation}
\end{theorem}

\begin{example}
\label{ex:triples}
  As an application of Theorem~\ref{thm:BirRegionFund}, here we
  obtain the birationality region when $n=1$.  Note that in this case
  the $(n+1)$-holomorphic chains are actually the holomorphic triples
  studied in~\cite{BGG} and that in that paper, although not
  explicitly defined, the birationality region was completely
  determined.
  First, we list the only indecomposable linear $2$-chains $\delta',
  \delta''$ satisfying $\Hom(\delta'',\delta')= 0 =
  \Ext(\delta'',\delta')$:
  \begin{enumerate}
  \item[(1)] $\delta'=\delta_{[1,1]}$, $\delta''=\delta_{[0,0]}$; 
  \item[(2)] $\delta'=\delta_{[0,1]}$, $\delta''=\delta_{[1,1]}$; 
  \item[(3)] $\delta'=\delta_{[0,0]}$, $\delta''=\delta_{[0,1]}$. 
  \end{enumerate}
  Applying Theorem~\ref{thm:BirRegionFund}
  to a dimension vector $r=(r_0, r_1)$, it follows that $\cV(r)$ is
  the set of pairs $(r',r'')$, formed by pairs of integers
  $r'=(r'_0,r'_1)$ and $r''=(r''_0,r''_1)$ with $r'+r''=r$,
  satisfying exactly one of the following conditions:
  \begin{enumerate}
  \item[(1)] $r'_0=r''_1=0$;  
  \item[(2)] $r''_0=0$, $r''_1\neq 0$ and $r'_0=r'_1$; 
  \item[(3)] $r'_0\neq 0$, $r'_1=0$ and $r''_0=r''_1$. 
  \end{enumerate}
  These conditions also follow from of~\cite[Lemma~4.5]{BGG}. 
  Fix now a type $t=(r_0,r_1;d_0,d_1)$. Using the previous
  description of $\cV(r)$, together with
  Definition~\ref{def:region-V}(3) and Proposition~\ref{prop:Rtilda},
  we immediately see that
\[
  \wt{R}(t)\cap \{(0,\alpha)\in\R^2 \} = \{ 0\} \times
  \left(\R\setminus\{\alpha_m(t),\alpha_M(t)\}\right),
\]
  where
  \begin{align*}
    \alpha_m(t) &:= \frac{d_0}{r_0}-\frac{d_1}{r_1},\\
    \alpha_M(t) &:= \frac{2 r_0}{| r_0-r_1 |} \alpha_m(t)
    =\left( 1+ \frac{r_0+r_1}{| r_0-r_1 |}\right)
    \left(\frac{d_0}{r_0}-\frac{d_1}{r_1}\right) .
  \end{align*}
  Here we set $\alpha_M(t)=+\infty$ when $r_0=r_1$, by convention.
  Hence, a connected component of $\wt{R}(t)\cap \{(0,\alpha)\in\R^2
  \}$ is given by the open interval $\alpha_m(t)< \alpha
  <\alpha_M(t)$. We have thus recovered~\cite[Lemma~4.6]{BGG}.
\end{example}

In the previous example we have seen that, when $n=1$, the parameter
region $R(t)$ is in fact a component of the birationality region
$\wt{R}(t)$.
For $n>1$, this is not generally true, as we will see
in~\secref{sec:extensions2}.  The following proposition clarifies the
relation between the parameter region $R(t)$ and the birationality
region $\wt{R}(t)$. Recall that $R(t)$ is bounded by the hyperplanes
$h_i$ defined in~\secref{sub:parameter-region} for $0 \leq i < n$ and
that, by Proposition~\ref{prop:Rtilda}, $\wt{R}(t)$ is the
complementary of the union $\cB(t)$ of the hyperplanes $\cB(t',t'')$
of Definition~\ref{def:region-V}.

\begin{proposition}
\label{prop:parameter-region-subset-birational-region}
  Fix a type $t=(r_j,j=0,\ldots,n; d_j, j=0,\ldots,n)$. 
  The hyperplanes $h_i$ are contained in the boundary $\cB(t)$ of
  $\wt{R}$, for all $0 \leq i < n$. 
\end{proposition}

  \begin{proof}
  Given a linear chain $V$ and an integer $k\geq 0$, let $kV$ be the
  direct sum of $k$ copies of $V$ if $k>0$ or $0$ if $k=0$.  Observe
  now that $\Hom (\delta_{[j,j]}, \delta_{[k,k]}) = 0 = \Ext^1
  (\delta_{[i,i]}, \delta_{[j,j]})$ for all $0\leq j<k\leq n$ (this
  corresponds to case (1) in Example~\ref{ex:triples}), so
\[
  \Hom \!\left(\bigoplus_{0\leq j\leq i}r_j\delta_{[j,j]},  \!\!\bigoplus_{i < j\leq n}r_j\delta_{[j,j]}\right) 
 = 0 
  = \Ext^1 \!\left(\bigoplus_{0\leq j\leq i}r_j\delta_{[j,j]},
  \!\!\bigoplus_{i < j\leq n} r_j\delta_{[j,j]}\right).
\]
  for all $0 \leq i < n$. Since the linear chains 
  \[
  \bigoplus_{0\leq j\leq i}r_j\delta_{[j,j]} \text{ and } \bigoplus_{i < j\leq n}r_j\delta_{[j,j]}
\]
  have dimension vectors 
\[
  r_{(\leq i)}:= (r_0,\ldots,r_i,0,\ldots,0)  \text{ and } 
  r_{(>i)}:= (0,\ldots,0,r_{i+1}, \ldots,r_n), 
\]
  respectively, it follows that $(r_{(>i)}, r_{(\leq i)})\in
  \cV(r)$, so $(t_{(>i)}, t_{(\leq i)})\in \cT(t)$, where 
\[
  t_{(\leq i)}:=(r_{(\leq i)}, d_{(\leq i)})  \text{ and } 
  t_{(> i)}:=(r_{(> i)}, d_{(> i)}), 
\]
  with 
\[
  d_{(\leq i)} := (d_0,\ldots, d_i,0,\ldots,0) \text{ and }
  d_{(>i)} := (0,\ldots,0, d_{i+1}, \ldots, d_n). 
\]
  By definition, this implies that $\cB(t_{(>i)}, t_{(\leq i)})$ is
  contained in the boundary $\cB(t)$. But $\cB(t_{(>i)}, t_{(\leq
  i)})$ is the set of points $\alpha$ on the hyperplane
\[
  \mu_\alpha(t_{(>i)}) =  \mu_\alpha(t_{(\leq i)})
\]
  Writing explicitly this equation, we see that $\cB(t_{(>i)},
  t_{(\leq i)})$ is in fact the hyperplane $h_i$ defined
  in~\secref{sub:parameter-region}.
\end{proof}


\subsection{Birationality of moduli spaces}
\label{sec:birationallity}

Let $\alpha_w, \alpha_w^+$ and $\alpha_w^-$ be defined as
in~\secref{sub:flip-loci}, where now $|\epsilon| > 0$ is small enough
so that $\alpha_w^-$ and $\alpha_w^+$ are in the same connected
component of $\wt{R}(t)$.

\begin{proposition}\label{codim}
Let $w$ be a wall contained in the region
$\wt R(t)\cap R_{2g-2}$.
Let ${\mathcal{S}}$ be a  family of
$\alpha_w$-semistable chains   of type $t$ all of which are pairwise
non-isomorphic, and such that ${\mathcal{S}}$
maps generically one-to-one in an open set in the moduli space
$\mathcal{M}_{\alpha_w^\pm}$.
Then the codimension of the strictly semistable locus (which we
assume non-empty) in ${\mathcal{S}}$ is at least $g-1$.
\end{proposition}

\begin{proof}

The codimension of the stricly semistable locus is at least
\begin{equation}\label{min}
\min \left\{-\sum_{j<i} \chi(t_j,t_i)+\frac{m(m-3)+2}{2} \right\},
\end{equation}
where the minimum is taken over all the numerically possible
types $t_i$ and $m$ that may occur for a Jordan-H\"older filtration of  a
strictly $\alpha_w$-semistable change of type $t$.
This follows  from substracting (\ref{dim-family}) to the dimension
of the moduli space of $\alpha_w^\pm$-stable chains of type $t$
which is
$1-\chi(t,t)$ (recall that under the hypotheses of the proposition
$h^2(C,C)$ vanishes), and using that $\chi(t,t)=\sum_{i,j}\chi(t_i,t_j)$
by Corollary \ref{cor:chi-relation}.

Now, from Theorem~\ref{thm:chi-bounds}, we have that
$-\chi(t_j,t_i)\geq g-1$. Hence, the codimension is at least
$$
\min \left\{\frac{m(m-1)}{2}(g-1)+\frac{m(m-3)+2}{2}\right\}.
$$
Clearly, the minimum is attained when $m=2$ giving the result.
\end{proof}

From Proposition \ref{codim} we immediately  obtain the following.

\begin{theorem}\label{thm:birational}
Let $w$ be a wall contained in the region
$\wt R(t)\cap R_{2g-2}$.  Then
$\mathcal{M}^s_{\alpha_w^{+}}$ and   $\mathcal{M}^s_{\alpha_w^{-}}$
are birationally equivalent. Moreover, if in addition,
$\gcd(r_1,...,r_n,d)=1$
$\alpha$-stability coincides with  $\alpha$-semistability,
by \rm  Corollary \ref{CoPrimalityII}, 
\it  and
hence $\mathcal{M}_{\alpha_w^{+}}$ and $\mathcal{M}_{\alpha_w^{-}}$
are birationally equivalent.
\end{theorem}

\begin{remark}
\label{rem:discussion:parameter-region-vs-birational-region}
We emphasize that the boundary $\cB(t)$ of the birational
region $\wt{R}(t)$ contains in general more hyperplanes than just the 
$h_i$ of Proposition~\ref{prop:parameter-region-subset-birational-region}.
However, in the examples considered in~\secref{sec:extensions2} we
will be able to bound the dimensions of the flip loci when the
parameter $\alpha$ crosses the hyperplanes $\cB(t',t'')$ which are
inside the parameter region $R(t)$. A natural question is whether
this can always be done in the analysis of moduli spaces
$\mathcal{M}^s_{\alpha}(t)$ of holomorphic $(n+1)$-chains for all
possible types $t$ and all the hyperplanes $\cB(t',t'')$ which are
inside $R(t)$.
\end{remark}

\section{Parameter regions for semistable 3-chains}
\label{sec:3-chains}                            
In this section, we study the region of possible stability  parameters for holomorphic
3-chains. Among other things, we will prove Conjecture
\ref{ParBounds} for $n=2$. To this end, we first study linear 3-chains over the field $k=\CC(X)$.
A general assumption will be $r_j>0$, $j=0,1,2$.
\subsection{3-Chains of $k$-vector spaces}
Recall that we have established in Theorem \ref{AlongALine} a
connection between the semistability of a holomorphic chain and the
semistability of the induced linear chain over $k:=\CC(X)$. Therefore, we now study linear 3-chains over $k$.
\begin{theorem}
\label{UnstableQuivers} {\rm i)} Suppose $(r_0,r_1,r_2)$ satisfies
$r_0>r_1\neq r_2$. If $V=(V_0,V_1,V_2;f_1,f_2)$ is an
$(\alpha_1,\alpha_2)$-semistable linear 3-chain of type $(r_0,r_1,\allowbreak r_2)$, then
$$
(\alpha_1,\alpha_2)=(0,0).
$$
\par
{\rm ii)} Assume that $r_0<r_1>r_2$ and that
$V=(V_0,V_1,V_2;f_1,f_2)$ is an $(\alpha_1,\alpha_2)$-semistable
3-chain of type $(r_0,r_1,r_2)$, then $(\alpha_1,\alpha_2)=0$ or
$r_0=r_2$, $f_1\circ f_2$ is an isomorphism, and
$$
(\alpha_1,\alpha_2)=\la\cdot (1,2)
$$
for some $\la\in \R_{\ge 0}$. If $V$ is a chain, such that $r_0=r_2$
and $f_1\circ f_2$ is an isomorphism, then $V$ is $(1,2)$-semistable.
\par
{\rm iii)} If, in $(r_0,r_1,r_2)$, $r_0\neq r_1=r_2$ and if
$V=(V_0,V_1,V_2;f_1,f_2)$ is an $(\alpha_1,\alpha_2)$-semistable
3-chain of type $(r_0,r_1,r_2)$, then $f_2$ is an isomorphism and
$$
(\alpha_1,\alpha_2)=\la\cdot (-1,1)
$$
for some $\la\in\R_{\ge 0}$. Conversely, if $V$ is a chain in which
$f_2$ is an isomorphism, then $V$ is $(-1,1)$-semistable.
\par
{\rm iv)} Assume $r_0=r_1=r_2$. If $V=(V_0,V_1,V_2;f_1,f_2)$ is an
$(\alpha_1,\alpha_2)$-semistable linear 3-chain of type $(r_0,r_1,r_2)$, then either
$(\alpha_1,\alpha_2)=(0,0)$ or $f_2$ is an isomorphism and
$(\alpha_1,\alpha_2)$ is a non-negative multiple of $(-1,1)$ or
$f_1$ is an isomorphism and $(\alpha_1,\alpha_2)=\la\cdot (2,1)$ for some
non-negative number $\la$ or both $f_1$ and $f_2$ are isomorphisms,
$\alpha_1\le 2\alpha_2$, and $\alpha_1+\alpha_2\ge 0$. Conversely, a chain in which $f_1$
and $f_2$ are both ismorphisms is $(\alpha_1,\alpha_2)$-semistable
for all $(\alpha_1,\alpha_2)$ with
$\alpha_1\le 2\alpha_2$ and $\alpha_1+\alpha_2\ge 0$, a
chain in which $f_2$ is an isomorphism is $(-1,1)$-semistable, and a
chain in which $f_1$ is an isomorphism is $(2,1)$-semistable.
\end{theorem}
\begin{proof}[Proof of Theorem \ref{UnstableQuivers}]
We will check the condition of semistability on several non-trivial
subobjects. The first one is $(V_0,0,0)$. Semistability yields
\begin{equation}
\label{FirstInEQ} \alpha_1 r_1+\alpha_2 r_2 \ge 0.
\end{equation}
From the subobject $(V_0,V_1,0)$, we get the condition
\begin{equation}
\label{SecondInEQ} -\alpha_1 r_1+\alpha_2 (r_0+r_1) \ge 0.
\end{equation}
\par
 i) Since we assume $r_0>r_1$, $(\im(f_1),V_1,V_2)$ is a
non-trivial subchain. If $r_0^\prime<r_0$ is the dimension of
$\im(f_1)$, this subobject yields
$(r_0-r_0^\prime)(\alpha_1r_1+\alpha_2r_2)\le 0$. We infer
\begin{equation}
\label{ThirdInEQ} \alpha_1 r_1+\alpha_2 r_2 \le 0.
\end{equation}
Observe that (\ref{FirstInEQ}) and (\ref{ThirdInEQ}) give
\begin{equation}
\label{FirstEQ} \alpha_1 r_1=-\alpha_2 r_2.
\end{equation}
\par
Now, assume $r_2>r_1$. Then, $(0,0,\ker(f_2))$ is a non-trivial
subchain. We see
\begin{equation}
\label{FourthInEQ} -\alpha_1 r_1+\alpha_2 (r_0+r_1) \le 0.
\end{equation}
Together with (\ref{SecondInEQ}), we conclude
$\alpha_1r_1=\alpha_2(r_0+r_1)$. Invoking (\ref{FirstEQ}), we find $\alpha_2=0$ and derive
the assertion of the proposition in the case $r_0>r_1<r_2$.
\par
Finally, we have to consider the case $r_2<r_1$. In this case, we
use the condition of semistability for the subchain
$(V_0,\im(f_2),V_2)$. This condition is
$\alpha_2(r_1-r_1^\prime)r_2\le \alpha_1(r_1-r_1^\prime)(r_0+r_2)$,
$r_1^\prime:=\dim(\im(f_2))$. Since $r_1^\prime\le r_2<r_1$, this
gives
\begin{equation}
\label{FifthInEQ} \alpha_1 (r_0+r_2)-\alpha_2 r_2 \ge 0.
\end{equation}
Together with (\ref{FirstInEQ}), this yields $\alpha_1\ge 0$. On the
other hand, (\ref{SecondInEQ}) shows that $\alpha_2\ge 0$. By
(\ref{ThirdInEQ}), we have $(\alpha_1,\alpha_2)=(0,0)$ as desired.
\par
 ii) If $r_0<r_1>r_2$, then there is the subchain
$(0,\ker(f_1),0)$ which gives the estimate
\begin{equation}
\label{30InEQ} \alpha_1(r_0+r_2)\le \alpha_2r_2.
\end{equation}
Inequality (\ref{FifthInEQ}) is also true, so that
$$
\alpha_1(r_0+r_2)=\alpha_2r_2.
$$
Note that $\alpha_2\ge 0$, by (\ref{FirstInEQ}) and
(\ref{SecondInEQ}), so that also $\alpha_1\ge 0$. One can now check
that a chain $V$ is $(0,\alpha_1,\alpha_2)$-semistable, if and only
if every subchain $(W_0,W_1,W_2)$ with $\dim(W_j)=r_j^\p$,
$j=0,1,2$, satisfies
\begin{equation}
\label{StabCond01} r_0r_2^\p\le r_0^\p r_2.
\end{equation}
If $f_1\circ f_2$ were not injective, then this condition would be
violated by $(0,f_2(\ker(f_1\circ f_2)),\ker(f_1\circ f_2))$. Thus, $f_1\circ f_2$ is
injective and $r_0\ge r_2$. Since the dual chain $V^\vee$ is
$(0,\alpha_2-\alpha_1,\alpha_2)$-semistable (compare Remark
\ref{quotientformulation}), we must also have $r_2\ge r_0$ and,
consequently, $r_0=r_2$. If we assume conversely that $r_0=r_2$ and
that $f_1\circ f_2$ is an isomorphism, then $r_2^\p \le r_0^\p$ for
every subchain, and (\ref{StabCond01}) is verified.
\par
 iii) First, assume that $r_0>r_1$. Then, (\ref{ThirdInEQ}) is
still satisfied. Thus, by (\ref{FirstEQ}) and $r_1=r_2$,
$\alpha_1=-\alpha_2$. Since $\alpha_2\ge 0$, we see that
$(\alpha_1,\alpha_2)$ is a non-negative multiple of $(-1,1)$ and, if it is a positive multiple, the
semistability condition becomes
\begin{equation}
\label{StabCond11} \dim(W_2)\le \dim(W_1)
\end{equation}
for all subchains $(0,W_1,W_2)$. This would be violated by
$(0,0,\ker(f_2))$, if $f_2$ were not injective. Conversely, if $f_2$
is injective, then (\ref{StabCond11}) will obviously be satisfied for
all subchains.
\par
For $r_0<r_1$, $W_1:=\ker(f_1)$ is a non-trivial subspace of $V_1$
of dimension, say, $r_1^\p$. Choose a subspace $W_2$ of $V_2$ of
dimension $r_1^\p$ which maps to $W_1$ under $f_1$. The subchain
$(0,W_1,W_2)$ yields
$$
(\alpha_1+\alpha_2)(r_0+2r_1)\le 2(\alpha_1+\alpha_2)r_1,
$$
i.e., (\ref{ThirdInEQ}) is again verified, and we may conclude as
before.
\par
 iv).
A chain $V$ is $(0,\alpha_1,\alpha_2)$-semistable if and only if it
is $\alpha^\p$-semistable for
$\alpha^\p=(\alpha_0^\p,\alpha_1^\p,\alpha_2^\p)$ with
$\alpha_0^\p:=-(\alpha_1+\alpha_2)/3$,
$\alpha_1^\p:=(2\alpha_1-\alpha_2)/3$, and
$\alpha_2^\p:=(-\alpha_1+2\alpha_2)/3$. Note that we have
$\alpha_0^\p+\alpha_1^\p+\alpha_2^\p=0$. We will require $\alpha'\neq 0$ in the following.
For $\alpha_0^\p=0$, i.e., $(\alpha_1,\alpha_2)=\la\cdot (-1,1)$ for
some $\la\in\R_{> 0}$, $\alpha^\p$-semistability is given by
(\ref{StabCond11}) for all subchains. As we have discussed before,
this is equvalent to $f_2$ being an isomorphism. For
$\alpha_2^\p=0$, that is $(\alpha_1,\alpha_2)=\la\cdot (2,1)$,
$\la> 0$, we find the condition
\begin{equation}
\label{StabCond13} \dim(W_0)\le \dim(W_1)
\end{equation}
for all subchains. This is equivalent to $f_1$ being an isomorphism.
In the case that $\alpha_1^\p=0$, i.e., $(\alpha_1,\alpha_2)=\la\cdot (1,2)$, $\la>0$, the condition of
$\alpha^\p$-semistability becomes
\begin{equation}
\label{StabCond14} \dim(W_2)\le \dim(W_0)
\end{equation}
for all subchains. This is equivalent to the fact that
$f_1\circ f_2$ is an isomorphism, so that both $f_1$ and
$f_2$ must be isomorphisms.
\par
Finally, we treat the case $\alpha_j'\neq 0$, $j=0,1,2$.
The conditions $\alpha_1+\alpha_2\ge 0$ and $\alpha_1\le 2\alpha_2$ result from (\ref{FirstInEQ})
and (\ref{SecondInEQ}), keeping in mind $r_0=r_1=r_2$.
If $f_1$ were not surjective, then (\ref{ThirdInEQ}) would hold and, thus, $\alpha_1=-\alpha_2$
and $\alpha_0'=0$, a contradiction. Similarly, we derive that $f_2$ must be an isomorphism.
To conclude, assume that $f_1$ and $f_2$ are isomorphisms, $\alpha_1\le 2\alpha_2$,
and $\alpha_1+\alpha_2\ge 0$. Let $(W_0,W_1,W_2)$ be a subchain with $\dim(W_j)=r_j'$, $j=0,1,2$.
Then, $r_0'\ge r_1'\ge r_2'$, and we may estimate as follows:
\begin{eqnarray*}
&&3(\alpha_0'r_0'+\alpha_1'r_1'+\alpha_2r_2')
\\
&=&
-(\alpha_1+\alpha_2)r_0'+(2\alpha_1-\alpha_2)r_1'+(-\alpha_1+2\alpha_2)r_2'
\\
&\le& -(\alpha_1+\alpha_2)r_1'+(2\alpha_1-\alpha_2)r_1'+(-\alpha_1+2\alpha_2)r_2'
\\
&=& (\alpha_1-2\alpha_2)r_1'+(-\alpha_1+2\alpha_2)r_2'
\\
&\le& (\alpha_1-2\alpha_2)r_2'+(-\alpha_1+2\alpha_2)r_2'=0.
\end{eqnarray*}
This shows that $V$ is a semistable chain.
\end{proof}
\begin{remark}
Observe that the proof shows that there never exists any stable
chain.
\end{remark}
\subsection{Parameter regions for rank maximal 3-chains}
\label{sub:Rank-Maximal}
A holomorphic chain will be called \it rank maximal\rm, if all the
homomorphisms, that is, $\phi_1$, $\phi_2$, and, in case $r_1\ge \max\{\,r_0,r_1\,\}$, also $\phi_1\circ\phi_2$,
have generically maximal rank, i.e., are either
injective or generically surjective. Note that being rank maximal is
an open property, so that describing the moduli spaces for rank
maximal chains will provide birational models for some components of
the moduli space of all chains. (Those components will be smooth in stable points by
Theorem \ref{thm:smoothdim}, v).) For rank maximal chains, one can use
the test objects analogous to those used in the proof of Proposition
\ref{UnstableQuivers} to find inequalities which limit the parameters
for which semistable rank maximal 3-chains might occur. The advantage is
that, for rank maximal chains, one obtains more explicit bounds on the
parameters.
We state the following more precise result.
\begin{proposition}
\label{prop:RankMaxChains}
Let $C=(E_0,E_1,E_2; \phi_1,\phi_2)$
be a holomorphic 3-chain of type $t=(r_0,r_1,r_2; d_0,d_1,d_2)$.
Define the inequalities
\begin{eqnarray}
\label{41InEQ}
\alpha_1r_1+\alpha_2r_2 &\ge& \frac{r_1+r_2}{r_0}\cdot d_0-d_1-d_2
\\
\label{42InEQ}
-\alpha_1r_1+\alpha_2(r_0+r_1) &\ge&  d_0+d_1-\frac{r_0+r_1}{r_2}\cdot d_2.
\end{eqnarray}
and, if applicable, the following additional inequalities:
\begin{enumerate}
\item If $r_0>r_1$ and $\phi_1$ is injective,
\begin{equation}
\label{43InEQ}
\alpha_1r_1+\alpha_2r_2\le
\frac{2r_1+r_2}{r_0-r_1}\cdot d_0-\frac{2r_0+r_2}{r_0-r_1}\cdot d_1-d_2.
\end{equation}
\item If $r_1< r_2$ and $\phi_2$ is generically surjective,
\begin{equation}
\label{44InEQ}
-\alpha_1r_1+\alpha_2(r_0+r_1)
\le  d_0+\frac{r_0+2r_2}{r_2-r_1}\cdot d_1-\frac{r_0+2r_1}{r_2-r_1}\cdot d_2.
\end{equation}
\item If $r_1> r_2$ and $\phi_2$ is injective,
\begin{equation}
\label{45InEQ}
-\alpha_1(r_0+r_2)+\alpha_2 r_2
\le  -d_0+\frac{r_0+2r_2}{r_1-r_2}\cdot d_1-\frac{r_0+2r_1}{r_1-r_2}\cdot d_2.
\end{equation}
\item If $r_0<r_1$ and $\phi_1$ is generically surjective,
\begin{equation}
\label{46InEQ}
\alpha_1(r_0+r_2)-\alpha_2 r_2
\le  \frac{2r_1+r_2}{r_1-r_0}\cdot d_0-\frac{2r_0+r_2}{r_1-r_0}\cdot d_1+d_2.
\end{equation}
\item $r_0<r_1> r_2$, $r_0<r_2$, and $\phi_1\circ\phi_2$ is generically surjective,
\begin{eqnarray}
&&\alpha_1(r_0-r_1+r_2)+\alpha_2 (r_0+r_1-r_2)\nonumber
\\
\label{47InEQ}
&\le&
\frac{-r_0+r_1+3r_2}{r_2-r_0}\cdot d_0+ 2d_1-\frac{3r_0+r_1-r_2}{r_2-r_0}\cdot d_2.
\end{eqnarray}
\end{enumerate}
Let $R(t)$ be the region described by all applicable inequalities. If $C$ is $(0,\alpha_1,\alpha_2)$-semistable,
then
$$
(\alpha_1,\alpha_2)\in R(t).
$$
\end{proposition}
\begin{proof}
The inequalities (\ref{41InEQ}) and (\ref{42InEQ}) have already been given
in Example \ref{ex:StandardBounds2}. Inequality (\ref{43InEQ}) follows from
the subchain $(\im(\phi_1)\cong E_1, E_1,E_2)$. For (\ref{44InEQ}), one uses the
subchain $(0,0,\ker(\phi_2))$ and the inequality $\deg(\ker(\phi_2))
=d_2-\deg(\im(\phi_2))\ge d_2-d_1$. The subchain
$(E_0, \im(\phi_2)\cong E_2, E_2)$ provides us with (\ref{45InEQ}). Inequality
(\ref{46InEQ}) is derived from the test object $(0,\ker(\phi_1),0)$ and
the fact $\deg(\ker(\phi_1))\ge d_1-d_0$. Last but not least,
the subchain $(0, \phi_2(\ker(\phi_1\circ\phi_2))\cong \ker(\phi_1\circ \phi_2), \ker(\phi_1\circ \phi_2))$
and the estimate $\deg(\ker(\phi_1\circ\phi_2))\ge d_2-d_0$ yield (\ref{47InEQ}).
\end{proof}
\begin{example}
\label{ex:Parallelo}
Let $C=(E_0,E_1,E_2; \phi_1,\phi_2)$
be a holomorphic 3-chain of type $t=(r_0,r_1,r_2; d_0,d_1,d_2)$, such that
$r_0<r_1> r_2$, $r_0<r_2$, and $\phi_1\circ\phi_2$ is generically surjective.
\begin{figure}[h]
\begin{center}
\begin{pspicture}(-2.5,.5)(2.5,6)
\psset{xunit=2cm, yunit=1cm}
\pspolygon*[linecolor=lblue](-2.5,.5)(-2.5,6)(2.5,6)(2.5,.5)
\pspolygon*[linecolor=lightorange](-2,1)(0,2.5)(2,5.5)(0,4)
\psline[linecolor=red, linewidth=1.5pt](-2.5,.625)(2.5,4.375)
\psline[linecolor=red, linewidth=1.5pt](-2.5,2.125)(2.5,5.875)
\psline[linecolor=red, linewidth=1.5pt](-2.33,.5)(1.33,6)
\psline[linecolor=red, linewidth=1.5pt](-1.33,.5)(2.33,6)
\psdot[linecolor=mix, linewidth=13pt](0,3.25)
\pstextpath[l](0,-.5){\psline[linestyle=none](0,2.65)(2.5,4.375)}{\red$-\alpha_1r_1+\alpha_2(r_0+r_1)=$}
\pstextpath[l](0,-1){\psline[linestyle=none](0,2.65)(2.5,4.375)}{\red$d_0+d_1-((r_0+r_1)/r_2)\cdot d_2$}
\pstextpath[r](0,.35){\psline[linestyle=none](-2.5,2.125)(0.3,4.2)}{\red$-\alpha_1r_1+\alpha_2(r_0+r_1)=K$}
\pstextpath[l](-.5,.125){\psline[linestyle=none](-2,1)(1.33,6)}
{\red\small$-\alpha_1(r_0+r_2)+\alpha_2r_2=L$}
\pstextpath[c](0,-.275){\psline[linestyle=none](0,2.5)(2.4,6)}
{\red\small$\alpha_1(r_0+r_2)-\alpha_2r_2=M$}
\uput[0](-.26,3.275){\bfseries\orange$R(t)$}
\uput[0](-.75, 1){$M:=\frac{2r_1+r_2}{r_1-r_0}\cdot d_0-\frac{2r_0+r_2}{r_1-r_0}\cdot d_1+d_2$}
\uput[0](-2.5, 5.5){$L:=-d_0+\frac{r_0+2r_2}{r_1-r_2}\cdot d_1-\frac{r_0+2r_1}{r_1-r_2}\cdot d_2$}
\end{pspicture}
\end{center}
\caption{The parameter region for rank maximal 3-chains with $r_0<r_1> r_2$ and $r_0<r_2$.}
\end{figure}
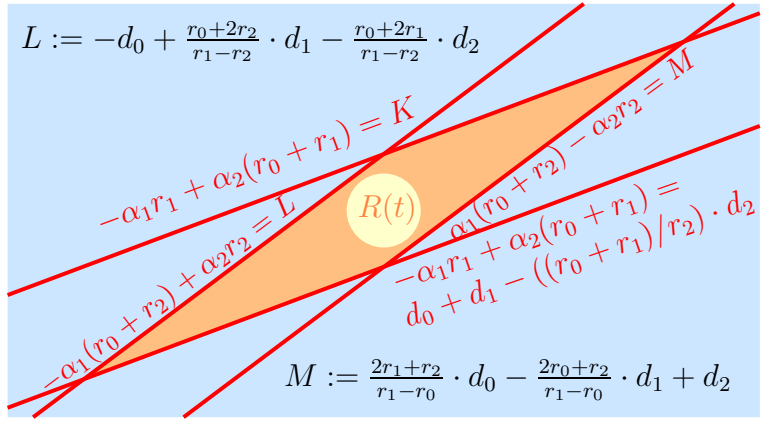
Then, obviously, $\phi_2$ must be injective and $\phi_1$ generically surjective.
If we add (\ref{45InEQ}) and (\ref{47InEQ}), we find an estimate
\begin{equation}
\label{eq:47aInEQ}
-\alpha_1r_1+\alpha_2(r_0+r_1)\le K=K(t).
\end{equation}
Note that
$$
\frac{r_1}{r_0+r_1}<1< \frac{r_0+r_2}{r_2}.
$$
It follows that the inequalities (\ref{42InEQ}), (\ref{eq:47aInEQ}), (\ref{45InEQ}),
and (\ref{46InEQ}) bound a region $R$ in the shape of a parallelogram (see Figure 2).
Thus, the possible stability parameters for rank maximal 3-chains of type$t=(r_0,r_1,r_2; d_0,d_1,d_2)$ with $r_0<r_1> r_2$ and $r_0<r_2$ live in the bounded region $R$.
(Note that the other inequalities may still ``cut away" some pieces from $R$.)
\end{example}
\subsection{Bounded parameter regions}
\label{sub:3-chains-bounds}
We now demonstrate
\begin{theorem}
\label{BoundedParameterRegion} {\rm Conjecture \ref{ParBounds}} holds true
for $n=2$.
\end{theorem}
\begin{proof}
Fix the type
$t=(r_0,r_1,r_2;d_0,d_1,d_2)$. The cases we have to consider are a) $r_0>r_1<r_2$,
b) $r_0>r_1>r_2$, and c) $r_0<r_1>r_2$ and $r_0<r_2$.
We would like to adapt the strategy in the proof of Proposition
\ref{UnstableQuivers}. However, if a test object contains a kernel
or an image of a map, then the semistability condition contains the
degree of that kernel or image which we do not know in general. So,
we have to modify some arguments.
\par
Let $C=(E_0,E_1,E_2;\phi_1,\phi_2)$ be an $\alpha$-semistable 3-chain.
From the Inequalities (\ref{41InEQ}) and (\ref{42InEQ}) which always hold,
we derive the estimate
$$
\alpha_2 \ge K_0:=\frac{d_0}{r_0}-\frac{d_2}{r_2}.
$$
\paragraph{\bfseries The case $r_0>r_1$}
We will give a bound on $\alpha_1$ under the condition that
$\phi_1$ is \sl not \rm injective. Suppose the rank of
$\ker(\phi_1)$ is $r_1^\prime$ and its degree is $d_1^\prime$. Then,
the semistability condition for the test object $(0,\ker(\phi_1),
0)$ reads
\begin{equation}
\label{9InEQ}
\begin{array}{rl}
&(r_0+r_1+r_2)d_1^\prime+\alpha_1r_1^\prime(r_0+r_1+r_2)
\\
\le &
r_1^\prime (d_0+d_1+d_2)+r_1^\prime (\alpha_1r_1+\alpha_2r_2).
\end{array}
\end{equation}
We invoke the condition arising from the subchain
$(\im(\phi_1),E_1,E_2)$, i.e.,
\begin{equation}
\label{10InEQ}
\begin{array}{rl}
&(r_0+r_1+r_2)(d_1+d_2)+(\alpha_1r_1+\alpha_2r_2)(r_0+r_1+r_2)
\\
&+(r_0+r_1+r_2)(d_1-d_1^\prime)
\\
\le& (r_1+r_2)(d_0+d_1+d_2)+(r_1+r_2)(\alpha_1r_1+\alpha_2r_2)
\\
&+(r_1-r_1^\prime)(d_0+d_1+d_2)+
(r_1-r_1^\prime)(\alpha_1r_1+\alpha_2r_2).
\end{array}
\end{equation}
Now, add (\ref{9InEQ}) and (\ref{10InEQ}) in order to find
$$
-\alpha_1(r_0+r_1+r_2)r_1^\prime\ge
K_1+(\alpha_1r_1+\alpha_2r_2)(r_0-r_1).
$$
Here, $K_1:=-(2r_1+r_2)d_0+(2r_0+r_2)d_1+(r_0-r_1)d_2$.
Together with (\ref{41InEQ}), we find
$$
-\alpha_1r_1^\prime\ge K_2:=-\frac{r_1}{r_0}\cdot d_0+d_1,
$$
i.e.,
$$
\alpha_1\le K_3:=-\min\biggl\{\, r_1^\prime=1,...,r_1\,\Big\vert\,
\frac{K_2}{r_1^\prime}\biggr\}.
$$
\par
If, on the other hand, $\phi_1$ is injective, then we have (\ref{43InEQ}).
Now, $-r_2\cdot$(\ref{42InEQ})$+(r_0+r_1)\cdot$(\ref{43InEQ})
provides the estimate
$$
\alpha_1 \le K_4:=\frac{2}{r_0-r_1}\left(d_0-\frac{r_0}{r_1}\cdot d_1\right).
$$
\paragraph{\bfseries The case $r_0>r_1<r_2$} First, we will derive a bound on
$\alpha_2$ under the condition that $\phi_2$ is not generically
surjective. Note that the latter is equivalent to the fact that
$\phi_2^\vee$ is not injective. Thus, we may look at the dual
holomorphic 3-chain $C^\vee$ which is
$(-\alpha_2,-\alpha_1,0)$-semistable, by Remark \ref{quotientformulation}, iv), and thus
$(0,\alpha_2-\alpha_1,\alpha_2)$-semistable, by Remark \ref{quotientformulation}, iii).
Our previous
computations may now be applied to find a bound
$$
\alpha_2-\alpha_1\le K_5:=\frac{2}{r_2-r_1}\left(\frac{r_2}{r_1}\cdot d_1-d_2\right).
$$
and, thus,
\begin{equation}
\label{11InEq}
 \alpha_2\le K_5+  \alpha_1\le K_5+K_4.
\end{equation}
If $\phi_2$ is generically surjective, we have (\ref{44InEQ}).
Since we have already bounded $\alpha_1$ from above, Inequality (\ref{44InEQ})
provides an upper bound for $\alpha_2$. All in all, we have found
lower and upper bounds for both $\alpha_1$ and $\alpha_2$ under the
assumption a) $r_0>r_1<r_2$.
\paragraph{\bfseries The case $r_1>r_2$}
We abbreviate $r:=
r_0+r_1+r_2$, $\mu:=(d_0+d_1+d_2)/r$, and $\mu_K:=d_2^\p/r_2^\p$,
$r_2^\p:=\rk(\ker(\phi_2))$ and $d_2^\p:=\deg(\ker(\phi_2))$.
From the subchain $(0,0,\ker(\phi_2))$, we get the condition
$$
\mu_K+\alpha_2\le \mu+\frac{\alpha_1r_1+\alpha_2r_2}{r},
$$
i.e.,
\begin{equation}
\label{20InEQ} \alpha_2-\frac{\alpha_1r_1+\alpha_2r_2}{r}\le
\mu-\mu_K.
\end{equation}
Next, we check the semistability condition for the subchain
$(E_0,\im(\phi_2),\allowbreak E_2)$. We find
$$
\frac{d_0+2d_2-r_2^\p\mu_K}{r_0+2r_2-r_2^\p}+\frac{\alpha_1(r_2-r_2^\p)+\alpha_2r_2}{r_0+2r_2-r_2^\p}
\le \mu+\frac{\alpha_1r_1+\alpha_2r_2}{r}.
$$
This may be rewritten as
$$
-r_2^\p\mu_K+\alpha_1(r_2-r_2^\p)+\alpha_2r_2\le -r_2^\p\mu+
\frac{r_0+2r_2-r_2^\p}{r}(\alpha_1r_1+\alpha_2r_2)+K_6,
$$
$K_6:=\mu(r_0+2r_2)-(d_0+2d_2)$, that is,
\begin{equation}
\label{21InEQ} \mu-\mu_K\le \frac{r_0+2r_2-r_2^\p}{r_2^\p
r}(\alpha_1r_1+\alpha_2r_2)+\frac{\alpha_1(r^\p_2-r_2)-\alpha_2r_2}{r_2^\p}+\frac{K_6}{r_2'}.
\end{equation}
We combine (\ref{20InEQ}) and (\ref{21InEQ}) and multiply by $r_2^\p
r$:
\begin{eqnarray*}
&&\alpha_2r_2^\p r-(\alpha_1r_1+\alpha_2r_2)r_2^\p
\\
&\le &
(r_0+2r_2-r_2^\p)(\alpha_1r_1+\alpha_2r_2)+\alpha_1(r_2^\p-r_2)r-\alpha_2r_2r+rK_6.
\end{eqnarray*}
We conclude
\begin{equation}
\label{48InEQ}
\alpha_2(r_2^\p r+r_2r-r_2(r_0+2r_2))\le \alpha_1(r_2^\p
r-r_2r+(r_0+2r_2)r_1)+rK_6.
\end{equation}
Observe $r-r_0-2r_2=r_1-r_2>0$, whence the coefficient of $\alpha_2$
is positive. Also,
$$
(r_0+2r_2)r_1-r_2r=(r_0+r_2)(r_1-r_2)>0,
$$
so that the coefficient of $\alpha_1$ is positive. (The value of $rK_6$
is $(r_2-r_1)d_0+(r_0+2r_2)d_1-(r_0+2r_1)d_2$ (compare (\ref{45InEQ})).)
\paragraph{\bfseries The case $r_0>r_1>r_2$}
Recall that we have already bounded $\alpha_1$ from above. Therefore, if $\phi_2$
is not injective, then (\ref{48InEQ}) provides an upper bound for $\alpha_2$, too.
If, on the other hand, $\phi_2$ is injective, Inequality (\ref{45InEQ}) holds.
This inequality also provides an upper bound of
$\alpha_2$ in terms of $\alpha_1$ and constants depending only on
the type $t$. Again, $\alpha_1$ and $\alpha_2$ are bounded both from above
and below, and we are done for case b).
\paragraph{\bfseries The case $r_0<r_1>r_2$; $r_0<r_2$}
By (\ref{48InEQ}), there is the inequality
$$
\alpha_2c_2\le \alpha_1c_1+K_7,
$$
with positive constants $c_1$ and $c_2$, if $\phi_2$ is not injective.
If $\phi_2$ is injective, we have Inequality (\ref{45InEQ}).
If $\phi_1$ is not generically surjective, then (\ref{48InEQ}) for the dual chain yields
$$
\alpha_2(r_0^\p r+r_0r-r_0(2r_0+r_2))\le (\alpha_2-\alpha_1)(r_0^\p
r-r_0r+(2r_0+r_2)r_1)+rK_8,
$$
$r_0':=\rk(\ker(\phi_1^\vee))$, i.e.,
\begin{equation}
\label{49InEQ}
\begin{array}{rl}
&\alpha_1(r_0^\p r-r_0r+(2r_0+r_2)r_1)
\\
=&\alpha_1(r_0^\p r+(r_0+r_2)(r_1-r_0))
\\
\le &
\alpha_2r_2(r_1-r_0)+rK_8.
\end{array}
\end{equation}
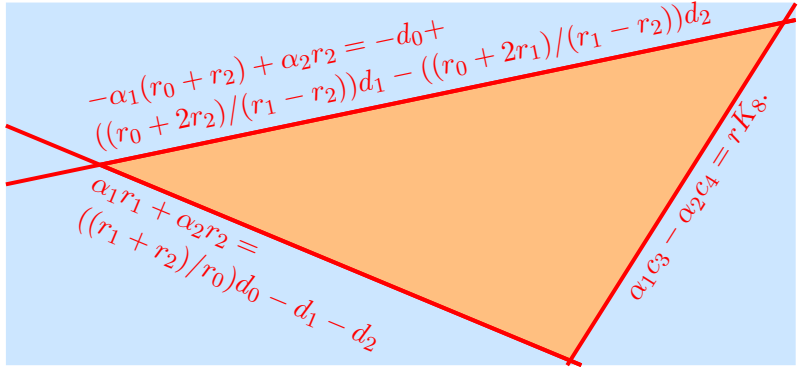
\begin{figure}[b]
\begin{center}
\begin{pspicture}(1,2)(2,6.5)
\psset{xunit=75cm, yunit=25cm}
\pspolygon*[linecolor=lblue](-.05,.06)(-.05,.253)(.09,.253)(.09,.06)
\pspolygon*[linecolor=lightorange](-.0333333,.166666)(.088,.2425)(.05,.0625)
\psline[linecolor=red, linewidth=1.5pt](-.05,.15625)(.09,.24375)
\pstextpath[l](0,.25){\psline[linecolor=red, linewidth=1.5pt](-.0333333,.166666)(.09,.24375)}
{\red\small$((r_0+2r_2)/(r_1-r_2)) d_1-((r_0+2r_1)/(r_1-r_2)) d_2$}
\pstextpath[l](0,.75){\psline[linestyle=none](-.0333333,.166666)(.09,.24375)}
{\red\small$-\alpha_1(r_0+r_2)+\alpha_2 r_2=-d_0+$}
\psline[linecolor=red, linewidth=1.5pt](-.05,.1875)(.052,.06)
\pstextpath[l](0,-.35){\psline[linestyle=none](-.0333333,.166666)(.052,.06)}
{\red\small$\alpha_1r_1+\alpha_2r_2 =$}
\pstextpath[l](0,-.85){\psline[linecolor=red, linewidth=1.5pt](-.0333333,.166666)(.052,.06)}
{\red\small$((r_1+r_2)/r_0) d_0-d_1-d_2$}
\pstextpath[c](0,-.35){\psline[linecolor=red, linewidth=1.5pt](.0495,.06)(.09,.2525)}
{\red\small$\alpha_1c_3-
\alpha_2c_4=rK_8.$}
\end{pspicture}
\end{center}
\caption{The triangle bounded by (\ref{41InEQ}), (\ref{45InEQ}), and (\ref{49InEQ}).}
\end{figure}
\\
Set $c_3:=(r_0^\p r+(r_0+r_2)(r_1-r_0))$ and $c_4:=r_2(r_1-r_0)$.
One easily checks
\begin{equation}
\label{eq:SomeIneqs}
\frac{c_3}{c_4}> \frac{r_0+r_2}{r_2} > \frac{c_1}{c_2}.
\end{equation}
Hence, if $\phi_1$ is not generically surjective, then (\ref{41InEQ}), (\ref{48InEQ}),
and (\ref{49InEQ}) bound a triangular region. The same goes for
(\ref{41InEQ}), (\ref{45InEQ}), and (\ref{49InEQ}) (cf.\ Figure 3).
In our argument, we may therefore assume that
$\phi_1$ is generically surjective, so that we have Inequality (\ref{46InEQ}).
By (\ref{eq:SomeIneqs}), (\ref{41InEQ}), (\ref{46InEQ}), and (\ref{48InEQ}) also bound a triangular region.
Thus, we may also assume that $\phi_2$ is injective. If $\phi_1\circ \phi_2$ is generically surjective,
we have Example \ref{ex:Parallelo}. The final case to consider is the one in which $\phi_1$ and $\phi_2$
have generically the maximal possible rank and $\phi_1\circ \phi_2$ has a cokernel of positive rank.
We use again the abbreviations $d:=d_0+d_1+d_2$, $r:=r_0+r_1+r_2$, and $\mu:=d/r$.
\begin{claim}
There is a constant $K_9=K_9(t)$, such that
$$
\mu-\mu\bigl(\ker(\phi_1\circ\phi_2)\bigr)\le K_9.
$$
\end{claim}
Assume, for the moment, this claim. Note that (\ref{45InEQ}) and (\ref{46InEQ}) yield
$$
\frac{\alpha_1(r_0+r_2)}{r_2}-K_{10}\le \alpha_2\le \frac{\alpha_1(r_0+r_2)}{r_2}+K_{11}.
$$
For a subchain $C'=(F_0,F_1,F_2)$ with $r_j':=\rk(F_j)$, $j=0,1,2$, $r':=r_0'+r_1'+r_2'$,
$d':=\deg(F_0)+\deg(F_1)+\deg(F_2)$, and
$\mu_{C'}:=d'/r'$, we then get from the condition of $(0,\alpha_1,\alpha_2)$-semistability
\begin{equation}
\label{eq:54InEQ}
\alpha_1(r_0r_2'-r_0'r_2)\le r'r_2(\mu-\mu_{C'})+r_2'r_2K_{10}+r'r_2\frac{K_{11}}{r}.
\end{equation}
If we apply this to the subchain $(0,\phi_2(\ker(\phi_1\circ\phi_2))\cong \ker(\phi_1\circ\phi_2),
\allowbreak \ker(\phi_1\circ \phi_2))$,
we find a bound
$$
\alpha_1\le K_{12}=K_{12}(t),
$$
using the above claim. This proves the boundedness of the parameter region.
\par
In order to establish the claim, we look at the dual chain $C^\vee=(\widetilde{E}_0,\widetilde{E}_1,\widetilde{E}_2;
\widetilde{\phi}_1,\widetilde{\phi}_2)$. If $\phi_1\circ \phi_2$ is not generically surjective,
then $\widetilde{\phi_1}\circ \widetilde{\phi_2}$ is not injective.
Set
$$
\widetilde{K}:=\widetilde{\phi}_2(\ker(\widetilde{\phi_1}\circ \widetilde{\phi_2}))
\cong \ker(\widetilde{\phi_1}\circ \widetilde{\phi_2}),
$$
because $\widetilde{\phi}_2$ is injective. If $d_K:=\deg(\ker(\phi_1\circ\phi_2))$,
$r_K:=\rk(\ker(\phi_1\circ\phi_2))$, and $\mu_{K}:=d_K/r_K$,
then
$$
\mu(\widetilde{K})=\frac{-d_0+d_2-r_K\mu_K}{r_0-r_2+r_K}.
$$
Since $C^\vee$ is $(0,\alpha_2-\alpha_1,\alpha_2)$-semistable, the subchain $(0,\widetilde{K},0)$
gives the estimate
$$
\mu(\widetilde{K})+\alpha_2-\alpha_1\le -\mu+\frac{-\alpha_1r_1+\alpha_2(r_0+r_1)}{r}.
$$
This may be rewritten as
\begin{eqnarray*}
r_K(\mu-\mu_K) &\le& (r_2-r_0)\mu+d_0-d_2+
\\
&&+\frac{r_0-r_2+r_K}{r}\bigl(\alpha_1(r_0+r_2)-\alpha_2r_2\bigr)
\\
&\le& (r_2-r_0)\mu+d_0-d_2+
\\
&&+\frac{r_0-r_2+r_K}{r}\Bigl(\frac{2r_1+r_2}{r_1-r_0}\cdot d_0-\frac{2r_0+r_2}{r_1-r_0}\cdot d_1+d_2\Bigr).
\end{eqnarray*}
For the second estimate, we have used the fact that $r_K\ge r_2-r_0$ and (\ref{46InEQ}).
The above inequality clearly settles the claim.
\end{proof}
\subsection{Unbounded parameter regions and the finiteness of the number of chambers}
\label{sub:3-chains-bounds-b}
If the region $R(t)$ of possible parameters for semistable chains of type $t$ is bounded,
then the local finiteness of the chamber decomposition (\ref{ChamberDecomposition})
implies that there are only finitely many chambers. In this section, we will show that also
if the parameter region is unbounded, there are only finitely many chambers. This is closely
related to the question whether the set of isomorphy classes vector bundles $F$ for which
there exist a parameter $\alpha$, an $\alpha$-semistable 3-chain $C=(E_0,E_1,E_2;\phi_1,\phi_2)$
of type $t$, and an index $j_0\in \{\, 0,1,2\,\}$ with $F\cong E_{j_0}$ is bounded.
We will first give some additional bounds for the possible parameter regions, then prove the above
boundedness statement, and, last but not least, derive the finiteness of the number of chambers.
\begin{proposition}
\label{prop:unbound}
{\rm i)} Assume $r_0\neq r_1=r_2$. Then, there exist a constant $K_{12}=K_{12}(t)$
and a bounded region $R_0(t)$, such that, for any $(0,\alpha_1,\alpha_2)$-semistable
3-chain $C=(E_0,E_1,E_2;\phi_1,\phi_2)$ of type $t$, one has either that $\phi_2$
is not injective and $(\alpha_1,\alpha_2)\in R_0(t)$ or $\phi_2$ is injective and
$(\alpha_1,\alpha_2)\in R_1(t)$. Here, $R_1(t)$ is the region bounded by {\rm (\ref{41InEQ})},
{\rm (\ref{42InEQ})}, and
\begin{equation}
\label{eq:51InEQ}
r_1(\alpha_1+\alpha_2)=\alpha_1r_1+\alpha_2r_2\le K_{12}.
\end{equation}
\par
{\rm ii)} Suppose $r_0=r_1=r_2$. Then, there are constants $K_{13}$ and $K_{14}$, depending
only on the type $t$ and a bounded region $R_0(t)$ with the following property:
Given any $(0,\alpha_1,\alpha_2)$-semistable
3-chain $C=(E_0,E_1,E_2;\phi_1,\allowbreak\phi_2)$ of type $t$, then either {\rm a)} neither $\phi_1$ nor $\phi_2$
is injective and $(\alpha_1,\alpha_2)\in R_0(t)$, or {\rm b)} $\phi_1$ is not injective but $\phi_2$ is
and $(\alpha_1,\alpha_2)\in R_1(t)$, or {\rm c)} $\phi_1$ is injective, $\phi_2$ isn't, and
$(\alpha_1,\alpha_2)\in R_2(t)$, or {\rm d)} both $\phi_1$ and $\phi_2$ are injective and
$(\alpha_1,\alpha_2)\in R_3(t)$. Here, $R_1(t)$ is the region limited by the inequalities
{\rm (\ref{41InEQ})}, {\rm (\ref{42InEQ})}, and
$$
\alpha_1r_1+\alpha_2r_2\le K_{13},
$$
$R_2(t)$ is the region confined by the restrictions
{\rm (\ref{41InEQ})}, {\rm (\ref{42InEQ})}, and
$$
-\alpha_1r_1+\alpha_2(r_0+r_1)\le K_{14},
$$
and $R_3(t)$ is bounded by {\rm (\ref{41InEQ})} and {\rm (\ref{42InEQ})} (see {\rm Figure 1}).
\par
{\rm iii)} Assume $r_0<r_1>r_2$, $r_0=r_2$. Then, there is a bounded region $R_0(t)$, such that,
for any $(0,\alpha_1,\alpha_2)$-semistable
3-chain $C=(E_0,E_1,E_2;\phi_1,\phi_2)$ of type $t$, $C$ fails to be rank maximal and $(\alpha_1,\alpha_2)\in R_0(t)$
or $C$ is rank maximal and $(\alpha_1,\alpha_2)$ lies in the region $R_1(t)$ bounded by {\rm (\ref{41InEQ})},
{\rm (\ref{45InEQ})}, and {\rm (\ref{46InEQ})}.
\end{proposition}
\begin{proof}
 i) We first treat the case that $\phi_2$ is not injective.
If $\phi_2$ fails to be injective, then (\ref{48InEQ}) holds true and yields
$\alpha_2\le \alpha_1+ K_6/r_2'$.
In the case $r_0>r_1$, we have
an upper bound on $\alpha_1$, from the proof of Theorem \ref{BoundedParameterRegion}.
Thus, $\alpha_2$ is also bounded from above.
As we have seen before, $\alpha_2$ is always bounded from below. Finally, (\ref{41InEQ}) provides
a lower bound for $\alpha_1$. If $r_0<r_1$, then the arguments used in the proof of the case $r_0<r_1>r_2$
in Theorem \ref{BoundedParameterRegion}\footnote{More precisely, one checks that either
(\ref{41InEQ}), (\ref{48InEQ}), and (\ref{49InEQ}) or (\ref{41InEQ}), (\ref{46InEQ}), and (\ref{48InEQ}) apply,
because $\phi_2$ is not injective.}
show that the parameter $(\alpha_1,\alpha_2)$
lives in a bounded triangular region.
\par
Now, we assume that $\phi_2$ is injective. If, in the case $r_0>r_1$, $\phi_1$ is injective, then
we are done by (\ref{43InEQ}).
We set $K:=\ker(\phi_1\circ \phi_2)$, and define $r_2':=\rk(\ker(\phi_1\circ \phi_2))$,
$\mu_K:=\mu(K)$, and $r:=r_0+r_1+r_2$.
The semistability condition for the subchain $(0,K\cong\phi_2(K),K)$ reads
\begin{equation}
\label{eq:51bInEQ}
\frac{r_0}{2r}\cdot (\alpha_1+\alpha_2)\le \mu-\mu_K.
\end{equation}
The semistability condition for the subchain $(\im(\phi_1\circ \phi_2), E_2\cong \phi_2(E_2),\allowbreak
E_2)$
may be written in the form
\begin{equation}
\label{eq:52InEQ}
-(\alpha_1r_1+\alpha_2r_2)\Bigl(\underbrace{1-\frac{3r_2-r_2'}{r}}_{=:c}\Bigr)\ge r_2'(\mu-\mu_K)+3d_2-3r_2\mu.
\end{equation}
Combining (\ref{eq:51bInEQ}), (\ref{eq:52InEQ}), and (\ref{41InEQ}) proves the claim, provided $c$ is non negative.
Observe $rc=r_0-r_2+r_2'=r_0-r_1+r_2'$. Therefore, $c$ is positive, if $r_0>r_1$.
If $r_0<r_1$, we use the fact $r_2'\ge r_2-r_0=r_1-r_0$ to conclude.
\par
 ii) If neither $\phi_1$ nor $\phi_2$ is injective, then (\ref{48InEQ}) and (\ref{49InEQ}) hold true,
so that $(\alpha_1,\alpha_2)$ is an element of the triangular region bounded by
(\ref{41InEQ}), (\ref{48InEQ}), and (\ref{49InEQ}). If $\phi_2$ is injective but $\phi_1$ is not, one may use
the arguments from the proof of i). Likewise, one gets the result in the case where
$\phi_1$ is injective but $\phi_2$ is not, by looking at the dual chain.
\par
 iii) If $\phi_1\circ\phi_2$ is not rank maximal, one may use the same arguments as in the proof of
Theorem \ref{BoundedParameterRegion}, Case $r_0<r_1>r_2$, $r_0<r_2$, to see that $(\alpha_1,\alpha_2)$
belongs to a bounded subset of $\R^2$.
\end{proof}
\begin{theorem}
\label{thm:GlobalBoundedness}
Fix the type $t$. Then, the set of isomorphy classes of vector bundles $F$ for which
there exist a parameter $\alpha$, an $\alpha$-semistable 3-chain $C=(E_0,E_1,E_2;\phi_1,\phi_2)$
of type $t$, and an index $j_0\in \{\, 0,1,2\,\}$ with $F\cong E_{j_0}$ is bounded.
\end{theorem}
\begin{proof}
The assertion is well known for one fixed parameter $\alpha$. By the local finiteness of the
chamber decomposition (\ref{ChamberDecomposition}), the theorem is also clear, if $\alpha$
is allowed to move in a bounded region. Thus, we are left with the cases of Proposition \ref{prop:unbound}.
\par
\paragraph{\bf Case i)}
Let $C=(E_0,E_1,E_2;\phi_1,\phi_2)$ be 3-chain of type $t$ which is semistable w.r.t.\ some parameter $\alpha$.
Again, we may exclude a bounded region, so that we may assume that $\phi_2$ is generically an isomorphism.
First, we look at a subbundle $F$ of $E_0$. The semistability condition for $(F,0,0)$ gives
$$
\mu(F)\le \frac{d_0+d_1+d_2}{r_0+r_1+r_2}+\frac{\alpha_1r_1+\alpha_2r_2}{r}\stackrel{(\ref{eq:51InEQ})}{\le}
\frac{d_0+d_1+d_2}{r_0+r_1+r_2}+\frac{K_{12}}{r}=:K_{15}.
$$
Since $K_{15}$ depends only on the type, this implies that $E_0$ moves in a bounded family.
Next, let $F$ be any subbundle of $E_2$ that is contained in the kernel of $\phi_1\circ\phi_2$.
We look at the subchain $(0,F\cong\phi_2(F), F)$. This gives
$$
\mu(F)\le \frac{d_0+d_1+d_2}{r_0+r_1+r_2}+\Bigl(\frac{2r_1}{r}-\frac{1}{2}\Bigr)\cdot
(\alpha_1+\alpha_2)
\le K_{16},
$$
by (\ref{41InEQ}).
Now, let $F$ be an arbitrary subbundle of $E_2$. Then, we find the extension
$$
\begin{CD}
0 @>>> F':=F\cap \ker(\phi_1\circ \phi_2) @>>> F @>>>
\\
@>>> F'':=(\phi_1\circ \phi_2)(F) @>>> 0.
\end{CD}
$$
Thus,
$$
\mu(F)=\frac{\rk(F')\mu(F')+\rk(F'')\mu(F'')}{\rk(F)} \le \max\{\, K_{15}, K_{16}\,\}.
$$
This shows that $E_2$ lives in a bounded family, too.
Finally, $E_1$ is given as an extension
$$
\begin{CD}
0 @>>> E_2 @>>> E_1 @>>> T @>>> 0,
\end{CD}
$$
with $T$ a torsion sheaf of length $d_1-d_2$. It follows easily that $E_1$ also belongs to a bounded family.
\paragraph{\bf Case ii)}
If $\phi_2$ is injective, but $\phi_1$ is not, then we may argue as in Case i). Similarly, we obtain
the result in the case that $\phi_1$ is injective, but $\phi_2$ is not.
Finally, we look at the case where both $\phi_1$ and $\phi_2$ are injective. For any subbundle
$F$ of $E_2$, we look at the subchain $(F\cong(\phi_1\circ\phi_2)(F), F\cong\phi_2(F), F)$.
This yields the condition
$$
\mu(F)\le \mu.
$$
This proves that $E_2$ is a member of a bounded family.
For $E_0$ and $E_1$, we find the analogous result by looking at the extensions
$$
\begin{CD}
0 @>>> E_2 @>>> E_{0,1} @>>> T_{0,1} @>>> 0,
\end{CD}
$$
with $T_{0,1}$ a torsion sheaf of length $d_{0,1}-d_2$.
\paragraph{\bf Case iii)}
If $C$ is a rank maximal chain, then we get Condition (\ref{eq:54InEQ}) for every subchain
$C'=(F_0,F_1,F_2)$ of $C$. Now, suppose that $F$ is a subbundle of $E_2$ and look at the subchain
$(F\cong(\phi_1\circ\phi_2)(F), F\cong\phi_2(F), F)$. Then, (\ref{eq:54InEQ}) yields the estimate
$$
\mu(F)\le \mu+\frac{1}{3}K_{10}+\frac{1}{r}K_{11}=:K_{17},
$$
so that the family of possible $E_2$'s is bounded. Since $E_0$ is an extension of a torsion sheaf
of length $d_0-d_2$ by $E_2$, the family of possible $E_0$'s is bounded, too.
Next, if $F$ is a subbundle of $E_1$ which is contained in the kernel of $\phi_1$, then
the condition for the subchain $(0,F,0)$ gives
$$
\mu(F)\le \mu+\frac{-\alpha_1(r_0+r_2)+\alpha_2r_2}{r}\stackrel{(\ref{45InEQ})}{\le} K_{18}.
$$
An arbitrary subbundle $F$ of $E_1$ is written as an extension
$$
\begin{CD}
0 @>>> F':=F\cap \ker(\phi_1) @>>> F @>>> F'':=\phi_1(F) @>>> 0.
\end{CD}
$$
We infer
$$
\mu(F)=\frac{\rk(F')\mu(F')+\rk(F'')\mu(F'')}{\rk(F)} \le \max\{\, K_{17}, K_{18}\,\}
$$
and settle the theorem.
\end{proof}
\begin{corollary}
Fix the type $t$. Then, there are only finitely many ``effective" chambers.
\end{corollary}
\begin{proof}
By Theorem \ref{thm:GlobalBoundedness},
there is a constant $d_\infty=d_\infty(t)$, such that, for any $(\alpha_1,\alpha_2)\in\R^2$,
any $\alpha$-semistable holomorphic chain $C=(E_0,E_1,\allowbreak E_2;\allowbreak \phi_1,\allowbreak \phi_2)$
of type $t$, $\alpha:=(0,\alpha_1,\alpha_2)$, any index $j_0\in\{\, 0,1,2\,\}$,
and any subbundle $F\subseteq E_{j_0}$, one has $\deg(F)\le d_\infty$.
One easily derives the following assertion.
\begin{lemma}
\label{lem:ChamberFinite}
Fix a constant $L$. Then, there is an integer $d_l=d_l(t,\allowbreak L)$, such that
for any $(\alpha_1,\alpha_2)\in\R^2$,
any $\alpha$-semistable holomorphic chain $C=(E_0,E_1,E_2;\phi_1,\phi_2)$
of type $t$, $\alpha:=(0,\alpha_1,\alpha_2)$,
and any subchain $C'=(F_0,F_1,F_2)$, one has
$$
\frac{\deg(F_0)+\deg(F_1)+\deg(F_2)}{r_0+r_1+r_2} < L,
$$
whenever there exists an index $j_0\in \{\,0,1,2\,\}$ with $\deg(F_{j_0})< d_l$.
\end{lemma}
We go again through the cases of Proposition \ref{prop:unbound}.
\paragraph{\bf Case i)}
As usual, we may assume that $\phi_2\neq 0$ for the chains we are dealing with.
Let $C=(E_0,E_1,E_2;\phi_1,\phi_2)$ be such a chain and $C'=(F_0,F_1,F_2)$ a subchain.
Note that (\ref{eq:51InEQ}) provides an upper bound for $\alpha_1+\alpha_2$ and that
``-$r_2$(\ref{42InEQ})+$(r_0+r_1)$(\ref{eq:51InEQ})" gives an upper bound for $\alpha_1$.
Since $\rk(F_1)\ge \rk(F_2)$, it is easy to find a constant $K_{19}=K_{19}(t)$ with
\begin{eqnarray*}
&&\mu(F_0\oplus F_1\oplus F_2)+\frac{\alpha_1\rk(F_1)+\alpha_2\rk(F_2)}{\rk(F_0\oplus F_1\oplus F_2)}
\\
&=&
\mu(F_0\oplus F_1\oplus F_2)+\frac{\alpha_1(\rk(F_1)-\rk(F_2))}{\rk(F_0\oplus F_1\oplus F_2)}
+ \frac{(\alpha_1+\alpha_2)\rk(F_2)}{\rk(F_0\oplus F_1\oplus F_2)}
\\
&\le&
\mu(F_0\oplus F_1\oplus F_2)+K_{19}.
\end{eqnarray*}
On the other hand, by (\ref{41InEQ}),
$$
\frac{d_0+d_1+d_2}{r_0+r_1+r_2}+\frac{\alpha_1r_1+\alpha_2r_2}{r_0+r_1+r_2}
\ge
K_{20}.
$$
With Lemma \ref{lem:ChamberFinite}, we see that, if there is one index $j_0\in\{\,0,1,2\,\}$ with
$\mu(F_{j_0})< d_l(t, K_{20}-K_{19})$, then the condition
$$
\mu(F_0\oplus F_1\oplus F_2)+\frac{\alpha_1\rk(F_1)+\alpha_2\rk(F_2)}{\rk(F_0\oplus F_1\oplus F_2)}
<\frac{d_0+d_1+d_2}{r_0+r_1+r_2}+\frac{\alpha_1r_1+\alpha_2r_2}{r_0+r_1+r_2}
$$
is satisfied. We, therefore, define
\begin{eqnarray*}
{\mathcal S}_{\rm eff}:=\Bigl\{\,(s_0,s_1,s_2;e)&\big|& 0\le s_j\le r_j,\ j=0,1,2,
\\
&& 0<s_0+s_1+s_2< r_0+r_1+r_2,
\\
&& 3d_l(t, K_{20}-K_{19})\le e\le 3d_\infty(t)\,\Bigr\}.
\end{eqnarray*}
This is a finite set. As in Section \ref{sub:chambers-walls}, we derive a decomposition of $\R^2$
into a (now finite) set of chambers. Together with the chamber decomposition of $R_0(t)$, we thus
obtain a decomposition of $\R^2$ into a finite set of locally closed chambers, such that the property
of Proposition \ref{ChamberDecompositionII} remains true.
\paragraph{\bf Case ii)}
The arguments of Case i) apply to such chains where one map fails to be a generic isomorphism.
Hence, we are reduced to study chains $C=(E_0,E_1,E_2;\phi_1,\phi_2)$ where both $\phi_1$
and $\phi_2$ are injective. Note that, by (\ref{41InEQ}) and (\ref{42InEQ}),
there is constant $K_{21}$ with
$$
\alpha_1\le 2\alpha_2+K_{21}.
$$
We may use a simple modification of the argument given at the end of the proof of Theorem \ref{UnstableQuivers}
to see that there is a constant $K_{21}$ which depends only on the type $t$, such that
$$
\mu(F_0\oplus F_1\oplus F_2)+\frac{\alpha_1\rk(F_1)+\alpha_2\rk(F_2)}{\rk(F_0\oplus F_1\oplus F_2)}
\le \mu(F_0\oplus F_1\oplus F_2)+K_{22}.
$$
The rest of the proof proceeds as before.
\paragraph{\bf Case iii)}
We look only at chains $C=(E_0,E_1,E_2;\phi_1,\phi_2)$ where $\phi_1\circ\phi_2$ is generically an
isomorphism.
Note that (\ref{45InEQ}) and the lower bound for
$\alpha_2$ provide us with a lower bound for $\alpha_1$. Using (\ref{eq:54InEQ})
and $\rk(F_2)\le \rk(F_0)$ for every subchain $C'=(F_0,F_1,F_2)$ of $C$, we find a constant $K_{23}$ with
$$
\alpha r_2(r_2'-r_0')=\alpha_1(r_0r_2'-r_0'r_2)\le K_{23}.
$$
For an appropriate choice of $L$, (\ref{eq:54InEQ}) will be satisfied with ``$<$", if  there is
$j_0\in\{\,0,1,2\,\}$ with $\mu(F_{j_0})< d_l(t, L)$. Hence, we may continue as before.
\end{proof}
\begin{example}
In Figure 4, we sketch the shape of the parameter region together with its chamber
structure away from some bounded region
in the cases $r_0\neq r_1=r_2$ and $r_0<r_1>r_2$, $r_0=r_2$. We have also marked
an ``extremal chamber" ${\mathcal C}_\infty$. The hope is that one can understand
the corresponding ``extremal moduli space" sufficiently well to start the investigation
of other moduli spaces via birational transformations, using the results from Section \ref{sec:wall-crossing}.
\begin{figure}[h]
\begin{center}
\begin{pspicture}(-6,-.5)(6,5)
\psset{xunit=1cm, yunit=1cm}
\pspolygon*[linecolor=lblue](-6,-.5)(-6,5)(6,5)(6,-.5)
\pspolygon*[linecolor=lightorange](-5,1)(-6,2)(-6,5)(-4.333,5)(-.333,1)\pspolygon*[linecolor=lightorange](.75,1)(2.75,5)(6,5)(6,3)(5,1)
\pspolygon[linecolor=lblue, fillstyle=hlines, hatchcolor=blue](-6,0)(-6,1)(-.25,1)(-.25,0)
\pspolygon[linecolor=lblue, fillstyle=vlines, hatchcolor=blue](6,0)(6,1)(.25,1)(.25,0)
\psdot[linecolor=mix, linewidth=14pt](-3.7,3.5)
\uput[0](-4.25,3.5){\scaleboxto(.8,.6){\orange ${\mathcal C}_\infty$}}
\uput[0](-4.25,-.25){$r_0\neq r_1=r_2$}
\uput[0](1.75,-.25){$r_0=r_2<r_1$}
\uput[0](-5.925,3.5){\red$\bullet\bullet\bullet$}
\uput[0](3.9,3.5){\red$\bullet\bullet\bullet$}
\psdot[linecolor=mix, linewidth=16pt](2.85,3.5)
\uput[0](2.25,3.5){\scaleboxto(.9,.7){\orange ${\mathcal C}_\infty$}}
\psline[linecolor=red, linewidth=1.5pt](-6,5)(-2,1)
\psline[linecolor=red, linewidth=1.5pt](-4.333,5)(-.333,1)
\psline[linecolor=red, linewidth=1.5pt](-6,3)(-4,1)
\psline[linecolor=red, linewidth=1.5pt](-6,2)(-5,1)
\psline[linecolor=red, linewidth=1.5pt](.75,1)(2.75,5)
\psline[linecolor=red, linewidth=1.5pt](5,1)(6,3)
\psline[linecolor=red, linewidth=1.5pt](4,1)(6,5)
\psline[linecolor=red, linewidth=1.5pt](2.5,1)(4.5,5)
\psline[linecolor=white, linewidth=5pt](0,-.5)(0,5)
\end{pspicture}
\end{center}
\caption{Unbounded parameter regions for 3-chains.}
\end{figure}
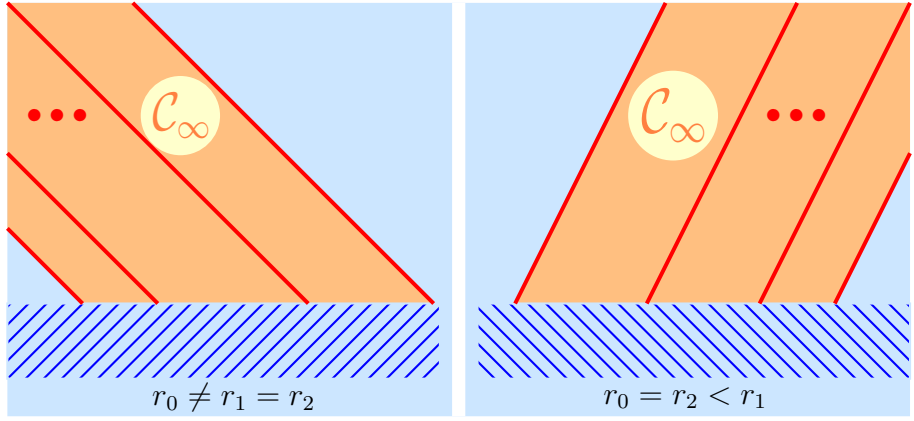
\end{example}
\subsection{Concluding remarks}
i) For any type $t$, we have found a region $R(t)\subset \R^2$, such that the existence of
an $(0,\alpha_1,\alpha_2)$-semistable chain of type $t$ implies $(\alpha_1,\alpha_2)\in R(t)$.
(Although the bounds we have found as well all the constants appearing are given by complicated
expressions, they can be explicitly determined.)
If we want that $R(t)$ has a non-empty interior, we find non-trivial restrictions on the type.
E.g., for $r_0>r_1<r_2$, we have
$$
\frac{d_0}{r_0}-\frac{d_1}{r_2}\le \alpha_2\le
\frac{2}{r_0-r_1}\Bigl(d_0-\frac{r_0}{r_1}d_1\Bigr)
+
\frac{2}{r_2-r_1}\Bigl(\frac{r_2}{r_1}d_1-d_2\Bigr).
$$
\par
ii) For augmented or decorated vector bundles, i.e., vector bundles together with a section in the
vector bundle associated by means of a homogeneous representation $\rho\colon \GL_r(\CC)\lra\GL(V)$,
the connection between the behaviour of the semistability concept for large parameters and the
invariant theory in $V$ has been understood in general in \cite{SchGlobal}.
To our knowlegde, we investigate here for the first time the analogous question for a reductive
group other than $\GL_r(\CC)$ (namely, $\GL_{r_0}(\CC)\times\GL_{r_1}(\C)\times\GL_{r_2}(\CC)$).
Our arguments are valid only for the special situation we are looking at, but in view of Theorem
\ref{UnstableQuivers}, the relationship between the shape of the region of possible stability parameters
and the invariant theory in $\Hom(\CC^{r_2},\CC^{r_1})\oplus \Hom(\CC^{r_1},\CC^{r_0})$ is clearly
perceptible. Thus, we get some feeling why Conjecture \ref{ParBounds} and some more general properties should
be true.


\section{Extremal moduli spaces for 3-chains}
\label{sec:extensions2}                          
This section serves as an illustration of the geometry of moduli
spaces for $3$-chains and its relation to other problems. We will
study a few specific types in which we have inserted ``ones". This
condition is used to grant that the chains we will consider are all
rank maximal, so that we have a good picture of the a priori
parameter region $R(t)$ and may exhibit a two-dimensional chamber
${\mathcal C}_\infty$ which yields the ``asymptotic moduli spaces".
We may expect that these moduli spaces are in a certain way the
easiest and are related to other well-known moduli spaces such as
moduli spaces of semistable vector bundles. On the other hand, we
have laid in previous chapters the foundations for studying other,
birationally equivalent moduli spaces via the ``flip-technology".
Although we discuss only very special types, it becomes clear how one
may in general relate the moduli spaces in the extremal chamber
${\mathcal C}_\infty$ to the moduli spaces on the nearby boundaries
that one may usually understand quite easily (compare Propositions
\ref{StabII} and \ref{prop:StabIV} below). This should give the
reader sufficient material to attack any special case she or he is
interested in.
\subsection{Generalities on moduli spaces for type $(m, 1, n;d_0,d_1,d_2)$}
\label{sub:large-moduli-n=1}
In the next two sections, we will describe the moduli spaces ${\mathcal M}_\alpha(t)$ for the type 
$t=(m, 1, n;d_0,d_1,d_2)$ with respect to the stability parameters $\alpha$ which lie in a certain 
``extremal" two-dimensional chamber. 
To this end, we first recall the results concerning the parameter region and the moduli spaces that 
we have already obtained.
\par
Let ${\mathcal C}$ be a two dimensional chamber. All the chains that we will have to consider will be automatically 
rank maximal in the sense of Section \ref{sub:Rank-Maximal}, because of Remark \ref{MapsAreNonZero}. 
Therefore, we only have to look at the relevant inequalities from Proposition \ref{prop:RankMaxChains} that bound the parameter region $R(t)$.
Let us remind the reader what these inequalities are and how they are obtained.
\par
\it Inequality \rm I.
Obviously, $(\im(\phi_1),E_1,E_2)$ is a subchain. Note that $\im(\phi_1)$ is isomorphic to $E_1$, so that it has degree $d_1$ and rank $1$. The condition of $\alpha$-semistability for this subchain reads
$$
(m-1)\alpha_1+(m-1)n\alpha_2\le A_{\rm I}:=(n+2)d_0-(2m+n)d_1+(1-m)d_2.
$$
\par
\it Inequality \rm II.
The condition of $\alpha$-semistability for the subchain $(E_0,\allowbreak 0,0)$ produces the inequality
$$
m\alpha_1+mn\alpha_2\ge A_{\rm II}:=(n+1)d_0-md_1-md_2.
$$
\par
\it Inequality \rm III.
Here, one checks $\alpha$-semistability for the subchain $(E_0,\allowbreak E_1,0)$.
This gives
$$
-n\alpha_1+(m+1)n\alpha_2\ge A_{\rm III}:= nd_0+nd_1-(m+1)d_2.
$$
\par
\it Inequality \rm IV. 
This inequality only applies, if $n>1$.
One uses the subchain $(0,0,\ker(\phi_2))$. Clearly, $\ker(\phi_2)$ has rank $n-1$ and degree $\deg(E_2)-\deg(\im(\phi_2))$. Since $\im(\phi_2)$ is a non-trivial subsheaf of the line bundle $E_1$, $\deg(\im(\phi_2))\le d_1$, i.e., $\deg(\ker(\phi_2))\ge d_2-d_1$. Thus, the condition of $\alpha$-semistability for the given subchain implies the necessary condition
$$
-(n-1)\alpha_1+(m+1)(n-1)\alpha_2\le A_{\rm IV}:=(n-1)d_0+(m+2n)d_1-(m+2)d_2.
$$
\begin{remark}
\label{rem:Parallelo}
i) If $n>1$, then the Inequalities I-IV bound a parallelogram. One checks that potential destabilizing objects are of the form $(F_0,0,0)$, $(0,0,F_2)$, and $(F_0,E_1,F_2)$ (we do not have to consider subchains of the form $(F_0,0,F_2)$, because the condition of $\alpha$-(semi)stability for such a subchain follows from those for the subchains
$(0,0,F_2)$ and $(F_0,0,0)$). One checks that the corresponding one-dimensional walls are parallel to one of the sides of the parallelogram.
\par
ii) For $n=1$, the inequalities I-III bound an ``open
parallelogram". We claim that all the one-dimensional walls are
defined by an equation of the form $\alpha_1+\alpha_2=c$.
Therefore, the parameter region with its chamber structure looks
as depicted in Figure 5. In fact, if
$C:=(E_0,E_1,E_2;\phi_1,\phi_2)$ is an $\alpha$-semistable but not
$\alpha$-stable chain, then we find an $\alpha$-destabilizing
subchain $(F_0,F_1,F_2)$. Setting
$C':=(F_0,F_1,F_2;\allowbreak\phi_{1|F_1},\phi_{2|F_2})$ and
$C'':=C/C'$, the chain $\widetilde{C}:=C'\oplus C''$ is still
$\alpha$-semistable. If we assume that $\alpha$ lies in the \sl
interior \rm of $R(t)$, then none of the homomorphisms in
$\widetilde{C}$ must be zero. If we replace $C$ by
$\widetilde{C}$, then we easily see that the wall containing
$\alpha$ is defined via a subchain $(0,0,F_2)$. From this, one
immediately arrives at our claim.
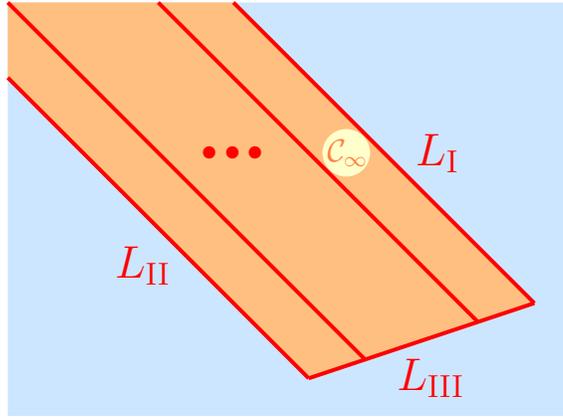
\begin{figure}[h]
\begin{center}
\begin{pspicture}(-6,-.5)(1.5,5)
\psset{xunit=1cm, yunit=1cm}
\pspolygon*[linecolor=lblue](-6,-.5)(-6,5)(1.5,5)(1.5,-.5)
\pspolygon*[linecolor=lightorange](-2,0)(-6,4)(-6,5)(-3,5)(1,1)
\psline[linecolor=red, linewidth=1.5pt](-2,0)(1,1)
\psline[linecolor=red, linewidth=1.5pt](-2,0)(-6,4)
\psline[linecolor=red, linewidth=1.5pt](-3,5)(1,1)
\psline[linecolor=red, linewidth=1.5pt](-1.25,.25)(-6,5)
\psline[linecolor=red, linewidth=1.5pt](.25,.75)(-4,5)
\uput[0](-3.6,3){\red$\bullet\bullet\bullet$}
\psdot[linecolor=mix, linewidth=8pt](-1.5,3)
\uput[0](-1.925,3){\orange${\mathcal C}_\infty$}
\uput[0](-.75,3){\Large\red${L}_{\rm I}$}
\uput[0](-4.75,1.5){\Large\red${L}_{\rm II}$}
\uput[0](-1,0){\Large\red${L}_{\rm III}$}
\end{pspicture}
\end{center}
\caption{The chamber structure for 3-chains of type $(m,1,1;d_0,d_1,d_2)$.}
\end{figure}
\end{remark}
Finally, we also add the following useful observation.
\begin{lemma}
\label{StabI}
{\rm i)} Suppose $n=1$. Then, $R(t)$ has a non-empty interior if and only
$$
d_1< \frac{d_0}{m},
$$
that is, the subsheaf $\im(\phi_1)$ does not destabilize $E_0$.
\par
{\rm ii)} If $n>1$, then the parameter region has a non-empty interior if and only if the condition in {\rm i)} holds and additionally
$$
(m+1)(nd_1-d_2)- md_2>0.
$$
\end{lemma}
\begin{proof}
i) The interior of $R(t)$ is obviously non-empty if and only if $m\cdot A_{\rm I}> (m-1)\cdot A_{\rm II}$ in the above notation. A few simplifications lead to the assertion.
\par
ii) The second condition arises from evaluating the inequality $n\cdot A_{\rm IV}> (n-1)\cdot A_{\rm III}$.
\end{proof}
\begin{remark}
\label{rem:StabII}
i) We note that the inequality in Lemma \ref{StabI} implies that $d_1<(d_0-d_1)/(m-1)$.
\par
ii) It will later be interesting to know that, in Lemma \ref{StabI}, we may choose $d_1$ and $d_2$ in such a way that the inequality holds and $(nd_1-d_2)$ is a prescribed value (which might also be negative).
\end{remark}
First, let us analyze the moduli spaces for parameters which do lie on the boundary of the parameter region. If the stability parameter $\alpha$ lies on one of the boundaries determined by Inequality II or III, then we have described the moduli space in Corollary \ref{cor:Start}. Thus, one of the remaining cases is the one when $(m-1)\alpha_1+(m-1)n\alpha_2=(n+2)d_0-(2m+n)d_1+(1-m)d_2$. We assume that the remaining Inequalities II, III, and IV are strict, so that both $\phi_2\not\equiv 0$ and $\phi_1\not\equiv 0$. Suppose $(E_0,E_1,E_2;\phi_1,\phi_2)$ is $\alpha$-semistable. By the definition of the boundaries, the subchain $(\im(\phi_1),E_1,E_2)$ becomes destabilizing. Standard arguments now show that
\begin{itemize}
\item $(\im(\phi_1)\cong E_1,E_1,E_2)$ is an $\alpha$-semistable
holomorphic $3$-chain of type $t^\p:=(m,1,1;d_0,d_1,d_2)$,
\item $Q_0:=E_0/\im(\phi_1)$ is a semistable vector bundle of degree
$d_0-d_1$ and rank $m-1$, and
\item $(E_0,E_1,E_2;\phi_1,\phi_2)$ is S-equivalent to the chain
$(\im(\phi_1)\oplus Q_0,\allowbreak E_1,E_2;\phi_1,\phi_2)$.
\end{itemize}
We infer.
\begin{proposition}
\label{StabII}
The natural morphism
\begin{eqnarray*}
\sigma_{\rm I}\colon {\mathcal M}_\alpha(t)&\lra&
{\mathcal M}_\alpha(t^\p) \times {\mathcal U}(n-1,d_0-d_1)
\\
\ [E_0,E_1,E_2; \phi_1,\phi_2]&\lma&\bigl([\im(\phi_1),E_1,E_2;\phi_1,\phi_2],[E_0/\im(\phi_1)]\bigr)
\end{eqnarray*}
is bijective, and there is also the inverse morphism
\begin{eqnarray*}
\tau_{\rm I}\colon {\mathcal M}_\alpha(t^\p)\times {\mathcal U}(n-1,d_0-d_1)
&\lra& {\mathcal M}_\alpha(t)
\\
\bigl([E_1,E_1,E_2;\phi_1,\phi_2],[Q_0]\bigr)&\lma&
[E_1\oplus Q_0,E_1,E_2; \phi_1,\phi_2].
\end{eqnarray*}
\end{proposition}
The other case to consider is the one when $-(n-1)\alpha_1+(m+1)(n-1)\alpha_2=(n-1)d_0+(m+2n)d_1-(m+2)d_2$.
We may assume that the remaining Inequalities I, II, and III are strict, i.e., $\phi_2\not\equiv 0$ and $\phi_1\not\equiv 0$.
\begin{proposition}
\label{prop:StabIV}
There are the bijective morphisms
\begin{eqnarray*}
\sigma_{\rm IV}\colon {\mathcal M}_\alpha({t})&\lra& {\mathcal U}(n-1,d_2-d_1)
\times {\mathcal M}_\alpha(m,1,1;d_0,d_1,d_1)
\\
\ [E_0,E_1,E_2; \phi_1,\phi_2]&\lma&\bigl([\ker(\phi_2)],
[E_0,E_1,E_1\cong E_2/\ker(\phi_2);\phi_1,\ol{\phi}_2]\bigr)
\end{eqnarray*}
and
$$
\begin{array}{c}
\tau_{\rm IV}\colon {\mathcal U}(n-1,d_2-d_1)
\times {\mathcal M}_\alpha(m,1,1;d_0,d_1,d_1)\ \lra\ {\mathcal M}_\alpha(t)
\\
\bigl([K_2],[E_0,E_1,E_1;\phi_1,\phi_2]\bigr)\ \lma\
[E_0,E_1,E_1\oplus K_2; \phi_1,\phi_2].
\end{array}
$$
\end{proposition}
\subsection{Extremal moduli spaces for type $(m, 1, 1;d_0,d_1,d_2)$}
Now, suppose $\alpha=(0,\alpha_1,\alpha_2)$ is such that $(\alpha_1,\alpha_2)$ lies in the chamber 
${\mathcal C}_\infty$ (i.e., the chamber in the interior of $R(t)$ which is adjacent to the line $L_{\rm I}$ 
on which Inequality I becomes an equality; see Figure 5). 
Let $\alpha^M=(0,\alpha_1^M,\alpha_2^M)$ be such that $(\alpha_1^M,\alpha_2^M)$ is an element of 
$\ol{\mathcal C}_\infty\cap L_{\rm I}$. 
Note that $\alpha_1+\alpha_2<\alpha_1^M+\alpha_2^M$, so that 
$\mu_\alpha(t)<\mu_{\alpha^M}(t)$. Let $(E_0,E_1,E_2;\phi_1,\phi_2)$ be an $\alpha$-semistable holomorphic chain of type $t=(m,1,1;d_0,d_1,d_2)$. We note the following properties.
\begin{proposition}
\label{bound102}
{\rm i)} The vector bundle $E_0$ does not possess a subbundle of slope $(d_0-d_1)/(m-1)$ or higher and is given by a non-split extension
$$
\begin{CD}
0 @>>> E_1 @>>> E_0 @>>> Q_0 @>>> 0
\end{CD}
$$
of a semistable vector bundle $Q_0$ of degree $d_0-d_1$ and rank $m-1$ by $E_1$.
\par
{\rm ii)} In {\rm i)}, we have
$$
\dim_\C\bigl({\rm Ext}^1(Q_0,E_1)\bigr)=d_0-md_1+(m-1)(g-1).
$$
\end{proposition}
\begin{proof}
i) By Proposition \ref{ChamberDecompositionII}, we know that $(E_0,E_1,E_2;\phi_1,\phi_2)$ is also $\alpha^M$-semi\-sta\-ble. By the results stated before Proposition \ref{StabII}, this implies that $\im(\phi_1)\cong E_1$ is a subbundle of $E_0$ and that $Q_0:= E_0/\im(\phi_1)$ is a semistable vector bundle of degree $d_0-d_1$ and rank $m-1$. Observe
\begin{equation}
\label{bound101}
\frac{d_0-d_1}{m-1}=\mu(Q_0)=\mu_{\alpha^M}(E_0,E_1,E_2;\phi_1,\phi_2).
\end{equation}
For a subbundle $F_0$ of $E_0$, we thus obtain
$$
\mu(F_0)\le \mu_{\alpha}(E_0,E_1,E_2;\phi_1,\phi_2)<\mu_{\alpha^M}(E_0,E_1,E_2;\phi_1,\phi_2)=\frac{d_0-d_1}{m-1}.
$$
This also implies that the extension is non-split.
\par
ii) Recall from Lemma \ref{StabI} that $d_1<\mu(E_0)$, whence $\mu(Q_0)>\mu(E_0)$.
Since $Q_0$ is semistable, this implies
$$
H^0(Q_0^\vee\otimes E_1)={\rm Hom}(Q_0,E_1)=\{0\}.
$$
Since ${\rm Ext}^1(Q_0,E_1)=H^1(Q_0^\vee\otimes E_1)$, the given formula is a consequence of the Riemann-Roch theorem.
\end{proof}
There is also a partial converse to Proposition \ref{bound102}.

\begin{proposition}
\label{StabIII}
Let $\alpha=(0,\alpha_1,\alpha_2)$ be a stability parameter with $(\alpha_1,\alpha_2)\in {\mathcal C}_\infty$. Then, 
a holomorphic chain $C=(E_0,E_1,E_2;\phi_1,\phi_2)$ with $E_1$ a line bundle of degree $d_1$, $D$ an effective divisor of degree $d_1-d_2$, $Q_0$ a {\bfseries stable} 
vector bundle of degree $d_0-d_1$ and rank $m-1$,
$$
\begin{CD}
0 @>>> E_1 @>\phi_1 >> E_0 @>>> Q_0 @>>> 0
\end{CD}
$$
a non-split extension, $E_2:=E_1(-D)$, and $\phi_2\colon E_2\subseteq E_1$ is $\alpha$-stable.
\end{proposition}
\begin{proof}
For any non-trivial subbundle $F_0\subsetneq E_0$, we have to check the stability condition for the subchain $(F_0,0,0)$,
and, if $E_1\subseteq F_0$, also for the subchains $(F_0,E_1,0)$ and $(F_0,E_1,E_2)$. 
\par
In the following, let $\alpha^M=(0,\alpha_1^M,\alpha_2^M)$ be such that $(\alpha_1^M,\alpha_2^M)$ is a point on $L_{\rm I}$ which lies in the interior of the region depicted
in Figure 1. Let $F_0$ be a subbundle of $E_0$. 
If $F_0\cap E_1=\{0\}$, we find
$$
\mu(F_0)< \mu(Q_0)=\frac{d_0-d_1}{m-1}=\mu_{\alpha^M}(E_0,E_1,E_2;\phi_1,\phi_2),
$$
because the extension is non-split. Otherwise, $E_1\subseteq F_0$, and we have the exact sequence
\begin{equation}
\label{Ext007}
\begin{CD}
0 @>>> E_1 @>>> F_0 @>>> F_0/E_1 @>>> 0
\end{CD}
\end{equation}
where
$$
\mu(F_0/E_1)\le \mu(Q_0)=\mu_{\alpha^M}(E_0,E_1,E_2;\phi_1,\phi_2).
$$
On the other hand,
$$
\mu(E_1)=d_1\stackrel{\rm Lemma\ \ref{StabI}}<\mu(E_0)<\mu_{\alpha^M}(E_0,E_1,E_2;\phi_1,\phi_2).
$$
These two facts imply again
$$
\mu(F_0)< \mu_{\alpha^M}(E_0,E_1,E_2;\phi_1,\phi_2).
$$
This strict inequality still holds for all $(\alpha_1,\alpha_2)$ in the interior of the parameter region that are close enough to $(\alpha_1^M,\alpha_2^M)$.
By definition of the chamber decomposition, it must then hold for all $(\alpha_1,\alpha_2)\in {\mathcal C}_\infty$.
\par
Next, let us look at a subchain of the type $(F_0,E_1,0)$, $F_0$ a subbundle of $E_0$. In this case, $F_0$ must contain $\im(\phi_1)$, 
so that we have again the extension (\ref{Ext007}). We claim that 
$$
\mu_{\alpha^M}(\im(\phi_1),E_1,0;\phi_1,0)< \mu_{\alpha^M}(E_0,E_1,E_2;\phi_1,\phi_2).
$$
Since $(\alpha_1^M,\alpha_2^M)$ lies above the line $L_{\rm III}$, we have 
\begin{equation}
\label{eq:WellKnown}
\mu_{\alpha^M}(E_0,E_1,0;\phi_1,0)< \mu_{\alpha^M}(E_0,E_1,E_2;\phi_1,\phi_2).
\end{equation}
The subchain $\ol C=(E_0,E_1,0;\phi_1,0)$ is an extension of the chain $\ol C'':=(Q_0,0,0;0,0)$ by the chain $\ol C':=(\im(\phi_1),E_1,0;\phi_1,0)$.
Hence,
$$
\mu_{\alpha^M}(\ol C)=\frac{\rk(\ol C')\mu_{\alpha^M}(\ol C')+\rk(\ol C'')\mu_{\alpha^M}(\ol C'')}{\rk(\ol C)}.
$$
Since $\mu_{\alpha^M}(\ol C'')=\mu_{\alpha^M}(C)$, our contention follows from (\ref{eq:WellKnown}).
By the stability of $Q_0$, we also have
$$
\mu(F_0/E_1) \le \mu(Q_0)=\mu_{\alpha^M}(E_0,E_1,E_2;\phi_1,\phi_2).
$$
Therefore,
\begin{eqnarray*}
&&
\mu_{\alpha^M}(F_0,E_1,0;\phi_1,0)
\\
&=&
\frac{d_1+d_1+\deg(F_0/E_1)+\alpha_1^M}{1+\rk(F_0)}
\\
&=& \frac{2 \mu_{\alpha^M}(\im(\phi_1),E_1,0;\phi_1,0)+(\rk(F_0)-1)\mu(F_0/E_1)}{1+\rk(F_0)}
\\
&<& \frac{2 \mu_{\alpha^M}(E_0,E_1,E_2;\phi_1,\phi_2)+(\rk(F_0)-1)
\mu_{\alpha^M}(E_0,E_1,E_2;\phi_1,\phi_2)}{1+\rk(F_0)}
\\
&=& \mu_{\alpha^M}(C).
\end{eqnarray*}
Again, we see that the same inequality holds for stability parameters $\alpha=(0,\alpha_1,\alpha_2)$ with $(\alpha_1,\alpha_2)\in {\mathcal C}_\infty$.
\par
For proper subchains $(F_0,E_1,E_2)$, we obtain the chain $(\im(\phi_1), E_1,\allowbreak  E_2;\allowbreak \phi_1,\phi_2)$ with 
$$
\mu_{\alpha^M}(\im(\phi_1), E_1, E_2;\phi_1,\phi_2)=\mu_{\alpha^M}(C)
$$
and the proper subbundle $F_0/E_1$ of $Q_0$ for which we have
$$
\mu(F_0/E_1)<\mu(Q_0),
$$
by the stability of $Q_0$. This enables us to conclude as before.
\end{proof}
Recall the following.
\begin{itemize}
\item For $(0,\alpha_1,\alpha_2)$ with $(\alpha_1,\alpha_2)\in
{\mathcal C}_\infty$, a holomorphic chain $(E_0,E_1,\allowbreak
E_2;\allowbreak\phi_1,\phi_2)$ of type ${t}$ is
$(0,\alpha_1,\alpha_2)$-semistable if and only if it is
$(0,\alpha_1,\alpha_2)$-stable. (This follows from Proposition
\ref{CoPrimalityI}.)
\item For $(\alpha^i_1,\alpha^i_2)\in
{\mathcal C}_\infty$, $i=1,2$, a holomorphic chain
$(E_0,E_1,E_2;\allowbreak\phi_1,\allowbreak\phi_2)$ of type $t$ is
$(0,\alpha^1_1,\alpha^1_2)$-semistable if and only if it is
$(0,\alpha^2_1,\alpha^2_2)$-semistable.
\end{itemize}
We may now describe the moduli spaces which belong to a stability
parameter in the chamber ${\mathcal C}_\infty$.
\begin{corollary}
\label{cor:DescMod1}
The moduli space ${\mathcal M}_\alpha(t)$ for
$\alpha=(0,\alpha_1,\alpha_2)$ and
$(\alpha_1,\allowbreak\alpha_2)\in {\mathcal C}_\infty$ is a
connected smooth projective variety of dimension
$$
d_0-(m-1)d_1-d_2+(m-1)m\bigl(g-1\bigr)+g.
$$
It is birationally equivalent to a ${\mathbb P}^N$-bundle over the product ${\mathcal J}^{d_1}\times X^{(d_1-d_2)}\times {\mathcal U}^{\rm s}(m-1,d_0-d_1)$
of the Jacobian of degree $d_1$ line bundles, the $(d_1-d_2)$-fold symmetric product of the curve, and the moduli space of stable vector bundles of rank $(m-1)$ and degree $(d_0-d_1)$, $N:=d_0-md_1+(m-1)(g-1)-1$.
\end{corollary}
\begin{proof}
The only thing that we have to prove is the irreducibility. The smoothness results from the fact that all $\alpha$-semistable 3-chains of type $t$ are $\alpha$-stable and Theorem \ref{thm:smoothdim}, v). The assertions about the dimension and the birational model are evident from Propositions \ref{bound102} and \ref{StabIII}. (Note that Theorem \ref{thm:smoothdim} also gives the dimension.) It suffices to exhibit an irreducible parameter space for all $\alpha$-semistable objects. The product ${\mathcal J}^{d_1}\times X^{(d_1-d_2)}$ parameterizes the pairs $(\phi_2\colon E_1(-D)\subseteq E_1)$. Moreover, it is well know that one can construct an irreducible variety ${\mathcal A}$ and a family $Q_{\mathcal A}$ on ${\mathcal A}\times X$ which contains any semistable vector bundle of rank $(m-1)$ and degree $(d_0-d_1)$. Using the theory of universal extensions \cite{Lange}, we may construct an affine bundle ${\mathcal B}$ over ${\mathcal J}^{d_1}\times {\mathcal A}$ and a vector bundle $E_{\mathcal B}$ that consists of all vector bundles which are extensions of a vector bundle $Q$ corresponding to a point $a\in {\mathcal A}$ by a line bundle of degree $d_1$. Thus, ${\mathcal D}:=X^{(d_1-d_2)}\times {\mathcal B}$ is an irreducible variety which carries a universal family of chains, such that any $\alpha$-semistable chain belongs to that family. Since $\alpha$-semistability is an open
condition, there is an open subvariety ${\mathcal D}^0$ that parameterizes exactly the $\alpha$-semistable chains. The irreducible variety ${\mathcal D}^0$ surjects onto
the moduli space ${\mathcal M}_\alpha(t)$.
\end{proof}
\begin{proposition}
\label{prop:BirRegion1}
Let $\alpha^i=(0,\alpha^i_1,\alpha^i_2)$, $i=1,2$, be two stability parameters, such that $(\alpha^i_1,\alpha^i_2)\in R(t)\cap R_{2g-2}$, $i=1,2$, that do not lie on any wall. Then, ${\mathcal M}_{\alpha^1}(t)$ and ${\mathcal M}_{\alpha^2}(t)$ are birationally equivalent smooth projective varieties of dimension $d_0-(m-1)d_1-d_2+(m-1)m(g-1)+g$ or empty.
\par
In particular, if $R_{2g-2}\cap {\rm Interior}(R(t))$ is non-empty, then these varieties are birationally equivalent 
to a ${\mathbb P}^N$-bundle over the product 
${\mathcal J}^{d_1}\times X^{(d_1-d_2)}\times {\mathcal U}^{\rm s}(m-1,d_0-d_1)$, $N:=d_0-md_1+(m-1)(g-1)-1$.
\end{proposition}
\begin{proof}
If $\alpha$ does not lie on any wall, then the notions of $\alpha$-stability and $\alpha$-semistability are equivalent,
by Corollary \ref{CoPrimalityII}. Therefore, Theorem \ref{thm:smoothdim} grants the smoothness of the moduli spaces
and determines the dimension. We will first study the region
${\rm Interior}(R(t))\cap \widetilde{R}(t)$.
\par
Let $\alpha^0=(0,\alpha_1^0,\alpha_2^0)$ be a parameter where
$(\alpha_1^0,\alpha_2^0)$ lies on a wall in the interior of $R(t)$
and $C':=(E_0',E_1',E_2';\phi_1',\phi_2')$ and $C'':=(E_0'',E_1'',E_2'';\allowbreak \phi_1'',\phi_2'')$ 
two $\alpha^0$-semistable 3-chains of type $t'$ and
$t''$, respectively, such that $t=t'+t''$ and the map $b$ defined
in (\ref{eq:extension-complex}) is an isomorphism. Since
$\widetilde{C}=C'\oplus C''$ will be an $\alpha^0$-semistable
chain and $\alpha^0$ lies in the interior of $R(t)$, no map in
$\widetilde{C}$ must be zero, by Remark \ref{MapsAreNonZero}.
We write $t'=(r_0',r_1',r_2';d_0',d_1',d_2')$ and $t''=(r_0'',r_1'',r_2'';d_0'',d_1'',d_2'')$.
Since we require the maps in $\widetilde{C}$ to be non-zero, there are two
possibilities: a) $(r_0',r_1',r_2')=(m',1,1)$ and $(r_0'',r_1'',r_2'')=(m'',0,0)$ or
b) $(r_0',r_1',r_2')=(m',0,0)$ and $(r_0'',r_1'',r_2'')=(m'',1,1)$. For a general point
$x\in X$, we let $V'$ and $V''$ be the restrictions of $C'$ and $C''$, respectively,
to $\{x\}$. These are $\C$-linear chains. With (\ref{eq:chi(V'',V')}), we compute
$$
\chi(V'',V')=m'm''
$$ 
in Case a). If $b$ is an isomorphism, we must have $\chi(V'', V')=0$, by Lemma \ref{lem:cohomology-linear-chains}. 
Since $m'$ and $m''$ are both non-zero, this is impossible.
\par
In Case b), we compute $\chi(V'',V')=m'm''-m'$. For this quantity to become zero, we must have $m''=1$
and $m'=m-1$. Thus, $V'$ has rank-type $(m-1,0,0)$ and $V''$ has rank type $(1,1,1)$. Note that these
types are \sl not excluded by Theorem \ref{thm:BirRegionFund}\rm. It remains to compute the dimension of the resulting flip
loci. Thus, let $C'=(E'_0,0,0;0,0)$ be an $\alpha^0$-semistable chain of type $(m-1,0,0;d_0',0,0)$ (which
means that $E'_0$ is a semistable vector bundle of rank $m-1$ and degree $d_0'$) and 
$C''=(E''_0, E_1, E_2; \phi_1'',\phi_2)$ an $\alpha^0$-semistable holomorphic chain of type $(1,1,1; d_0'',d_1,d_2)$.
We have to compute $\dim_\C({\rm Ext}^1(C'',C'))$. Note that ${\rm Hom}(C'',C')=\{0\}$ and that
${\mathbb H}^2(C'',C')=\{0\}$, by Proposition \ref{prop:h2-vanishing}. Therefore, Proposition \ref{prop:chi(C'',C')}
gives
$$
\dim_\C\bigl({\rm Ext}^1(C'',C')\bigr)=(m-1)(d_0''-d_1).
$$
Thus, the space of isomorphy classes of $\alpha^0$-semistable chains $C$ which are non-split extensions of a
chain $C''$ by a chain $C'$ as above has dimension
\begin{equation}
\label{eq:DimEst1}
(m-1)^2\bigl(g-1)+g+md_0''+(1-m)d_1-d_2.
\end{equation}
By assumption, we have
$$
\mu_{\alpha^0}(C'')=\mu_{\alpha^0}(C')
$$ 
that is
\begin{equation}
\label{eq:WallDefine}
(m-1)(\alpha^0_1+\alpha^0_2)=3d_0-(m+2)d_0''+(1-m)d_1+(1-m)d_2.
\end{equation}
Since $(\alpha_1^0,\alpha_2^0)$ is supposed to lie in the interior of the region $R(t)$, we find
$$
\frac{m}{m-1}\bigl(3d_0-(m+2)d_0''+(1-m)d_1+(1-m)d_2\bigr)>2d_0-md_1-md_2
$$
which amounts to 
$$
md_0''< d_0.
$$
Plugging this information into (\ref{eq:DimEst1}) and using the dimension formula Theorem \ref{thm:smoothdim}, iv)
(which applies for the same reasons as before), we see that the locus of chains $C$ which are extensions of the
type we have considered has codimension at least $g$ in the relevant moduli spaces, whence we may forget
about these extensions.
\par
The arguments given in Proposition \ref{codim} deal with the remaining cases (in which the map $b$ cannot be an
isomorphism).
\end{proof}
\begin{remark}
\label{rem:WallDefine2}
The minimal possible value for $d_0''$ is $d_1$. In that case, (\ref{eq:WallDefine}) describes exactly
the line $L_{\rm I}$. 
\end{remark}
\subsection{Extremal moduli spaces for type $(m, 1, n;d_0,d_1,d_2)$ where $n>1$}
Under this assumption, the chamber structure looks like the one depicted in Figure 6 (compare Remark \ref{rem:Parallelo}). 
In this section, we set out to describe the moduli spaces for stability parameters $\alpha=(0,\alpha_1,\alpha_2)$, such that
$(\alpha_1,\alpha_2)\in {\mathcal C}:={\mathcal C}_\infty$. It turns out that the geometry of these moduli spaces
is closely related to the geometry of Brill-Noether loci. 
\par
As before, we first describe the necessary conditions that have to be fulfilled by $\alpha$-semistable objects. To this end,
we fix stability parameters $\alpha^{M_{\rm I}}=(0,\alpha^{M_{\rm I}}_1,\alpha^{M_{\rm I}}_2)$ and $\alpha^{M_{\rm IV}}
=(0,\alpha^{M_{\rm IV}}_1, \alpha^{M_{\rm IV}}_2)$. The assumption is that $(\alpha^{M_{\rm I}}_1,\alpha^{M_{\rm I}}_2)$
is contained in the relative interior of $\ol{\mathcal C}\cap L_{\rm I}$ and that 
$(\alpha^{M_{\rm IV}}_1,\alpha^{M_{\rm IV}}_2)$ is contained in the relative interior of 
$\ol{\mathcal C}\cap L_{\rm IV}$. 
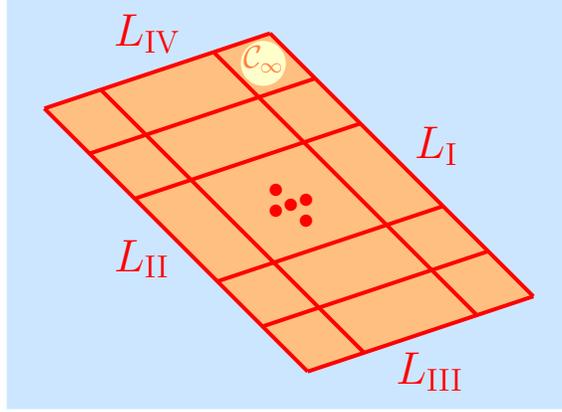
\begin{figure}[h]
\begin{center}
\begin{pspicture}(-6,-.5)(1.5,5)
\psset{xunit=1cm, yunit=1cm}
\pspolygon*[linecolor=lblue](-6,-.5)(-6,5)(1.5,5)(1.5,-.5)
\pspolygon*[linecolor=lightorange](-2,0)(-5.5,3.5)(-2.5,4.5)(1,1)
\psline[linecolor=red, linewidth=1.5pt](-2,0)(1,1)
\psline[linecolor=red, linewidth=1.5pt](-2,0)(-5.5,3.5)
\psline[linecolor=red, linewidth=1.5pt](-5.5,3.5)(-2.5,4.5)
\psline[linecolor=red, linewidth=1.5pt](-4.9,2.9)(-1.9,3.9)
\psline[linecolor=red, linewidth=1.5pt](-4.3,2.3)(-1.3,3.3)
\psline[linecolor=red, linewidth=1.5pt](-3.2,1.2)(-0.2,2.2)
\psline[linecolor=red, linewidth=1.5pt](-2.6,.6)(0.4,1.6)
\psline[linecolor=red, linewidth=1.5pt](-2.5,4.5)(1,1)
\psline[linecolor=red, linewidth=1.5pt](-1.25,.25)(-4.75,3.75)
\psline[linecolor=red, linewidth=1.5pt](.25,.75)(-3.25,4.25)
\uput[0](-2.5,2.2){\red$\bullet$}
\uput[0](-2.7,2.125){\red$\bullet$}
\uput[0](-2.3,2.275){\red$\bullet$}
\uput[0](-2.7,2.4){\red$\bullet$}
\uput[0](-2.3,2){\red$\bullet$}
\psdot[linecolor=mix, linewidth=7.5pt](-2.59,4.1)
\uput[0](-3.025,4.15){\orange${\mathcal C}_\infty$}
\uput[0](-.75,3){\Large\red${L}_{\rm I}$}
\uput[0](-4.75,1.5){\Large\red${L}_{\rm II}$}
\uput[0](-1,0){\Large\red${L}_{\rm III}$}
\uput[0](-4.75,4.5){\Large\red${L}_{\rm IV}$}
\end{pspicture}
\end{center}
\caption{The chamber structure for 3-chains of type $(m,1,n;d_0,d_1,d_2)$.}
\end{figure}
\begin{proposition}
\label{bound103}
{\rm i)} The vector bundle $E_0$ does not possess a subbundle of slope $(d_0-d_1)/(m-1)$ or higher and is given by a non-split extension
$$
\begin{CD}
0 @>>> E_1 @>>> E_0 @>>> Q_0 @>>> 0
\end{CD}
$$
of a semistable vector bundle $Q_0$ of degree $d_0-d_1$ and rank $m-1$ by $E_1$.
\par
{\rm ii)} In {\rm i)}, we have
$$
\dim_\C\bigl({\rm Ext}^1(Q_0,E_1)\bigr)=d_0-md_1+(m-1)(g-1).
$$
\par
{\rm iii)} Set $K_2:=\ker(\phi_2)$. Then, $K_2$ is a semistable vector bundle of rank $n-1$ and degree $d_2-d_1$.
The vector bundle $E_2$ is given as a non-split extension
$$
\begin{CD}
0 @>>> K_2 @>>> E_2 @>>> E_1 @>>> 0.
\end{CD}
$$
\end{proposition}
\begin{proof}
Part i) and ii) are checked in the same fashion as their counterparts in Proposition \ref{bound102} (i.e., by
using that an $\alpha$-semistable chain is also $\alpha^{M_{\rm I}}$-semistable).
\par
In order to establish iii), we use that an $\alpha$-semistable chain is $\alpha^{M_{\rm IV}}$-semistable, too.
Then, the assertion becomes a straightforward consequence of Proposition \ref{prop:StabIV}.
\end{proof}
\begin{remark}
\label{rem:VarExt}
We emphasize that we do not have a formula for computing $\dim({\rm
Ext}^1(E_1,K_2))=h^0(B_2)$, $B_2:=K_2^\vee\otimes
E_1\otimes\omega_X$, in Part iii) of the above proposition. The
degree of $B_2$ is $b_2:=nd_1-d_2+(n-1)(2g-2)$.  In Remark
\ref{rem:StabII}, we have already stressed that $b_2$ can take on any
prescribed value. The understanding of our moduli spaces thus rests
on our understanding of the spaces of global sections of the
semistable vector bundle $B_2$. This gives an interesting link
between our moduli problem and Brill-Noether theory that we shall
exploit below.
\end{remark}
In the next step, we will demonstrate a partial converse to Proposition \ref{bound103}.
\begin{proposition}
\label{prop:StabV}
Assume $\alpha=(0,\alpha_1,\alpha_2)$ is a stability parameter, such
that $(\alpha_1,\alpha_2)$ belongs to the chamber ${\mathcal
C}_\infty$. Suppose that $C:=(E_0,E_1,E_2;\phi_1,\phi_2)$ is a
holomorphic 3-chain of type $t$, such that $E_1$ is a line bundle of
degree $d_1$, $Q_0$ a {\bfseries stable} vector bundle of degree
$d_0-d_1$ and rank $m-1$,
$$
\begin{CD}
0 @>>> E_1 @>\phi_1 >> E_0 @>>> Q_0 @>>> 0
\end{CD}
$$
a non-split extension, $K_2$ a {\bfseries stable} vector bundle of rank $n-1$ and degree $d_2-d_1$, and 
$$
\begin{CD}
0 @>>> K_2 @>>> E_2 @>\phi_2>> E_1 @>>> 0
\end{CD}
$$
a non-split extension. Then, $C$ is $\alpha$-stable.
\end{proposition}
\begin{proof}
First, we look at a subchain $(0,0,F_2)$ where $F_2$ is a subbundle of $K_2$. The definition of the wall $L_{\rm IV}$
implies that, for a stability parameter $\widetilde{\alpha}=(0,\widetilde{\alpha}_1,\widetilde{\alpha}_2)$,
such that $(\widetilde{\alpha}_1,\widetilde{\alpha}_2)$ lies in the interior of the parameter region $R(t)$,
one has 
$$
\frac{d_2-d_1}{n-1}+\widetilde{\alpha}_2<\mu_{\widetilde{\alpha}}(C).
$$
Since $K_2$ is, by assumption, a stable vector bundle, we also have 
$$
\mu(F_2)+\widetilde{\alpha}_2\le \frac{d_2-d_1}{n-1}+\widetilde{\alpha}_2< \mu_{\widetilde{\alpha}}(C),
$$
so that the chain $(0,0,F_2)$ is not destabilizing for any stability parameter as in the proposition.
\par
With the methods of the proof of Proposition \ref{StabIII}, one also checks that no subchain of the form $(F_0,0,0)$
is destabilizing. This also implies that subchains of the form $(F_0,0,F_2)$ aren't destabilizing either.
\par 
The remaining case to study is that of a subchain of the shape $(F_0,E_1,F_2)$. We begin with the following 
construction. We have the subchain $C':=(0,0,K_2;0,0)$ of $C$ and form the quotient chain $\widetilde{C}:=C/C'$. 
Note that $\widetilde{C}=(E_0,E_1,E_1; {\phi}_1,\widetilde{\phi}_2)$. This 3-chain possesses the subchain 
$C'':=(E_1\cong \im(\phi_1), E_1, E_1; {\phi}_1,\widetilde{\phi}_2)$. Let $C''':=\widetilde{C}/C''$.
The chain $C'''$ is given as $(Q_0,0,0;0,0)$. Let $\alpha^M:=(0,\alpha_1^M,\alpha_2^M)$ be the stability parameter that is characterized by the condition that $(\alpha_1^M,\alpha_2^M)$ is the point of intersection of the
lines $L_{\rm I}$ and $L_{\rm IV}$. By construction, we have
\begin{equation}
\label{eq:Gleich100}
\mu_{\alpha^M}(C')=\mu_{\alpha^M}(C'')=\mu_{\alpha^M}(C''')=\mu_{\alpha^M}(t).
\end{equation}
The chains $C'$ and $C'''$ are $\alpha$- and $\alpha^M$-stable, because the vector bundles $Q_0$ and $K_2$ are 
stable. We claim that the chain $C''$ is $\alpha^M$-stable, too. For this, we have to check the two subchains 
$C_1:=(E_1,E_1,0)$ and $C_2:=(E_1,0,0)$. The condition
$$
\mu_{\alpha^M}(C_i)<\mu_{\alpha^M}(t)=\mu_{\alpha^M}(C''),\q i=1,2,
$$
follows by applying a trick similar to the one used in the proof of Proposition \ref{StabIII}, because $C_1$ and $C_2$ may also be viewed as subchains of $C$, 
and the fact that 
$(\alpha_1^M,\alpha_2^M)$ lies in the interior of the intersection of the two half spaces defined by the Inequalities II and III
(see Figure 1) implies that the two subchains $(E_0,0,0)$ and $(E_0,E_1,0)$ do not destabilize $C$.
\par
Now, we return to a subchain $\ol{C}=(F_0,E_1,F_2)$ of $C$. We write
the vector bundle $F_2$ as the extension
$$
\begin{CD}
0 @>>> \widetilde{K}_2:=F_2\cap K_2 @>>> F_2 @>>> Q_2 @>>> 0,
\end{CD}
$$
and the vector bundle $F_0$ as the extension
$$
\begin{CD}
0 @>>> E_1 @>>> F_0 @>>> \widetilde{Q}_0 @>>> 0.
\end{CD}
$$
With these constructions, we find the subchain $\ol{C}':=(0,0,\widetilde{K}_2)$ of $C'$, the subchain $\ol{C}'':=
(E_1,E_1,Q_2)$ of $C''$, and the subchain $\ol{C}''':=(\widetilde{Q}_0,0,0)$ of $C'''$. 
\par
Exploiting the $\alpha^{M}$-stability of $C'$, $C''$, and $C'''$ and the fact that one of the subchains $\ol C'$, $\ol C''$, and $\ol C'''$ will be a proper one, we find
\begin{eqnarray*}
&& \mu_{\alpha^M}(\ol C)
\\
&=& \frac{1}{\rk(\ol{C})}\Bigl(\rk(\ol{C}')\mu_{\alpha^M}(\ol{C}')+\rk(\ol{C}'')
\mu_{\alpha^M}(\ol{C}'')+\rk(\ol{C}''') \mu_{\alpha^M}(\ol{C}''')\Bigr)
\\
&<& \frac{1}{\rk(\ol{C})}\Bigl(\rk(\ol{C}')\mu_{\alpha^M}(C')+\rk(C'')
\mu_{\alpha^M}(C'')+\rk(C''') \mu_{\alpha^M}(C''')\Bigr)
\\
&=&\mu_{\alpha^M}(t).
\end{eqnarray*}
From this inequality, it is clear that $\mu_\alpha(\ol{C})<\mu_\alpha(t)$ holds for all $\alpha$, such that
$(\alpha_1,\alpha_2)$ is sufficiently close to $(\alpha_1^M,\alpha_2^M)$.
By the definition of the walls, the same must be true for all $\alpha$, such that $(\alpha_1,\alpha_2)$
lies in the chamber ${\mathcal C}_\infty$.
\end{proof}
Next, we will include some observations regarding the relationship with Brill-Noether theory. 
(A survey on Brill-Noether theory which was also very helpful to the authors is \cite{Mercat}.)
Let $\alpha=(0,\alpha_1,\alpha_2)$ be a stability parameter, such that $(\alpha_1,\alpha_2)\in {\mathcal C}_\infty$.
Note that a holomorphic 3-chain $C$ of type $t$ is $\alpha$-stable, if and only if it is $\alpha$-semistable.
By Theorem \ref{thm:smoothdim}, the moduli space ${\mathcal M}_\alpha(t)$ is smooth, and, if non-empty, it has
dimension 
$$
\begin{array}{c}
\underbrace{(m-1)^2\bigl(g-1\bigr)+(d_0-md_1)+(m-1)\bigl(g-1\bigr)+g}_{=:f}
\\
+
\\
\underbrace{\Bigl((n-1)^2\bigl(g-1\bigr)+1\Bigr)}_{=:h}+
\underbrace{\Bigl(nd_1-d_2+(n-1)\bigl(g-1\bigr)-1\Bigr)}_{=:c}.
\end{array}
$$
Recall the notation from Remark \ref{rem:VarExt}. In that notation, we define the morphism
\begin{eqnarray*}
\Phi\colon {\mathcal M}_\alpha(t)&\lra& {\mathcal U}(n-1,b_2)
\\
C=(E_0,E_1,E_2;\phi_1,\phi_2) &\lma& B_2.
\end{eqnarray*}
Call a vector bundle $E$ on $X$ \it special\rm, if both $h^0(E)\neq 0$ and $h^1(E)\neq 0$. Using the map
$\Phi$ and the dimension formula for ${\mathcal M}_\alpha(t)$, we find the following result.
\begin{theorem}[Laumon \cite{Laumon}]
Let $r>0$ and $l$ be integers. Then, the generic vector bundle $E$ of rank $r$ and degree $l$ on $X$ is non-special. 
\end{theorem}
\begin{proof}
By Serre duality, we may assume that $l\le r(g-1)$. For these values, a vector bundle $E$ of rank $r$ and degree
$l$ is special, if and only if $h^0(E)>0$. Now, we may pick $m$, $n$, $d_0$, $d_1$, $d_2$ in such a way
that $r=n-1$, $l=b_2=b_2(g,n,d_1,d_2)$ (see Remark \ref{rem:VarExt}), and $d_0>md_1$. Our considerations, in particular
Proposition \ref{prop:StabV}, imply that the intersection of the image of $\Phi$ with ${\mathcal U}^{\rm s}(r,l)$
consists exactly of the set of isomorphism classes of special stable vector bundles of rank $r=n-1$ and degree $l=b_2$. One easily verifies that the dimension of 
any fiber of $\Phi$ over a special stable vector bundle is at least $f$. Therefore, the image has dimension at most 
$h+c$. Now, $h$ is the dimension of ${\mathcal U}^{\rm s}(r,l)$
and $c$ is negative. This proves that the generic stable vector bundle is non-special. Since the generic vector bundle
is stable, we are done.
\end{proof}
\begin{remark}
The above result is Corollary 1.7 in \cite{Laumon}. Of course, our
proof is neither easier nor more natural than the one of Laumon, but
it is a nice illustration of the strength of our results. Moreover,
the morphism $\Phi$ relates the geometry of the moduli space
${\mathcal M}_\alpha(t)$ to the geometry of the Brill-Noether locus
inside ${\mathcal U}(n-1,b_2)$. This enables us to derive fundamental
properties of our moduli spaces from the basic results in
Brill-Noether theory.
\end{remark}
Following Laumon, we let $W^l_{X,r}$ be the closed substack of vector
bundles that do possess global sections inside the stack ${\rm
Bun}^l_{X,r}$ of all vector bundles of rank $r$ and degree $l$ on
$X$.
\begin{theorem}
\label{thm:stack1}
{\rm i)} The stack $W^l_{X,r}$ is irreducible.
\par
{\rm ii)} For every $r$ and every $l> 0$, there exist stable vector
bundles of rank $r$ and degree $l$ with global sections.
\end{theorem}
\begin{proof}
i) This is \cite{Laumon}, Corollary 5.2.
\par
ii) For $l>n(g-1)$, the existence of global sections follows from
the theorem of Riemann-Roch. In the remaining range $0<l\le
n(g-1)$, Sundaram \cite{Sundaram} shows that there exist stable
vector bundles of rank $r$ and degree $l$ with global sections. (Of
course, for $l=0$, there are no stable vector bundles with global
sections.)
\end{proof}
\begin{remark}
If $l\le r(g-1)$, then the generic vector bundle in $W^l_{X,r}$
has precisely one global section.  This follows from \cite{Laumon},
Lemma 2.6.
\end{remark}
\begin{proposition}
\label{prop:DescMod2}
Assume that the type $t=(m,1,n;d_0,d_1,d_2)$ is such that $b_2>0$
(see {\rm Remark \ref{rem:VarExt}}) and $d_0>md_1$ and that
$\alpha=(0,\alpha_1,\alpha_2)$ is a stability parameter with
$(\alpha_1,\alpha_2)\in {\mathcal C}_\infty$.  Then, the moduli space
${\mathcal M}_\alpha(t)$ is a smooth connected projective variety of
dimension
$$
\bigl(g-1\bigr)\bigl(m^2+1+n^2-m-n\bigr)+(d_0-md_1)+(nd_1-d_2)+1.
$$ 
\end{proposition}
\begin{proof}
Proposition \ref{prop:StabV} and Theorem \ref{thm:stack1} grant that
${\mathcal M}_\alpha(t)$ is non-empty.  By Theorem
\ref{thm:smoothdim}, the moduli space is smooth of the indicated
dimension. It remains to check that it is also connected.  Since the
semistable vector bundles form an open substack of the stack of all
vector bundles, the stack of semistable vector bundles of rank $r$
and degree $l$ with global sections is still irreducible, by Theorem
\ref{thm:stack1}, i).  Using this and the techniques from the proof
of Corollary \ref{cor:DescMod1}, one easily constructs a connected
parameter space for $\alpha$-stable chains of type $t$.
\end{proof}
\begin{remark}
i) We remind the reader that the conditions $b_2\ge 0$ and $d_0>md_1$
are both necessary for the non-emptiness of the moduli spaces. The
first one, because a semistable vector bundle of negative degree
never possesses global sections, and the second one, because
otherwise the interior of the parameter region is empty, by Remark
\ref{StabI}.  Observe that the existence problem for $b_2=0$ has been
left open.
\par
ii) Let $l>0$ and $d$ be integers and ${\mathcal W}^l_{X,r}\subset
{\mathcal U}^{\rm s}(r,l)$ the locus of stable vector bundles with
global sections.  Under the assumption that $0<l\le r(g-1)$,
Hoffmann \cite[Example 5.11]{Hoff} has recently shown that ${\mathcal
W}^l_{X,r}$ is birationally equivalent to ${\mathbb P}^s\times
{\mathcal J}^0(X)$. Thus, the morphism $\Phi$ introduced before will
help in a more detailed investigation of the geometry of the moduli
spaces of holomorphic chains.
\end{remark}
To conclude this example, we will again determine a region of
parameters, such that all moduli spaces associated to parameters in
that region will be birationally equivalent.
\begin{proposition}
\label{prop:BirRegion2}
Let $\alpha^i=(0,\alpha^i_1,\alpha^i_2)$, $i=1,2$, be two stability
parameters, such that $(\alpha^i_1,\alpha^i_2)\in R(t)\cap R_{2g-2}$,
$i=1,2$, that do not lie on any wall.  Then, ${\mathcal
M}_{\alpha^1}(t)$ and ${\mathcal M}_{\alpha^2}(t)$ are birationally
equivalent smooth projective varieties which either have dimension
$$
\bigl(g-1\bigr)\bigl(m^2+1+n^2-m-n\bigr)+(d_0-md_1)+(nd_1-d_2)+1
$$ 
or are empty.
\end{proposition}
\begin{proof}
We use the same discussion as at the beginning of the proof of
Proposition \ref{prop:BirRegion1}.  Note that we may always replace
$C'$ by an $\alpha^0$-stable subchain, i.e., we can directly require
$C'$ to be $\alpha^0$-stable.
\par
In the respective notation we distinguish the two cases a)
$(r_0',r_1',r_2')=(m',1,n')$ with $m'n'>0$ and
$(r_0'',r_1'',r_2'')=(m'',0,n'')$ with either $m''\neq 0$ or $n''\neq
0$ and b) $(r_0',r_1',r_2')=(m',0,n')$ with either $m'\neq 0$ or
$n'\neq 0$ and $(r_0',r_1',r_2')=(m'',1,n'')$ with $m''n''>0$.
\par
In Case a), the $\C$-linear chain $\delta_{[0,2]}$ will appear in the
decomposition of $V'$

into indecomposable chains. We easily see
$\Hom(\delta_{[0,0]},\allowbreak \delta_{[0,2]})\cong \C$.  Thus,
$\delta_{[0,0]}$ cannot occur in the decomposition of $V''$, so that
$V''=\delta_{[2,2]}^{\oplus n''}$.  The formula for $\chi(V'',V')$
yields the value $n'n''-n''$. We want this to be zero, so that we
conclude $n'=1$ and $n''=n-1$. The rank type of $V'$ is thus
$(m,1,1)$ and the rank type of $V''$ is $(0,0,n-1)$. This is again a
possibility which is allowed by Theorem \ref{thm:BirRegionFund}.
\par
Let $C'=(E_0,E_1,E_2';\phi_1,\phi_2')$ be an $\alpha^0$-stable chain
of type $(d_0,d_1,d_2';\allowbreak m,1,1)$ and $C''=(0,0,E_2'';0,0)$
an $\alpha^0$-semistable holomorphic chain of type
$(0,0,d_2'';0,0,n-1)$.  We have to evaluate the dimension of the
locus of chains $C$ which might be obtained as non-split extensions
of a chain such as $C''$ by a chain such as $C'$. By the dimension
formula Theorem \ref{thm:smoothdim}, iv), the moduli space for
$\alpha^0$-stable chains of type $(d_0,d_1,d_2'; m,1,1)$ has
dimension $d_0-(m-1)d_1-d_2+(m-1)m(g-1)+g$.  The moduli space
of chains of type $(0,0,d_2'';0,0,n-1)$ agrees with the moduli space
of semistable vector bundles of rank $n-1$ and degree $d_2''$ and has
dimension $(n-1)^2(g-1)$. Again, we find $\Hom(C'',C')=\{0\}
={\mathbb H}^2(C'',C')$ (Proposition \ref{prop:h2-vanishing}, iv), so
that Proposition \ref{prop:chi(C'',C')} gives $\dim_\C({\rm
Ext}^1(C'',C'))=(n-1)(d_1-d_2')$. All in all, the dimension of the
locus we wish to describe is
\begin{equation}
\label{eq:DimensionCount}
\bigl(g-1\bigr)(m^2+1+n^2-m-2n+1)+(d_0-md_1)+n(d_1-d_2'')+1.
\end{equation}
We have $\mu_{\alpha^0}(C')=\mu_{\alpha^0}(C)$ which gives
\begin{equation}
\label{eq:WallDefine3}
-(n-1)\alpha_1^0+(m+1)(n-1)\alpha_2^0=(n-1)(d_0+d_1)-(m+2)d_2+(m+n+1)d_2'.
\end{equation}
Since $\alpha^0$ is supposed to lie in the interior of the parameter region $R(t)$, we must
have
$$
(n-1)(d_0+d_1)-(m+2)d_2+(m+n+1)d_2'>\frac{n-1}{n}\bigl(nd_0+nd_1-(m+1)d_2\bigr)
$$
that is
$$
nd_2'>d_2.
$$
Together with Formula (\ref{eq:DimensionCount}), this shows that
the codimension of the flip locus in question is at least
$(n-1)(g-1)+1$. It, therefore, may be neglected.
\par
Case b) follows immediately from Case a) by passing to the \sl dual
chains \rm (see Remark \ref{quotientformulation}, iv).  The remaining
flip loci are covered by Proposition \ref{codim}.
\end{proof}
\begin{remark}
\label{rem:WallDefine4}
In (\ref{eq:WallDefine3}), the maximal value of $d_2'$ is $d_1$.
Plugging $d_1$ into that formula gives the equation for the line
$L_{\rm IV}$. Likewise, one obtains the equation for the line $L_{\rm
I}$ (this is hidden in the argument with the dual chains).
\end{remark}
\subsection{Moduli spaces for type $(1, m, 1;d_0,d_1,d_2)$}
\label{sub:large-moduli-1-m-1}
In the examples above, we have merely used vector bundle techniques
in order to analyze the asymptotic moduli spaces. Here, we will
discuss an elementary example where we make use of extensions \sl in
the category of holomorphic chains \rm in order to gather interesting
information on the geometry of the asymptotic moduli spaces.
\par
Let us recall the relevant inequalities from Proposition
\ref{prop:RankMaxChains} that bound the parameter region $R(t)$.
\par
\it Inequality \rm I. 
Using the test object $(E_0,0,0)$, one arrives at the inequality
$$
m\alpha_1+\alpha_2\ge A_{\rm I}:=(m+1)d_0-d_1-d_2.
$$
\par
\it Inequality \rm II.
Testing $\alpha$-semistability with the subchain $(E_0,E_1,0)$ yields
the inequality
$$
-m\alpha_1+(m+1)\alpha_2\ge A_{\rm II}:=d_0+d_1-(m+1)d_2.
$$
\par
\it Inequality \rm III.
The subchain $(E_0, E_2\cong \im(\phi_2), E_2)$ provides us with the
inequality
$$
-2\alpha_1+\alpha_2\le A_{\rm III}:= -d_0+\frac{3}{m-1}d_1-\frac{2m+1}{m-1}d_2.
$$
\par
\it Inequality \rm IV.
Finally, one finds the inequality 
$$
-2\alpha_1+\alpha_2\ge A_{\rm IV}:=-\frac{2m+1}{m-1}d_0+\frac{3}{m-1}d_1-d_2
$$
with the subchain $(0,\ker(\phi_1),0)$.
\begin{remark}
\label{rem:Parallelo2}
i) If we wish that the interior of the parameter region $R(t)$
becomes non-empty, then we must choose the numerical data in such a
way that $A_{\rm III}>A_{\rm IV}$ holds true. This condition simply
amounts to $d_0>d_2$.
\par
ii) Inequality III and IV bound a strip in the
$(\alpha_1,\alpha_2)$-plane. The remaining inequalities provide the
lower bounds for the parameter region $R(t)$.  Figure 7 illustrates
the shape of the resulting domain.
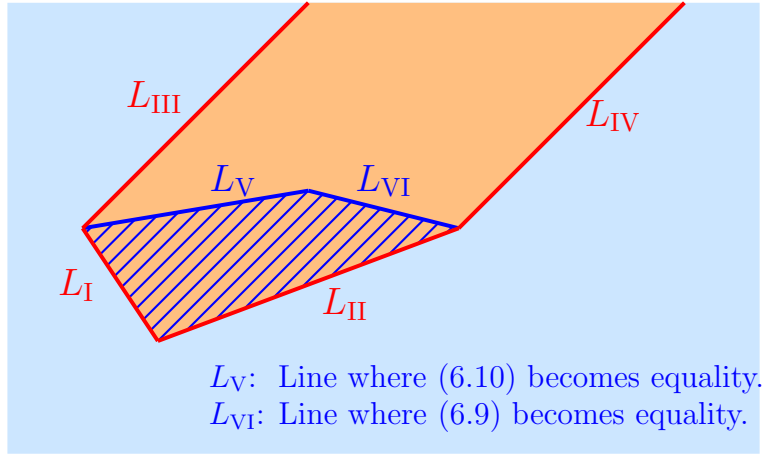
\begin{figure}[h]
\begin{center}
\begin{pspicture}(5,1)(4.5,6.5)
\psset{xunit=2cm, yunit=1cm}
\pspolygon*[linecolor=lblue](0,.5)(0,6.5)(5,6.5)(5,.5)
\pspolygon*[linecolor=lightorange](2,6.5)(.5,3.5)(1,2)(3,3.5)(4.5,6.5)
\pspolygon[linecolor=lblue, fillstyle=hlines, hatchcolor=blue](.5,3.5)(1,2)(3,3.5)(2,4)
\psline[linecolor=blue, linewidth=1.5pt](3,3.5)(2,4)
\psline[linecolor=blue, linewidth=1.5pt](.5,3.5)(2,4)
\uput[90](1.5,3.75){\large\blue $L_{\rm V}$}
\uput[90](2.5,3.75){\large\blue $L_{\rm VI}$}
\psline[linecolor=red, linewidth=1.5pt](0.5,3.5)(2,6.5)
\uput[180](1.25,5.25){\large\red $L_{\rm III}$}
\psline[linecolor=red, linewidth=1.5pt](3,3.5)(4.5,6.5)
\uput[0](3.75,5){\large\red $L_{\rm IV}$}
\psline[linecolor=red, linewidth=1.5pt](3,3.5)(1,2)
\uput[0](2,2.5){\large\red $L_{\rm II}$}
\psline[linecolor=red, linewidth=1.5pt](.5,3.5)(1,2)
\uput[0](.25,2.75){\large\red $L_{\rm I}$}
\uput[0](1.25,1.5){\blue $L_{\rm V}$: \hskip .1cm Line where (\ref{eq:RankMaximal2}) becomes equality.}
\uput[0](1.25,1){\blue $L_{\rm VI}$: Line where (\ref{eq:RankMaximal}) becomes equality.}
\end{pspicture}
\end{center}
\caption{The parameter region for type $(r_0,r_1,r_2;d_0,d_1,d_2)$.}
\end{figure}
\par
iii) Let us examine when the map $\phi_1\circ\phi_2\colon E_2\lra
E_0$ is non-zero. If $\phi_1\circ\phi_2\equiv 0$, then we obtain the
subchain $(0,\ker(\phi_1),E_2)$. Noting that $\deg(\ker(\phi_1))\ge
d_1-d_0$, the condition of $\alpha$-semistability applied to the
given subchain yields the inequality
\begin{equation}
\label{eq:RankMaximal}
(m-2)\alpha_1+2\alpha_2\le 2(m+1)d_0-2d_1-2d_2.
\end{equation}
Alternatively, one may use the subchain $(0,E_2\cong\phi_2(E_2),
E_2)$. This leads to the inequality
\begin{equation}
\label{eq:RankMaximal2}
(-m+2)\alpha_1+m\alpha_2\le 2d_0+2d_1-2(m+1)d_2.
\end{equation}
The reader may check that Inequality I-IV together with either
Inequality (\ref{eq:RankMaximal}) or Inequality
(\ref{eq:RankMaximal2}) cuts out a bounded region in the
$(\alpha_1,\alpha_2)$-plane. Therefore, if the converse to either
(\ref{eq:RankMaximal}) or (\ref{eq:RankMaximal2}) holds, then
$\phi_1\circ \phi_2$ is non-trivial (whence a generic isomorphism).
\end{remark}
In Figure 4, we have already sketched the parameter region $R(t)$
away from some bounded subregion and selected a chamber ${\mathcal
C}_\infty$. It is formally defined to be the (unique) unbounded
two-dimensional chamber whose closure intersects the line $L_{\rm
III}$ where Inequality III becomes an equality.  This definition
involves that the converse to both (\ref{eq:RankMaximal}) and
(\ref{eq:RankMaximal2}) is verified for the elements of ${\mathcal
C}_\infty$.
\begin{remark}
In this example, one might also declare the unbounded two-dimensional
chamber whose closure intersects the line $L_{\rm IV}$ where
Inequality IV becomes equality to be the chamber ${\mathcal
C}_\infty$. The reader may verify that techniques analogous to those
presented in the following lead to a description of the relevant
moduli spaces.
\end{remark}
\begin{proposition}
\label{prop:Asympt1m1}
Let $\alpha=(0,\alpha_1,\alpha_2)$ be a stability parameter, such
that $(\alpha_1,\alpha_2)\in {\mathcal C}_\infty$, and let
$C:=(E_0,E_1,E_2;\phi_1,\phi_2)$ be an $\alpha$-semistable
holomorphic chain of type $(1,m,1;d_0,d_1,d_2)$. Then, $C$ is a
non-split extension of the holomorphic chain
$$
\begin{CD}
C'':= C/C': 0 @>>> Q_1:=E_1/\phi_2(E_2)@>>> 0
\end{CD} 
$$
by the chain
$$
\begin{CD}
C':E_2 @>\phi_2>> E_2\cong \im(\phi_2) @>\cdot D>> E_2(D)\cong E_0.
\end{CD}
$$
Here, $D$ is an effective divisor of degree $d_0-d_2$, and $Q_1$
is a semistable vector bundle of rank $m-1$ and degree $d_1-d_2$.
\par
Furthermore,
$$
\dim_\C\bigl({\rm Ext}^1(C'',C')\bigr)=(m-1)(d_0-d_2).
$$
\end{proposition}
\begin{proof}
Everything apart from the formula for $\dim_\C({\rm Ext}^1(C'',C'))$
follows from the definition of the line $L_{\rm III}$ and the
adjacency of ${\mathcal C}_\infty$ to that line.
\par
First, note that Proposition \ref{prop:h2-vanishing}, iv), grants
that ${\mathbb H}^2(C'',C')=\{0\}$. Next, we obviously have ${\rm
Hom}(C'',C')=\{0\}$.  Therefore, we find $\dim_\C({\rm
Ext}^1(C'',C'))=-\chi(C'',C')$. Finally, we use Proposition
\ref{prop:chi(C'',C')} to compute $\chi(C'',C')$.
\end{proof}
Again, we prove a partial converse to Proposition \ref{prop:Asympt1m1}.
\begin{proposition}
\label{prop:2Asympt1m1}
Assume that $\alpha=(0,\alpha_1,\alpha_2)$ is a stability parameter
with $(\alpha_1,\alpha_2)\in {\mathcal C}_\infty$. Let $E_2$ be a
line bundle of degree $d_2$, $D$ an effective divisor of degree
$d_0-d_2$, and $Q_1$ a {\bfseries stable} vector bundle of rank $m-1$
and degree $d_1-d_2$. Set $E_0:=E_2(D)$ and define the chain
$$
\begin{CD}
C'': 0 @>>> Q_1@>>> 0
\end{CD} 
$$
as well as the chain
$$
\begin{CD}
C':E_2 @>{\id}_{E_2}>> E_2@>\cdot D>>  E_0.
\end{CD}
$$
Any non-split extension $C$ of the chain $C''$ by the chain $C'$
is an $\alpha$-stable holomorphic 3-chain of type
$t=(1,m,1;d_0,d_1,d_2)$.
\end{proposition}
\begin{proof}
We fix a stability parameter $\alpha^M=(0,\alpha_1^M,\alpha_2^M)$,
such that $(\alpha_1^M,\allowbreak\alpha^M_2)$ lies on the line
$L_{\rm III}$, in the interior of the region described by the
remaining Inequalities I, II, and IV, and in the closure of
${\mathcal C}_\infty$.  Let $C$ be a chain as in the statement of the
proposition. We start with the investigation of some special
subchains. First, we recall that the subchain $C_1:=(E_0,0,0)$ is
neither $\alpha$- nor $\alpha^M$-destabilizing, because
$(\alpha_1,\alpha_2)$ and $(\alpha_1^M,\alpha_2^M)$ both lie in the
interior of the region from Figure 1. Next, we look at the subchain
$C_2:=(E_0,E_2,E_2)$.  The definition of the line $L_{\rm III}$ and
the fact that $(\alpha_1,\alpha_2)$ lies below that line imply that
$C_2$ does not $\alpha$-destabilize $C$ either. (However, we have
$\mu_{\alpha^M}(C_2)=\mu_{\alpha^M}(C)$.)  As the third subchain, we
define $C_3:=(E_0,E_2,0)$ and we ask whether the inequality
$$
\frac{d_0+d_2+\alpha_1}{2}<\frac{d_0+d_1+d_2+m\alpha_1+\alpha_2}{m+2}
$$
is verified, that is, if 
$$
(m-2)\alpha_1+2\alpha_2> md_0-2d_1+md_2
$$
holds true. As we have pointed out after Remark
\ref{rem:Parallelo2}, the definition of the chamber ${\mathcal
C}_\infty$ involves the inequality
$$
(m-2)\alpha_1+2\alpha_2> 2(m+1)d_0-2d_1-2d_2.
$$
One checks that the necessary condition $d_0>d_2$ implies
$$
2(m+1)d_0-2d_1-2d_2>md_0-2d_1+md_2,
$$
so that we may conclude that $C_3$ isn't an $\alpha$-destabilizing
subchain for $C$. The same conclusion applies with respect to
$\alpha^M$.
\par
Now, let $\ol C=(F_0,F_1,F_2)$ be any subchain of $C$. Then, we write
the vector bundle $F_1$ as the extension
$$
\begin{CD}
0 @>>> E_2\cap F_1@>>> F_1@>>> \widetilde{Q}_1 @>>> 0.
\end{CD}
$$
First, we suppose that $F_0$ and $F_2$ are both trivial. Then,
$E_2\cap F_1=\{0\}$, so that we may view $F_1$ as a subsheaf of
$Q_1$. It is a proper subsheaf. Otherwise, we would have
$E_1=E_2\oplus Q_1$. But, since $F_1$ obviously agrees with
$\ker(\phi_1)$, this implies $C=C'\oplus C''$, contradicting our
assumption.  Therefore, since $Q_1$ is stable,
$$
\mu_{\alpha^M}(\ol C)=\mu(F_1)+\alpha_1^M <
\mu(Q_1)+\alpha_1^M=\mu_{\alpha^M}(C')=\mu_{\alpha^M}(C).
$$
In the remaining cases, we may write the chain $\ol C$ as an
extension of the chain $(0,\widetilde{Q}_1,0;0,0)$ by one of the
subchains $C_1$, $C_2$, or $C_3$. Since $Q_1$ is a stable bundle, we
find
$$
\mu(\widetilde{Q}_1)+\alpha_1^M\le \mu(Q_1)+\alpha_1^M=\mu_{\alpha^M}(C').
$$
Using the results on the subchains $C_1$, $C_2$, and $C_3$ (and
the fact that $\widetilde{Q}_1$ will be a proper subbundle of $Q_1$
in the case of the chain $C_2$), one checks that
$$
\mu_{\alpha^M}(\ol C)< \mu_{\alpha^M}(C)
$$
is verified. As in the proof of Proposition \ref{StabIII}, we
easily derive the assertion of the proposition.
\end{proof}
Our discussions yield the following information on the moduli spaces
which belong to a stability parameter in the chamber ${\mathcal
C}_\infty$.
\begin{corollary}
\label{cor:DescMod2}
The moduli space ${\mathcal M}_\alpha(t)$ for
$\alpha=(0,\alpha_1,\alpha_2)$ and $(\alpha_1,\allowbreak\alpha_2)\in
{\mathcal C}_\infty$ is a connected smooth projective variety of
dimension
$$
(m-1)^2\bigl(g-1\bigr)+g+m(d_0-d_2).
$$
It is birationally equivalent to a ${\mathbb P}^N$-bundle over the
product ${\mathcal J}^{d_2}\times X^{(d_0-d_2)}\times {\mathcal
U}^{\rm s}(m-1,d_1-d_2)$ of the Jacobian of degree $d_0$ line
bundles, the $(d_0-d_2)$-fold symmetric product of the curve, and the
moduli space of stable vector bundles of rank $(m-1)$ and degree
$(d_1-d_2)$, $N:=(m-1)(d_0-d_2)-1$.
\end{corollary}
\begin{proof}
The smoothness follows again from the fact that all
$\alpha$-semistable 3-chains of type $t$ are $\alpha$-stable and
Theorem \ref{thm:smoothdim}, v).  Proposition \ref{prop:Asympt1m1}
and \ref{prop:2Asympt1m1} establish the assertions about the
dimension and the birational geometry. The irreducibility of the
moduli space may be proved along the same lines as Corollary
\ref{cor:DescMod1}.
\end{proof}
\begin{remark}
The reader may try to find the result by just using extension
techniques for vector bundles. That approach doesn't seem to work
properly.
\end{remark}
\begin{proposition}
\label{prop:BirRegion3}
Let $\alpha^i=(0,\alpha^i_1,\alpha^i_2)$, $i=1,2$, be two stability
parameters, such that $(\alpha^i_1,\alpha^i_2)\in R(t)\cap R_{2g-2}$,
$i=1,2$, that do not lie on any wall.  Then, ${\mathcal
M}_{\alpha^1}(t)$ and ${\mathcal M}_{\alpha^2}(t)$ are birationally
equivalent smooth projective varieties which either have dimension
$$
(m-1)^2\bigl(g-1\bigr)+g+m(d_0-d_2).
$$
or are empty.
\end{proposition}
\begin{proof}
We use the set-up described at the beginning of the proof of Proposition \ref{prop:BirRegion1}.
This time, there are the cases a) $(r_0',r_1',r_2')=(0,m',1)$ and $(r_0'',r_1'',r_2'')=(1,m'',0)$
with $m'm''>0$, b) $(r_0',r_1',r_2')=(1,m',0)$ and $(r_0'',r_1'',r_2'')=(0,m'',1)$
with $m'm''>0$, c) $(r_0',r_1',r_2')=(1,m',1)$ and $(r_0'',r_1'',r_2'')=(0,m'',0)$
with $m'm''>0$, and d) $(r_0',r_1',r_2')=(0,m',0)$ and $(r_0'',r_1'',r_2'')=(1,m'',1)$
with $m'm''>0$.
\par
Recall that we would like to first determine the cases when the map
$b$ from (\ref{eq:extension-complex}) may be an isomorphism. This
requires $\chi(V'',V')=0$ and $\Hom(V'',\allowbreak V')=\{0\}$ (Lemma
\ref{lem:cohomology-linear-chains}).
In Case a), we have $\chi(V'',V')=m'm''$. This is never zero, so that
this case needs not be considered.
In Case b), we have $\chi(V'',V')=m'm''-m$. This is non-zero except
for the case $m=4$ and $m'=2=m''$.
In that case, the decomposition of both $V'$ and $V''$ contains the
linear chain $\delta_{[1,1]}$, whence $\Hom(V'',V')\neq \{0\}$, and
that excludes this case, too.
\par
In Case c), we have $\chi(V'',V')=0$ and $\Hom(V'',V')=\{0\}$, if and
only if $m'=1$ and $m''=m-1$.  The $\alpha^0$-semistable chains $C'$
of type $(1,1,1; d_0,d_1',d_2)$ have a $(g+d_0-d_2)$-dimensional
moduli space whereas the $\alpha^0$-semistable chains of type
$(0,m-1,0;d_0,d_1'',d_2)$ have an $((m-1)^2(g-1)+1)$-dimensional
moduli space. For these chains, we compute that $\Hom(C'',C')=\{0\}$,
${\mathbb H}^2(C'',C')=\{0\}$, by Proposition
\ref{prop:h2-vanishing}, and that $\dim_\C({\rm
Ext}^1(C'',C'))=(m-1)(d_0-d_1')$, by Proposition
\ref{prop:chi(C'',C')}. The dimension of the locus of chains $C$ that
may be written as an extension of a chain $C''$ as above by a chain
$C'$ as above thus has dimension
$$
(m-1)^2\bigl(g-1\bigr)+g+(d_0-d_2)+(m-1)(d_0-d_1').
$$
We obviously have $d_2\le d_1'\le d_0$. If $d_1'=d_2$, then
$\alpha^0$ has to lie on the line $L_{\rm III}$, and we exclude that.
Thus, we see that the flip locus under investigation has proper
codimension.
\par
Case d) reduces immediately to Case c), by passing to the dual chains
(Remark \ref{quotientformulation}, iv).  Together with Proposition
\ref{codim}, these considerations imply our contention.
\end{proof}
\begin{remark}
\label{rem:WallDefine5}
Again, we emphasize that the walls ${\mathcal B}(t',t'')$ from
Definition \ref{def:region-V} comprise the walls $L_{\rm III}$ and
$L_{\rm IV}$.
\end{remark}
\subsection{Concluding remarks regarding the parameter region and the
birationality region} 
In Definition \ref{def:region-V}, a certain region $\widetilde{R}(t)$
in the plane $\R^2$ has been defined, so that non-empty moduli spaces
belonging to parameters in the same connected component of the region
$R_{2g-2}\cap \widetilde{R}(t)$ are birationally equivalent. In the
proofs of Proposition \ref{prop:BirRegion1}, \ref{prop:BirRegion2},
and \ref{prop:BirRegion3}, we have seen that the definition of
$\widetilde{R}(t)$ is not optimal. Indeed, there are many walls in
${\mathcal B}(t',t'')$, such that the moduli spaces in two open
chambers adjacent to such a wall are still birationally equivalent.
It would be interesting to know, if the relevant computations may be
performed in greater generality in order to arrive at better results 
(cf. Remark~\ref{rem:discussion:parameter-region-vs-birational-region}).
\par
Another interesting feature, at which we pointed in Remark
\ref{rem:WallDefine2}, \ref{rem:WallDefine4}, and
\ref{rem:WallDefine5}, is that the walls that bound the parameter
region belong to those defined in Proposition \ref{StandardBounds}
and, more generally, to those of the form ${\mathcal B}(t',t'')$.  This
observation might help to find the a priori region for parameters
with non-empty moduli spaces in more general situations, i.e., for
chains of greater length or eventually more general quivers as
considered in~\cite{AG2, Schquiv}.

\end{document}